\def\DateTime{27/July/2004, 16:00(JP)}
\def\Version{Version $1.0a$}

\def\no{\if01}
\def\iftwelvept{\no}
\def\ifdedicatory{\no}
\def\ifusepdf{\no}
\def\ifpsfont{\no}

\iftwelvept
\documentclass[leqno,12pt]{amsart}
\else
\documentclass[leqno]{amsart}
\fi
\usepackage{amssymb}
\usepackage{amscd}
\usepackage{latexsym}
\usepackage{verbatim}
\usepackage[all]{xy}

\ifusepdf
\usepackage{hyperref}
\else\fi
\ifpsfont
\usepackage[T1]{fontenc}
\usepackage{times}
\else\fi

\iftwelvept
\else\fi

\theoremstyle{plain}
\newtheorem{Theorem}{Theorem}[section]
\newtheorem{Proposition}[Theorem]{Proposition}
\newtheorem{Lemma}[Theorem]{Lemma}
\newtheorem{Corollary}[Theorem]{Corollary}

\newtheorem{Claim}{Claim}[Theorem]

\theoremstyle{definition}

\newtheorem{Remark}[Theorem]{Remark}

\renewcommand{\theTheorem}{\arabic{section}.\arabic{Theorem}}
\renewcommand{\theClaim}{\arabic{section}.\arabic{Theorem}.\arabic{Claim}}
\renewcommand{\theequation}{\arabic{section}.\arabic{Theorem}.\arabic{Claim}}

\def\rom{\textup}
\newcommand{\ZZ}{{\mathbb{Z}}}
\newcommand{\QQ}{{\mathbb{Q}}}

\newcommand{\PP}{{\mathbb{P}}}
\newcommand{\NN}{{\mathbb{N}}}
\newcommand{\AAA}{{\mathbb{A}}}
\newcommand{\OO}{{\mathcal{O}}}

\newcommand{\Rat}{\operatorname{Rat}}
\newcommand{\SDRat}{\operatorname{SDRat}}
\newcommand{\Hom}{\operatorname{Hom}}
\newcommand{\Ext}{\operatorname{Ext}}
\newcommand{\Sing}{\operatorname{Sing}}

\newcommand{\Ker}{\operatorname{Ker}}
\newcommand{\Coker}{\operatorname{Coker}}
\newcommand{\Spec}{\operatorname{Spec}}
\newcommand{\Gal}{\operatorname{Gal}}
\newcommand{\Supp}{\operatorname{Supp}}
\newcommand{\codim}{\operatorname{codim}}
\newcommand{\mult}{\operatorname{mult}}
\newcommand{\gcm}{\operatorname{gcm}}
\newcommand{\trdeg}{\operatorname{tr.deg}}
\newcommand{\rank}{\operatorname{rk}}
\newcommand{\lformal}{[\![}
\newcommand{\rformal}{]\!]}

\newcommand{\Sym}{\operatorname{Sym}}

\newcommand{\Proof}{{\sl Proof.}\quad}
\newcommand{\QED}{{\unskip\nobreak\hfil\penalty50\quad\null\nobreak\hfil
{$\Box$}\parfillskip0pt\finalhyphendemerits0\par\medskip}}
\newcommand{\rest}[2]{\left.{#1}\right\vert_{{#2}}}

\begin{document}

\title
{Kobayashi-Ochiai's theorem for log schemes}
\author{Isamu Iwanari and Atsushi Moriwaki}
\address{Department of Mathematics, Faculty of Science,
Kyoto University, Kyoto, 606-8502, Japan}
\email{iwanari@math.kyoto-u.ac.jp}
\email{moriwaki@math.kyoto-u.ac.jp}
\date{\DateTime, (\Version)}
\ifdedicatory
\dedicatory{Dedicated to Professor Masaki Maruyama on his 60th birthday}
\else\fi
\begin{abstract}
Kobayashi-Ochiai's theorem says us that the set of dominant rational maps
to a complex variety of general type is finite.
In this paper, we give a generalization of it
in the category of log schemes.
\end{abstract}


\maketitle

\setcounter{tocdepth}{1}
\tableofcontents

\section*{Introduction}
\renewcommand{\theTheorem}{\Alph{Theorem}}

In the paper \cite{KoOchi},
Kobayashi and Ochiai proved that
the set of dominant rational maps to a complex
variety of general type is finite.
This result was generalized to the case over
a field of positive characteristic by
Dechamps and Menegaux \cite{DesMene}.
Furthermore, Tsushima \cite{Tsushima} established
finiteness for open varieties over a field of
characteristic zero.
In this paper, we consider their generalization in the
category of log schemes.
As we know, logarithmic geometry is a general framework
to cover compactification and singularities in degeneration.
The most typical example of these mixed phenomena
is a logarithmic structure on a semistable variety.
Actually, we deal with a log rational map on a semistable variety
with a logarithmic structure.
The following is the main theorem of this paper:

\begin{Theorem}
\label{thm:intro:main}
Let $k$ be an algebraically closed field and $M_k$
a fine log structure of $\Spec(k)$.
Let $X$ and $Y$ be proper semistable varieties
over $k$, and let $M_X$ and $M_Y$ be fine
log structures of $X$ and $Y$ over $M_k$
respectively such that
\[
 (X, M_X) \to (\Spec(k), M_k)
\quad\text{and}\quad
 (Y, M_Y) \to (\Spec(k), M_k)
\]
are log smooth and integral.
We assume that $(Y, M_Y)$ is of log general type, that is,
$\det(\Omega^1_{Y/k}(\log(M_Y/M_k)))$ is a big line bundle on $Y$
\rom{(}see Conventions and 
terminology \rom{\ref{subsub:big:line:bundle}}\rom{)}.
Then, the set of all log rational maps
\[
 (\phi, h) : (X, M_X) \dasharrow (Y, M_Y)
\]
over $(\Spec(k), M_k)$
with the following properties \rom{(1)} and \rom{(2)} is finite:
\begin{enumerate}
\renewcommand{\labelenumi}{(\arabic{enumi})}
\item
$\phi : X \dasharrow Y$ is a rational map
defined over a dense open set $U$ 
with $\codim(X \setminus U) \geq 2$,
and $(\phi, h) : (U, \rest{M_X}{U}) \to (Y, M_Y)$
is a log morphism over $(\Spec(k), M_k)$.

\item
For any irreducible component $X'$ of $X$,
there is an irreducible component $Y'$ of $Y$
such that $\phi(X') \subseteq Y'$ and
the induced rational map $\phi' : X' \dasharrow Y'$
is dominant and separable.
\end{enumerate}
\end{Theorem}

As an immediate corollary of the above theorem,
we have the following:

\begin{Corollary}
Let $X$ be a proper semistable variety over $k$ and
$M_X$ a fine log structure of $X$
over $M_k$
such that $(X, M_X) \to (\Spec(k), M_k)$
is log smooth and integral.
If $(X, M_X)$ is of log general type, then
the set of automorphisms of  $(X, M_X)$
over $(\Spec(k), M_k)$ is finite.
\end{Corollary}

Here let us give a sketch of the proof of Theorem~\ref{thm:intro:main}.
For this purpose,
we need to deal with the classical case and
the non-classical case.
In the case where $M_k = k^{\times}$ and $X$ and $Y$ are smooth over $k$
(the classical case),
we can use the similar arguments as in \cite{DesMene}. 
Actually, we prove it under the weaker conditions
(cf. Theorem~\ref{thm:finite:dominant:rat:map:general:type:smooth}).
However, if
$M_k$ is not trivial (the non-classical case), 
we have to determine a local description of a log structure.
Indeed, we have the following theorem:

\begin{Theorem}
Let $X$ be a semistable variety over $k$ and $M_X$
a fine log structure of $X$ over $M_k$ such that
$(X, M_X) \to (\Spec(k), M_k)$
is log smooth and integral.
Let us take a fine and sharp monoid $Q$ with $M_k = Q \times k^{\times}$.
For a closed point $x \in X$, 
there is a good chart
$(Q \to M_k,\ P \to M_{X, \bar{x}},\ Q \to P)$ of 
$(X, M_X) \to (\Spec(k), M_k)$ at $x$, namely,
\begin{enumerate}
\renewcommand{\labelenumi}{(\alph{enumi})}
\item
$Q \to M_k/k^{\times}$ and $P \to M_{X, \bar{x}}/\OO^{\times}_{X,\bar{x}}$
are bijective.

\item
The diagram
\[
 \begin{CD}
 Q @>>> P \\
 @VVV @VVV \\
 M_k @>>> M_{X,\bar{x}}
 \end{CD}
\]
is commutative.

\item
$k \otimes_{k[Q]} k[P] \to \OO_{X,\bar{x}}$ is smooth. 
\end{enumerate}
Moreover, using the good chart $(Q \to M_k, P \to M_{X,\bar{x}}, Q \to P)$,
we can determine the local structure in the following ways:
\begin{enumerate}
\renewcommand{\labelenumi}{(\arabic{enumi})}
\item
If $\mult_x(X) = 1$,
then $Q \to P$ splits and $P \simeq Q \times \NN^r$ for some $r$.

\item
If $\mult_x(X) = 2$, then we have one of the following:
\begin{enumerate}
\renewcommand{\labelenumii}{(\arabic{enumi}.\arabic{enumii})}
\item
If $Q \to P$ does not split, then $P$ is of semistable type over $Q$.

\item
If $Q \to P$ splits, then $\operatorname{char}(k) \not= 2$ and
$\widehat{\OO}_{X, x}$ is canonically isomorphic to 
$k\lformal X_1, \ldots, X_n\rformal /(X_1^2 - X_2^2)$.
\end{enumerate}

\item
If $\mult_x(X) \geq 3$, then $Q \to P$ does not split and
$P$ is of semistable type over $Q$.
\end{enumerate}
For the definition of a monoid of semistable type, see
\S\rom{\ref{sec:monoid:semistable:type}}.
\end{Theorem}

By using the above local structure result,
we can see the uniqueness of a log morphism over the fixed
scheme morphism, namely, we have the following:

\begin{Theorem}
\label{thm:intro:unqueness:log:hom:semistable}
Let $X$ and $Y$ be semistable varieties over $k$
and let $M_X$ and  $M_Y$ be fine log structures of $X$ and $Y$ over $M_k$
respectively
such that $(X, M_X)$ and $(Y, M_Y)$ are 
log smooth and integral over $(\Spec(k), M_k)$.
Let $\Supp(M_Y/M_k)$ be the union of $\Sing(Y)$
and the boundaries of the log structure of $M_Y$ over $M_k$.
Let $\phi : X \to Y$ be a morphism over $k$ such that 
$\phi(X') \not\subseteq \Supp(M_Y/M_k)$
for any irreducible component $X'$ of $X$.
If $(\phi, h) : (X, M_X) \to (Y, M_Y)$ and
$(\phi, h') : (X, M_X) \to (Y, M_Y)$ are morphisms of log schemes over
$(\Spec(k), M_k)$, then $h = h'$.
\end{Theorem}

By virtue of this theorem, the non-classical case can be reduced
to the classical case, so that we complete the proof of the theorem. 

Finally, we would like to express our sincere thanks
to Prof. Kazuya Kato for telling us the fantastic finiteness problem.

\renewcommand{\thesubsubsection}{\arabic{subsubsection}}

\bigskip
\subsection*{Conventions and terminology}
Here we will fix several conventions and terminology for this paper.

\subsubsection{}
\label{subsubsec:log:structure}
Throughout this paper, we work within the logarithmic
structures in the sense of J.-M Fontaine, L. Illusie, 
and K. Kato. For the details, we refer to \cite{KatoLog}.
All log structures on schemes are considered 
with respect to the etale topology.
We often denote the log structure on a scheme
$X$ by $M_{X}$ and the quotient $M_{X}/\mathcal{O}_{X}^{\times}$
by $\overline{M}_{X}$.

\subsubsection{}
\label{subsub:natural:number}
We denote by $\NN$ the set of natural integers. Note that $0 \in \NN$.
For $I = (a_1, \ldots, a_n) \in \NN^n$,
we define $\Supp(I)$ and $\deg(I)$ to be
\[
\Supp(I) = \{ i \mid a_i > 0 \}
\quad\text{and}\quad
\deg(I) = \sum_{i=1}^n a_i.
\]
The $i$-th entry of $I$ is denoted by $I(i)$, i.e., $I(i) = a_i$.
For $I, J \in \NN^n$,
a partial order $I \geq J$ is defined by
$I(i) \geq J(i)$ for all $i = 1, \ldots, n$.
The non-negative number $g$ with
$g \ZZ  = \ZZ I(1) + \cdots + \ZZ I(n)$ is denoted by
$\gcm(I)$.

\subsubsection{}
\label{subsub:monoid}
Here let us briefly recall some generalities on monoids.
All monoids in this paper are commutative with the unit element.
The binary operation of a monoid is often written additively.
We say a monoid $P$ is {\em finitely generated} if there are
$p_1, \ldots, p_n$ such that 
$P = \NN p_1 + \cdots + \NN p_r$.
Moreover, $P$ is said to be {\em integral} if
$x + z = y + z$ for $x, y, z \in P$, then $x = y$.
An integral and finitely generated monoid is said to be {\em fine}.
We say $P$ is {\em sharp} if $x + y = 0$ for $x, y \in P$,
then $x = y = 0$. For a sharp monoid $P$, an element $x$ of $P$
is said to be {\em irreducible} if
$x = y + z$ for $y, z \in P$, then either $y = 0$ or $z = 0$.
It is well known that if $P$ is fine and sharp,
then there are only finitely many irreducible elements and
$P$ is generated by irreducible elements
(cf. Proposition~\ref{prop:fine:sharp:irreducible}).
If $k$ is a field and $P$ is a sharp monoid,
then $M = \bigoplus_{x \in P \setminus \{ 0\}} k \cdot x$ forms
the maximal ideal of $k[P]$. This $M$ is called
{\em the origin of $k[P]$}. 
An integral monoid $P$ is said to be {\em saturated} 
if $n x \in P$ for $x \in P^{gr}$ and $n > 0$,
then $x \in P$, where $P^{gr}$ is the Grothendieck group
associated with $P$.
A homomorphism $f : Q \to P$ of monoids is said to be {\em integral} if
$f(q) + p = f(q') + p'$ for $p, p' \in P$ and $q, q' \in Q$,
then there are $q_1, q_2 \in Q$ and $p'' \in P$
such that $q + q_1 = q' + q_2$, $p = f(q_1) + p''$
and $p' = f(q_2) + p''$.
Note that an integral homomorphism of sharp monoids is injective.
Moreover, we say an injective homomorphism
$f : Q \to P$ {\em splits} if there is a submonoid $N$ of $P$ with
$P = f(Q) \times N$.
Finally, let us recall {\em congruence relation}.
A congruence relation  on a monoid
$P$ is a subset 
$S \subset P \times P$ which is both a submonoid
and a set-theoretic equivalence relation. We say
that a subset $T \subset S$ {\em generates the
congruence relation} $S$ if $S$ is the smallest
congruence relation on $P$ containing $T$.
Let $S$ be an equivalent relation on $P$.
It is easy to see that
$P \rightarrow P/S$ gives rise a
structure of a monoid on $P/S$ if and only if
$S$ is a congruence relation.

\subsubsection{}
\label{subsub:pushout}
Let $P$ and $Q$ be monoids and let
$f : \NN \to P$ and $g : \NN \to Q$ be homomorphisms with
$p = f(1)$ and $q = g(1)$.
Let $P \times_{\NN} Q$ be the pushout of $f : \NN \to P$ and $g : \NN \to Q$:
\[
\begin{CD}
\NN @>>> Q \\
@VVV @VVV \\
P @>>> P \times_{\NN} Q
\end{CD}
\]
We denote this pushout $P \times_{\NN} Q$ by $P \times_{(p,q)} Q$.

\subsubsection{}
\label{subsub:monomial}
Let $k$ be a field and $R$ be either the ring of polynomials of $n$-variables
over $k$, or the ring of formal power series of $n$-variables over $k$,
that is, $R = k[X_1, \ldots, X_n]$ or $k\lformal X_1, \ldots, X_n\rformal $.
For $I \in \NN^n$, we denote the monomial $X_1^{I(1)} \cdots X_n^{I(n)}$
by $X^I$.

\subsubsection{}
\label{subsub:log:monomial}
Let $P$ be a monoid, $p_1, \ldots, p_n \in P$ and $I \in \NN^n$.
For simplicity, $\sum_{i=1}^n I(i) p_i$ is often denoted by
$I \cdot p$.

\subsubsection{}
\label{subsubsec:underlining:log:structure}
Let $(X, M_X)$ be a log scheme and $\alpha : M_X \to \OO_X$ the
structure homomorphism. Then, 
$\alpha(M_X) \setminus \{\text{zero divisors of $\OO_X$}\}$
give rise to a log structure because
\[
 \OO_X^{\times} \subseteq \alpha(M_X)
 \setminus \{\text{zero divisors of $\OO_X$}\}.
\]
$\alpha(M_X) \setminus \{\text{zero divisors of $\OO_X$}\}$
is called {\em the underlining log structure} of $M_X$ 
and is denoted by $M^u_X$.
Let $f : (X, M_X) \to (Y, M_Y)$ be a morphism of log schemes such that
one of the following conditions is satisfied:
\begin{enumerate}
\renewcommand{\labelenumi}{(\arabic{enumi})}
\item
$X \to Y$ is flat.

\item
$X$ and $Y$ are integral schemes and $X \to Y$ is a dominant
morphism.
\end{enumerate}
Then we have the induced morphism
$f^{u} : (X, M^{u}_X) \to (Y, M^{u}_Y)$.

\subsubsection{}
\label{subsubsec:rational:map}
Let $X$ and $Y$ be reduced noetherian schemes.
Let $\phi : X \dasharrow Y$ be a rational map.
We say $\phi$ is {\em dominant} (resp. {\em separably dominant})
if for any irreducible component $X'$ of $X$, there
is an irreducible component $Y'$ of $Y$ such that
$\phi(X') \subseteq Y'$ and
the induced rational map $\phi' : X' \dasharrow  Y'$
is dominant (resp. dominant and separable).
Moreover, we say $\phi$ is {\em defined in codimension one} if
there is a dense open set $U$ of $X$ such that
$\phi$ is defined over $U$ and $\codim(X \setminus U) \geq 2$.

Let $f : X \to T$ and $g : Y \to T$ be
morphisms of reduced noetherian schemes.
A rational map $\phi : X \dasharrow Y$ is called
a {\em relative rational map}
if there is a dense open set $U$ of $X$ such that
$\phi$ is defined on $U$, $\phi : U \to Y$ is a morphism over $T$
(i.e., $f = g \cdot \phi$) and $X_t \cap U \not= \emptyset$ for all
$t \in T$.

\subsubsection{}
\label{subsub:semistable:variety}
Let $k$ be an algebraically closed field and
$X$ a reduced algebraic scheme over $k$.
We say $X$ is a {\em semistable variety} if
for any closed point $x \in X$,
the completion $\widehat{\OO}_{X, x}$ at $x$ is isomorphic to
the ring of the type $k\lformal X_1, \ldots, X_n \rformal /(X_1 \cdots X_l)$.

\subsubsection{}
\label{subsub:big:line:bundle}
Let $k$ be an algebraically closed field.
Let $X$ be a proper reduced algebraic scheme over $k$ and
$H$ a line bundle on $X$.
We say $H$ is {\em very big} if
there is a dense open set $U$ of $X$ such that
$H^0(X, H) \otimes \OO_X \to H$ is surjective on $U$ and
the induced rational map $X \dasharrow  \PP(H^0(X, H))$
is birational to the image.
Moreover, $H$ is said to be {\em big} if $H^{\otimes m}$ is very big
for some positive integer $m$.

\renewcommand{\theTheorem}{\arabic{section}.\arabic{Theorem}}
\renewcommand{\thesubsubsection}{\arabic{section}.\arabic{subsection}.\arabic{subsubsection}}

\section{Existence of a good chart on a generalized semistable variety}
\label{sec:exist:good:chart}
Let $k$ be an algebraically closed field and
$X$ an algebraic scheme over $k$.
We say $X$ is a {\em generalized semistable variety}
if, for any closed point $x$ of $X$, the completion $\hat{\OO}_{X, x}$
of $\OO_{X, x}$ is isomorphic to a ring of the following type:
\[
k\lformal T_1, \ldots, T_e \rformal/(T^{A_1}, \ldots, T^{A_l}),
\]
where $A_1, \ldots, A_l$ are elements of  $\NN^e \setminus \{ 0 \}$
such that
$A_i(j)$ is either $0$ or $1$ for all $i, j$.
Note that a generalized semistable variety is a reduced
scheme (cf. Lemma~\ref{lem:reduced:gen:semi:stable}).

Let $M_k$ and $M_X$ be fine log structures on $\Spec(k)$ and $X$
respectively. We assume that
$(X, M_k)$ is log smooth
and integral over $(\Spec(k), M_k)$.
Since the map $x \mapsto x^n$ on $k$ is surjective
for any positive integer $n$,
we can see that $M_k \to \overline{M}_k$ splits.
Thus, there are a fine and sharp monoid $Q$ and
a chart $\pi_Q : Q \to M_k$ such that
$Q \to M_k \to \overline{M}_k$ is bijective.

Next, let us choose a closed point $x$ of $X$.
In the case where $X$ is a generalized semistable variety,
we would like to construct a chart
$\pi_P : P \to M_{X, \bar{x}}$ together with a homomorphism
$f : Q \to P$ such that $P \to 
M_{X, \bar{x}} \to \overline{M}_{X, \bar{x}}$ is
bijective, the natural morphism $X \to
\Spec(k) \times_{k[Q]} \Spec(k[P])$ is smooth and
the following diagram is commutative:
\[
 \begin{CD}
  Q @>{f}>> P \\
  @V{\pi_Q}VV @VV{\pi_P}V \\
  M_k @>>> M_{X, \bar{x}}.
 \end{CD}
\]
Then, a triple $(Q \to M_k, P \to M_{X, \bar{x}},
Q \to P)$ is called {\em a good chart of
$(X, M_X) \to (\Spec(k), M_k)$ at $x$}.
For this purpose, we need to see the following theorem.

\begin{Theorem}
\label{thm:log:smooth:coprime:torsion}
Let $\mu : (X, M_X) \to (Y, M_Y)$ be a log smooth and integral
morphism of fine log schemes.
Let $x \in X$ and $y = \mu(x)$.
Let $k$ be the algebraic closure of the residue field at $x$ and
$\eta : \Spec(k) \to X \overset{\mu}{\longrightarrow} Y$ the induced morphism.
If $X \times_{Y} \Spec(k)$ is a generalized semistable variety over $k$,
then the torsion part of $\Coker(\overline{M}^{gr}_{Y, \bar{y}} \to
\overline{M}^{gr}_{X, \bar{x}})$ is a finite group of order invertible in
$\OO_{X, \bar{x}}$.
\end{Theorem}

\Proof
Let us begin with the following lemma.

\begin{Lemma}
\label{lem:chart:at:residue:field}
Let $(X, M_X)$ be a log scheme with a fine log structure. Then, we have the following:
\begin{enumerate}
\renewcommand{\labelenumi}{(\arabic{enumi})}
\item
Let $\pi : P \to \rest{M_X}{U}$ be a local chart of $M_X$ on an etale neighborhood $U$.
Then, for $x \in U$, the natural map $P/\pi^{-1}(\OO^{\times}_{X, \bar{x}}) \to\overline{M}_{X, \bar{x}}$
is bijective.

\item
Let $k$ be a separably closed field and 
$\eta: \Spec(k) \to X$ a geometric point.
Then, the natural homomorphism $\overline{M}_{X, \bar{x}} \to \overline{\eta^*(M_X)}$
is an isomorphism, where $x$ is the image of $\eta$.
\end{enumerate}
\end{Lemma}

\Proof
(1) The surjectivity of $P/\pi^{-1}(\OO^{\times}_{X, \bar{x}}) \to\overline{M}_{X, \bar{x}}$
is obvious. Let us assume that $\pi(a) \equiv \pi(b) \mod \OO^{\times}_{X, \bar{x}}$.
Then, there is $u \in \OO^{\times}_{X, \bar{x}}$ with $\pi(a) = \pi(b) \cdot u$.
Since $\pi : P \to \rest{M_X}{U}$ is a chart, we have the natural isomorphism
\[
P \times_{\pi^{-1}(\OO^{\times}_{X, \bar{x}})} \OO^{\times}_{X,\bar{x}} \overset{\sim}{\longrightarrow}
M_{X,\bar{x}}.
\]
Thus, there are $\alpha, \beta \in \pi^{-1}(\OO^{\times}_{X, \bar{x}})$ such that
\[
(a, 1) + (\alpha, \pi(\alpha)^{-1}) = (b, u) + (\beta, \pi(\beta)^{-1}).
\]
In particular, $a + \alpha = b + \beta$. Thus, $x \equiv y \mod \pi^{-1}(\OO^{\times}_{X, \bar{x}})$.

\medskip
(2) Let $P \to M_X$ be a local chart around $x$ and $\alpha : P \to \OO_X$ the induced
homomorphism. Note that $M_X$ is isomorphic to the associated log structure $P^a$.
Let $\alpha' : P \to k$ be a homomorphism given by the compositions:
\[
P \overset{\alpha}{\longrightarrow} \OO_{X, \bar{x}} \to \kappa(\bar{x}) \hookrightarrow k,
\]
where $\kappa(\bar{x})$ is the residue field at $\bar{x}$.
Then, by \cite[(1.4.2)]{KatoLog}, $\eta^*(M_X)$ is
the associated log structure of $\alpha' : P \to k$.
Therefore, we get the following commutative diagram:
\[
\begin{CD}
P @= P \\
@VVV @VVV \\
\overline{M}_{X,\bar{x}} @>>> \overline{\eta^*(M_X)}.
\end{CD}
\]
On the other hand,
\begin{align*}
a \in \alpha^{-1}(\OO^{\times}_{X, \bar{x}})  &
 \Longleftrightarrow \alpha(a) \in \OO^{\times}_{X, \bar{x}}
 \Longleftrightarrow \text{$\alpha(a) \not= 0$ in $\kappa(\bar{x})$} \\
 & \Longleftrightarrow \alpha'(a) \not= 0 
 \Longleftrightarrow  a \in {\alpha'}^{-1}(k^{\times}).
\end{align*}
Therefore, $\alpha^{-1}(\OO^{\times}_{X, \bar{x}}) = {\alpha'}^{-1}(k^{\times})$.
Thus, (1) implies (2).
\QED

\bigskip
\addtocounter{Theorem}{1}
Let us go back to the proof of Theorem~\ref{thm:log:smooth:coprime:torsion}.
We denote $X \times_Y \Spec(k)$ by $X'$.
Then, we have the following commutative diagram:
\[
\begin{CD}
X @<{\tilde{\eta}}<< X' \\
@V{\mu}VV @VV{\mu'}V \\
Y @<{\eta}<< \Spec(k).
\end{CD}
\]
Note that the natural morphism 
$\eta' : \Spec(k) \to X'$ gives rise to a section of
$\mu' : X' \to \Spec(k)$.
Let $x'$ be the image of $\eta'$. We consider the natural commutative diagram:
\[
\begin{CD}
\overline{M}_{X,\bar{x}} @>>> \overline{\tilde{\eta}^*(M_X)}_{X', \bar{x}'} @>>> \overline{{\eta'}^*(\tilde{\eta}^*(M_X))} \\
@AAA @AAA @AAA \\
\overline{M}_{Y,\bar{y}} @>>> \overline{\eta^*(M_Y)} @= \overline{\eta^*(M_Y)} 
\end{CD}
\]
By (2) of Lemma~\ref{lem:chart:at:residue:field},
\[
\overline{M}_{Y,\bar{y}} \to \overline{\eta^*(M_Y)}\quad\text{and}\quad
\overline{\tilde{\eta}^*(M_X)}_{X', \bar{x}'} \to \overline{{\eta'}^*(\tilde{\eta}^*(M_X))}
\]
are bijective.
Moreover, since ${\eta'}^*(\tilde{\eta}^*(M_X)) = 
(\tilde{\eta} \cdot \eta')^*(M_X)$, the composition
\[
\overline{M}_{X,\bar{x}}  \to \overline{\tilde{\eta}^*(M_X)}_{X', \bar{x}'} 
\to \overline{{\eta'}^*(\tilde{\eta}^*(M_X))}
\]
is also bijective. Thus, we can see that
\[
\overline{M}_{X,\bar{x}} \to
\overline{\tilde{\eta}^*(M_X)}_{X', \bar{x}'} 
\]
is an isomorphism.
Moreover, $(X', \tilde{\eta}^*(M_X)) \to (\Spec(k), \eta^*(M_Y))$
is smooth and integral.
Thus, we may assume that $Y = \Spec(k)$, $X$ is a generalized semistable variety
over $k$ and $x$ is a closed point of $X$.

\bigskip
Clearly, we may assume that $p = \operatorname{char}(k) > 0$.
We can take a fine and sharp monoid $Q$ with
$M_k = Q \times k^{\times}$.
Let $f : Q \to M_{X, \bar{x}}$ and
$\bar{f} : Q \to \overline{M}_{X, \bar{x}}$ be the canonical
homomorphisms.

Let us choose $t_1, \ldots, t_r \in M_{X, \bar{x}}$ such that
$d\log(t_1), \ldots, d\log(t_r)$ form a free basis
of $\Omega^1_{X/k, \bar{x}}(\log(M_X/M_k))$.
Then, in the same way as in \cite[(3.13)]{KatoLog},
we have the following:
\begin{enumerate}
\renewcommand{\labelenumi}{(\roman{enumi})}
\item
If we set $P_1 = \NN^r \times Q$ and
a homomorphism $\pi_1 : P_1 \to M_{X, \bar{x}}$
by 
\[
 \pi_1(a_1, \ldots, a_r, q) = a_1 t_1 + \cdots + a_r t_r + f(q),
\]
then there is a fine monoid $P$ such that
$P \supseteq P_1$, $P^{gr}/P_1^{gr}$ is a finite group of order
invertible in $\OO_{X, \bar{x}}$ and that
$\pi_1 : P_1 \to M_{X, \bar{x}}$ extends
to the surjective homomorphism $\pi : P \to M_{X, \bar{x}}$.
Moreover, $P$ gives a local chart around $x$.
Here we have the natural homomorphism
$h : Q \to P_1 \hookrightarrow P$. Then, the following diagram is commutative:
\[
\begin{CD}
Q @>{h}>> P \\
@VVV @VV{\pi}V \\
M_k @>>> M_{X, \bar{x}}.
\end{CD}
\]

\item
The natural morphism $g : X \to \Spec(k) \times_{\Spec(k[Q])} \Spec(k[P])$
is etale around $x$.
\end{enumerate}

Let $\bar{p}_1, \ldots, \bar{p}_e$ be all irreducible elements of $\overline{M}_{X, \bar{x}}$
not lying in the image $Q \to \overline{M}_{X, \bar{x}}$.
Let us choose $p_1, \ldots, p_e \in M_{X, \bar{x}}$ such that
the image of $p_i$ in $\overline{M}_{X,\bar{x}}$ is $\bar{p}_i$.
Let $\alpha : M_X \to \OO_X$ be the canonical homomorphism.
We set $z_i = \alpha(p_i)$ for $i=1, \ldots, e$. Then, we have the following:

\begin{Claim}
\label{Claim:prop:log:smooth:coprime:torsion:0}
$z_i \not= 0$ in $\OO_{X, \bar{x}}$ for all $i$.
\end{Claim}

Since $\bar{\pi} : P \to M_{X, \bar{x}} \to \overline{M}_{X,\bar{x}}$ 
is surjective, there are $p'_1, \ldots, p'_r \in P$ with 
$\bar{\pi}(p'_i) = \bar{p}_i$.
Let us choose $u_1, \ldots, u_a \in P$ such that
the kernel of $P^{gr} \to \overline{M}^{gr}_{X, \bar{x}}$
is generated by $u_1, \ldots, u_a$.
Note that 
$\pi(u_i) \in \OO^{\times}_{X, \bar{x}}$ and
$P$ is generated by $p'_1, \ldots, p'_r$, $u_1, \ldots, u_a$ and
$h(q)$ ($q \in Q$).
Let us consider a non-trivial congruence relation
\[
 I \cdot p' + J \cdot u + h(q)  = 
    I' \cdot p' + J' \cdot u + h(q'),
\]
where $I, I' \in \NN^r$, 
$J, J' \in \NN^a$, 
$q, q' \in Q$,
$\Supp(I) \cap \Supp(I') = \emptyset$ and
$\Supp(J) \cap \Supp(J') = \emptyset$
(See Conventions and terminology \ref{subsub:log:monomial}).
Let $\phi : k[Z_1, \ldots, Z_r, U_1, \ldots, U_a] \to
       k \otimes_{k[Q]} k[P]$
be the natural surjective homomorphism
given by $\phi(Z_i) = 1 \otimes p'_1$ and $\phi(U_j) = 1 \otimes u_j$.
Then, the kernel of $\phi$ is generated by elements of the type
\[
 \beta(q) \cdot Z^{I} \cdot U^{J} - 
    \beta(q') \cdot Z^{I'} \cdot U^{J'},
\]
where
\[
 \beta(q) = \begin{cases}
 1 & \text{if $q = 0$} \\
 0 & \text{if $q \not= 0$}.
 \end{cases}
\]
Here note that 
$I \cdot \bar{p} + \bar{f}(q) = 
I' \cdot \bar{p} +  \bar{f}(q')$
and $\bar{p}_i$'s are irreducible.
Thus, 
\[
 \beta(q) \cdot Z^{I} \cdot U^{J} - 
    \beta(q') \cdot Z^{I'} \cdot U^{J'}
\]
is equal to either
\[
 \pm Z^{I} \cdot U^{J}\quad  (\deg(I) \geq 2)
\]
or
\[
 Z^{I} \cdot U^{J} - Z^{I'} \cdot U^{J'} \quad
(\deg(I) \geq 2, \ \deg(I') \geq 2).
\]
Therefore,
\[
 \Ker(\phi) \subseteq (Z_1, \dots, Z_r)^2.
\]
Now let us consider a natural homomorphism
\[
 g^* : R = 
  k[Z_1, \ldots, Z_r, U_1, \ldots, U_a]/\ker(\phi) \to \OO_{X, \bar{x}}.
\]
Note that $g^*(\bar{Z}_i) = v_i \cdot z_i$ and 
$g^*(\bar{U}_j) = \alpha(\pi(u_j))$,
where $v_i \in \OO_{X, \bar{x}}$ and $\bar{Z}_i$ and $\bar{U}_j$
are the classes of $Z_i$ and $U_j$ in 
$k[Z_1, \ldots, Z_r, U_1, \ldots, U_a]/\ker(\phi)$ respectively.
Let $y = g(\bar{x})$. Then, since $\pi(u_j)$'s are
units, we can set $y = (\underbrace{0,\ldots,0}_r, c_1, \ldots, c_a )$,
where $c_1, \ldots, c_a \in k^{\times}$. 
Since $g$ is etale, $g^* : R_y \to \OO_{X,\bar{x}}$ is injective.
Thus, if $z_i = 0$, then
$Z_i \in \Ker(\phi)k[Z_1, \ldots, Z_r, U_1, \ldots, U_a]_y$.
This is a contradiction because
\[
 \Ker(\phi)k[Z_1, \ldots, Z_r, U_1, \ldots, U_a]_y \subseteq
 (Z_1, \ldots, Z_r)^2 k[Z_1, \ldots, Z_r, U_1, \ldots, U_a]_y.
\]

\bigskip
Note that $M_{X, \bar{x}}$ is generated by 
$p_1, \ldots, p_e$, $\OO^{\times}_{X,\bar{x}}$ and
the image of $Q$ in $M_{X,\bar{x}}$, so that, from now on,
we always choose $t_1, \ldots, t_r$ from elements of the following types:
\[
p_i u\  \ (u \in \OO^{\times}_{X, \bar{x}}, \ i=1, \ldots, e)
\quad\text{and}\quad
v \ \ (v \in \OO^{\times}_{X, \bar{x}}).
\]
We set $x_i = \alpha(t_i)$ for $i=1, \ldots, r$.

\begin{Claim}
\label{Claim:prop:log:smooth:coprime:torsion:1}
\begin{enumerate}
\renewcommand{\labelenumi}{(\alph{enumi})}
\item
$x_1^{a_1} \cdots x_r^{a_r} \not= 0$ for any non-negative
integers $a_1, \ldots, a_r$.

\item
If $x_1^{a_1} \cdots x_r^{a_r} = x_1^{a'_1} \cdots x_r^{a'_r}$
for non-negative integers $a_1, \ldots, a_r, a'_1, \ldots, a'_r$,
then $(a_1, \ldots, a_r) = (a'_1, \ldots, a'_r)$.
\end{enumerate}
\end{Claim}

Let $T_i$ be an element of $k \otimes_{k[Q]} k[P]$ arising from
$e_i = (0, \ldots, 1, \ldots, 0) \in \NN^r$ 
($i$-th standard basis of $\NN^r$), namely, $T_i = 1 \otimes e_i$.
As in the previous claim, let us choose $u_1, \ldots, u_a \in P$ such that
the kernel of $P^{gr} \to \overline{M}^{gr}_{X, \bar{x}}$
is generated by $u_1, \ldots, u_a$.
Let $P'$ be the submonoid of $P^{gr}$ generated by
$\pm e_1, \ldots, \pm e_r, \pm u_1, \ldots, \pm u_a$ and $P$.

First, let us see that $\bar{f} : Q \to \bar{\pi}(P')$ is integral.
We consider an equation
\[
p - I \cdot \bar{e} + \bar{f}(q) = p' - I' \cdot \bar{e} + \bar{f}(q'),
\]
where $p, p' \in \overline{M}_{X, \bar{x}}$, $q, q' \in Q$ and $I, I' \in \NN^r$.
Then,
\[
p + I' \cdot \bar{e} + \bar{f}(q) = p' + I \cdot \bar{e} + \bar{f}(q').
\]
Thus, since $Q \to \overline{M}_{X,\bar{x}}$ is integral,
there are $q_1, q_2 \in Q$ and $x \in P$ such that
\[
\begin{cases}
q + q_1 = q' + q_2,  \\
p + I' \cdot \bar{e} = \bar{f}(q_1) + x \\
p' + I \cdot \bar{e} = \bar{f}(q_2) + x.
\end{cases}
\]
Therefore,
\[
\begin{cases}
p - I \cdot \bar{e} = \bar{f}(q_1) + x - (I+I') \cdot \bar{e} \\
p' - I' \cdot \bar{e} = \bar{f}(q_2) + x - (I+I') \cdot \bar{e}.
\end{cases}
\]

Next let us see that the natural homomorphism 
$\nu : Q \times \ZZ^r \to P'$ given by $\nu(q,I) = f(q) + I \cdot e$
is integral.
For this purpose,
let us consider an equation
\[
x + \nu(q, I) = x' + \nu(q', I'),
\]
where $x, x' \in P'$, $q, q' \in Q$ and $I, I' \in \ZZ^r$.
Then, in $\bar{\pi}(P')$, we have
\[
\bar{x} + I \cdot \bar{e} + \bar{f}(q) = \bar{x}' + I' \cdot \bar{e} + \bar{f}(q').
\]
Thus, there are $q_1, q_2 \in Q$, $y \in P'$ and $J, J' \in \ZZ^a$ such that
\[
\begin{cases}
q + q_1 = q' + q_2 \\
x + I \cdot e  = \nu(q_1, 0) + J \cdot u + y \\
x' + I' \cdot e = \nu(q_2, 0) + J' \cdot u + y.
\end{cases}
\]
Therefore, using the equation $x + \nu(q, I) = x' + \nu(q', I')$,
we can see that $J \cdot u + y = J' \cdot u + y$.
Thus,
\[
x = \nu(q_1, -I) + z\quad\text{and}\quad
x' = \nu(q_2,-I') + z
\]
for some $z \in P'$ and 
\[
\nu(q, I) + \nu(q_1, -I) = \nu(q+ q_1, 0) = \nu(q'+ q_2, 0) = \nu(q', I') + \nu(q_2, -I').
\]

Therefore, by \cite[Proposition~(4.1)]{KatoLog},
$k[P']$ is flat over $k[Q \times \ZZ^r]$.
Moreover, since 
\[
 k \otimes_{k[Q]} k[P'] \simeq (k \otimes_{k[Q]} 
k[Q \times \ZZ^r]) \otimes_{k[Q \times \ZZ^r]} k[P'],
\]
the following diagram
\[
\begin{CD}
\Spec(k \otimes_{k[Q]} k[P']) @>>> \Spec(k[P']) \\
@VVV @VVV \\
\Spec(k \otimes_{k[Q]} k[Q \times \ZZ^r]) @>>>
\Spec(k[Q \times \ZZ^r])
\end{CD} 
\]
is Cartesian. Therefore,
\[
 \Spec(k \otimes_{k[Q]} k[P']) \to
 \Spec(k \otimes_{k[Q]} k[Q \times \ZZ^r]) = \Spec(k[\ZZ^r])
\]
is flat. In particular,
\[
\beta : k[\ZZ^r] = k \otimes_{k[Q]} k[Q \times \ZZ^r] 
\to k \otimes_{k[Q]} k[P']
\]
is injective because $k[\ZZ^r]$ is a integral domain.
Further, $\beta(Y_i) = T_i$ for $i=1, \ldots, r$,
where $k[\ZZ^r] = k[Y_1^{\pm}, \ldots, Y_r^{\pm}]$.

Let $U$ be an etale neighborhood at $x$ and
$V$ a non-empty open set of $\Spec(k \otimes_{k[Q]} k[P])$
such that $V = g(U)$ and
$g : U \to V$ is etale.
Moreover, we set $W = \Spec(k \otimes_{k[Q]} k[P'])$.
Then, $W$ is an open set of $\Spec(k \otimes_{k[Q]} k[P])$, i.e.,
\[
 W = \left\{ t \in \Spec(k \otimes_{k[Q]} k[P]) 
\mid T_i(t) \not = 0 \ \forall i\ \ 
(1 \otimes u_j)(t) \not= 0 \ \forall j\right\}.
\]
Let $\overline{W}$ be the closure of $W$.
Note that
\begin{multline*}
 \Spec(k \otimes_{k[Q]} k[P]) = \\
\overline{W} \cup \{ T_1 = 0\} \cup
\cdots \cup \{ T_r = 0 \} \cup \{ 1 \otimes u_1 = 0 \} \cup \cdots
\cup \{ 1 \otimes u_a = 0\}.
\end{multline*}
Moreover, if we set $y = g(\bar{x})$,
then $(1 \otimes u_j)(y) \not= 0$ for all $j$
because $\pi(u_j) \in \OO^{\times}_{X, \bar{x}}$.
Note that the local ring $(k \otimes_{k[Q]} k[P])_{y}$
is reduced because 
$g^* : (k \otimes_{k[Q]} k[P])_{y} \to \OO_{X,\bar{x}}$
is etale.
Therefore, if $y \not\in \overline{W}$, then
$T_i = 0$ in $(k \otimes_{k[Q]} k[P])_{y}$.
This contradicts to
Claim~\ref{Claim:prop:log:smooth:coprime:torsion:0} because
$g^*(T_i) = x_i$.
Thus, $y \in \overline{W}$.
Let us consider
\[
 \gamma : k[\ZZ^r] \overset{\beta}{\longrightarrow} 
 \OO_W \to \OO_{W \cap V} \overset{g^*}{\longrightarrow} 
 \OO_{g^{-1}(W \cap V)}.
\]
Then, $\gamma(Y_i) = x_i$. Further, $\gamma$ is injective 
because $\beta$ and $g^*$ are injective and
$k[\ZZ^r]$ is an integral domain.
Thus, we get the claim.
 
\bigskip
Here we choose $t_1, \ldots, t_r  \in M_{X, \bar{x}}$
with the following properties:
\begin{enumerate}
\renewcommand{\labelenumi}{(\arabic{enumi})}
\item
$t_i$ is equal to either $p_j u$ ($u \in \OO_{X, \bar{x}}$) or
a unit $v$ for all $i$.

\item
$d\log(t_1), \ldots, d\log(t_r)$ form a free basis
of $\Omega^1_{X/k, \bar{x}}(\log(M_X/M_k))$.

\item
If we replace
the non-unit $t_i \not\in \OO^{\times}_{X, \bar{x}}$ by
a unit $t'_i \in \OO^{\times}_{X, \bar{x}}$, then
\[
 d\log(t_1), \ldots, d\log(t'_i), \ldots, d\log(t_r)
\]
do not form a free basis
of $\Omega^1_{X/k, \bar{x}}(\log(M_X/M_k))$.
\end{enumerate}

\begin{Claim}
\label{Claim:prop:log:smooth:coprime:torsion:2}
For a non-unit $t_i$ and $u \in \OO^{\times}_{X, \bar{x}}$, 
\[
 d\log(t_1), \ldots, d\log(t_i u), \ldots, d\log(t_r)
\]
form a free basis
of $\Omega^1_{X/k, \bar{x}}(\log(M_X/M_k))$.
\end{Claim}

We set $d\log(u) = f_1 d\log(t_1) + \cdots + f_r d\log(t_r)$.
If $f_i \in \OO^{\times}_{X, \bar{x}}$, then $d\log(t_i)$
belongs to a submodule generated by 
\[
d\log(u), d \log(t_1), \ldots, d\log(t_{i-1}), d\log(t_{i+1}), \ldots,  d\log(t_r).
\]
Thus, $d\log(u), d\log(t_1), \cdots, d\log(t_{i-1}), d\log(t_{i+1}), \cdots, d\log(t_r)$
form a basis, so that $f_i$ belongs to the maximal
ideal of $\OO_{X, \bar{x}}$.
Therefore,
\[
d\log(t_i u) = (1 + f_i) d\log(t_i) + \sum_{j\not=i} f_j d\log(t_j).
\]
and $1 + f_i \in \OO^{\times}_{X,\bar{x}}$.
Thus, we get the claim.

\bigskip
Renumbering $t_1, \ldots, t_r$, we may assume that
\[
 \{ t_1, \ldots, t_s \} = \{ t_i \mid 
 \text{$t_i $ is not a unit}\}
\]

\begin{Claim}
\label{Claim:prop:log:smooth:coprime:torsion:3}
Let $a_1, \ldots, a_s, a'_1, \ldots, a'_s$ be non-negative integers
such that either $a_i$ or $a'_i$ is zero for all $i$. 
For $u \in \OO^{\times}_{X, \bar{x}}$,
if 
\[
x_1^{a_1} \cdots  x_s^{a_s} = u x_1^{a'_1} \cdots x_s^{a'_s},
\]
then $a_1 = \cdots = a_s = a'_1 = \cdots = a'_s = 0$ 
and $u = 1$.
\end{Claim}

We assume the contrary. 
Let us choose a non-negative integer $k$ such that
$a_i = p^k b_i$ and $a'_i = p^k b'_i$ for all $i$ and that
\[
 \gcm(b_1, \ldots, b_s, b'_1, \ldots, b'_s)
\]
is prime to $p$.
Then, by Lemma~\ref{lem:pth:root:unit}, there is 
$v \in \OO^{\times}_{X, \bar{x}}$ with
\[
 x_1^{a_1} \cdots  x_s^{a_s} = v^{p^k} x_1^{a'_1} \cdots x_s^{a'_s}.
\]
Moreover by our construction, replacing $v$ by $v^{-1}$,
we can find $b'_i$ prime to $p$.
Thus, there is $v' \in \OO^{\times}_{X, \bar{x}}$ with
${v'}^{b'_i} = v$.
Hence if we replace $t_i$ by $v't_i$, then
we have $x_1^{a_1} \cdots x_s^{a_s} = x_1^{a'_1} \cdots x_s^{a'_s}$.
Therefore, by Claim~\ref{Claim:prop:log:smooth:coprime:torsion:1} and
Claim~\ref{Claim:prop:log:smooth:coprime:torsion:2},
$a_1 = a'_1, \ldots, a_s = a'_s$, which implies that
$a_1 = \cdots = a_s = a'_1 = \cdots = a'_s = 0$.
This is a contradiction.

\medskip
\begin{Claim}
\label{Claim:prop:log:smooth:coprime:torsion:4}
$t_1, \ldots, t_s$ are linearly independent over $\ZZ$
in $\Coker(Q^{gr} \to \overline{M}^{gr}_{X, \bar{x}})$.
\end{Claim}

We assume that a non-trivial relation
$a_1 t_1 + \cdots + a_s t_s = 0$ ($a_1, \ldots, a_s \in \ZZ$)
in $\Coker(Q^{gr} \to \overline{M}^{gr}_{X, \bar{x}})$.
Let $\bar{t}_i$ be the class of $t_i$ in $\overline{M}_{X, \bar{x}}$.
Then,
$a_1 \bar{t}_1 + \cdots + a_s \bar{t}_s = \bar{f}(q)$ 
for some $q \in Q^{gr}$.
Renumbering $t_1, \ldots, t_s$, we may assume that
$a_1, \ldots, a_l > 0$ and $a_{l+1}, \ldots, a_s \leq 0$.
Thus, we have
\[
 b_1 \bar{t}_1 + \cdots + b_l \bar{t}_l + \bar{f}(q_1) =
 b_{l+1} \bar{t}_{l+1} + \cdots + b_s \bar{t}_s + \bar{f}(q_2)
\]
for some $q_1, q_1 \in Q$, 
where $b_1 = a_1, \ldots, b_l = a_l$ and
$b_{l+1} = -a_{l+1}, \ldots, b_s = -a_s$.
Since $\bar{f}$ is integral, there are $q_3, q_4 \in Q$,
$x \in M_{X, \bar{x}}$ and 
$u, u' \in \OO^{\times}_{X,\bar{x}}$ with
\[
 \begin{cases}
  q_1 + q_3 = q_2 + q_4 \\
  b_1 t_1 + \cdots + b_l t_l = f(q_3) + x + u\\
  b_{l+1} t_{l+1} + \cdots + b_s t_s = f(q_4) + x + u'.
 \end{cases}
\]
Thus, if $q_3 \not= 0$, then
$x_1^{b_1} \cdots x_s^{b_s} = 0$, which contradicts to
Claim~\ref{Claim:prop:log:smooth:coprime:torsion:1}.
Therefore, $q_3 = 0$. In the same way, $q_4 = 0$.
Thus, we get
\[
 b_1 t_1 + \cdots + b_l t_l = 
 b_{l+1} t_{l+1} + \cdots + b_s t_s + v_0
\]
for some $v_0 \in \OO^{\times}_{X, \bar{x}}$.
Thus, $x_1^{b_1} \cdots x_l^{b_l} = v_0
x_{l+1}^{b_{l+1}} \cdots x_{s}^{b_s}$.
Therefore, by Claim~\ref{Claim:prop:log:smooth:coprime:torsion:3},
$b_1 = \cdots = b_l = b_{l+1} = \cdots = b_s = 0$.
This is a contradiction.

\bigskip
Let 
$\lambda : P^{gr} \to \overline{M}^{gr}_{X,\bar{x}}$
be the natural surjective homomorphism and
\[
 \lambda' : \Coker(Q^{gr} \to P^{gr}) \to 
\Coker(Q^{gr} \to \overline{M}^{gr}_{X,\bar{x}})
\]
the induced homomorphism.
Then, by using Claim~\ref{Claim:prop:log:smooth:coprime:torsion:4},
if we set
\[
 T = \Coker\left( \ZZ t_1 \oplus \cdots \oplus \ZZ t_r \to
 \Coker(Q^{gr} \to P^{gr})\right)
\]
and
\[
 T' = \Coker\left( \ZZ t_1 \oplus \cdots \oplus \ZZ t_s \to 
 \Coker(Q^{gr} \to \overline{M}^{gr}_{X,\bar{x}})\right),
\]
then we have the following commutative diagram:
\[
 \begin{CD}
 0 @>>> \ZZ t_1 \oplus \cdots \oplus \ZZ t_r @>>> 
 \Coker(Q^{gr} \to P^{gr}) @>>> T @>>> 0 \\
 @. @VV{\text{projection}}V @VV{\lambda'}V @VVV @. \\
 0 @>>> \ZZ t_1 \oplus \cdots \oplus \ZZ t_s @>>> 
 \Coker(Q^{gr} \to \overline{M}^{gr}_{X,\bar{x}})
 @>>> T' @>>> 0 \\
 @. @VVV @VVV @VVV @. \\
 @. 0 @. 0 @. 0 @. 
 \end{CD}
\]
Here $T$ is a torsion group of order prime to $p$.
Therefore, we get our assertion.
\QED

\begin{Lemma}
\label{lem:pth:root:unit}
Let $X$ be a generalized semistable variety over an algebraically closed field $k$
of characteristic $p > 0$
and $x$ a closed point of $X$.
Let $\OO_{X, \bar{x}}$ be the local ring at $x$ in the etale topology.
Let $H$ and $G$ be elements of $\OO_{X, \bar{x}}$ and
$u \in \OO^{\times}_{X, \bar{x}}$.
If $H^{p^k} u = G^{p^k}$, then there is $v \in \OO^{\times}_{X, \bar{x}}$
with $(H v)^{p^k} = G^{p^k}$. 
\end{Lemma}

\Proof
By Artin's approximation theorem, it is sufficient to find $v$
in $\hat{\OO}_{X, \bar{x}}$.
Since $X$ is a generalized semistable variety, we can set
\[
 \hat{\OO}_{X, \bar{x}} = 
 k\lformal T_1, \ldots, T_e \rformal/(T^{I_1},  \ldots,  T^{I_l}),
\]
where $I_1, \ldots, I_l \in \NN^e \setminus \{ 0\}$.
We set 
\[
\Omega =  \bigcup_{i=1}^l (I_i + \NN^e)\quad\text{and}\quad
\Sigma = \NN^e \setminus \bigcup_{i=1}^l (I_i + \NN^e)
\]
and $\Sigma_k = \{ I \in \Sigma \mid 
p^k | I(i) \ \forall i \}$.
Then, any elements of $\hat{\OO}_{X, \bar{x}}$ can be uniquely
written as a form
\[
 \sum_{I \in \Sigma} \alpha_I T^I.
\]
We set $u = \sum_{I \in \Sigma} a_I T^I$ and
$H = \sum_{I \in \Sigma} b_I T^I$.
Moreover, we set
\[
 u' = \sum_{I \in \Sigma_k} a_I T^I\quad\text{and}\quad
 u'' =  \sum_{I \not\in \Sigma_k} a_I T^I.
\]
Then, $u = u' + u''$ and there is a unit $v$ with
$v^{p^k} = u'$.
Thus, $H^{p^k}u'' = (G - H v)^{p^k}$. 
Therefore,
\[
  (G - H v)^{p^k} = \left(\sum_{I \in \Sigma}b_I^{p^k}T^{p^k I}\right)
\left(\sum_{I \not\in \Sigma_k} a_I T^I\right).
\]
Even if we delete the terms $T^J$ with $J \in \Omega$,
the left hand side of the above equations consists
of the terms $T^J$ with $J \in \Sigma_k$ and
the right hand side does not contain the terms $T^J$ with $J \in \Sigma_k$.
Thus, $(G - H v)^{p^k} = 0$.
\QED

\begin{Corollary}
\label{cor:good:pair:chart}
Let $X$ be a generalized semistable variety over an
algebraically closed field $k$.
Let $M_k$ and $M_X$ be fine log structures on $\Spec(k)$ and $X$
respectively. We assume that
$(X, M_X)$ is log smooth
and integral over $(\Spec(k), M_k)$.
Let $Q$ be a fine and sharp monoid with $M_k \simeq Q \times k^{\times}$
and $\pi_Q : Q \to M_k$ the composition of
$Q \to Q \times k^{\times}$ \rom{(}$q \mapsto (q,1)$\rom{)} 
and $Q \times k^{\times} \overset{\sim}{\longrightarrow} M_k$.
Then,
there is a fine and sharp monoid $P$ together with
homomorphisms $\pi_P : P \to M_{X, \bar{x}}$ and
$f : Q \to P$ such that a triple
$(\pi_Q : Q \to M_k, \ \pi_P : P \to M_{X, \bar{x}}, \ f : Q \to P)$
is a good chart of $(X, M_X) \to (\Spec(k), M_k)$ at $x$, namely,
the following properties are satisfied:
\begin{enumerate}
\renewcommand{\labelenumi}{(\arabic{enumi})}
\item
The diagram
\[
\begin{CD}
Q @>{f}>>  P \\
@V{\pi_Q}VV @VV{\pi_P}V \\
M_k @>>> M_{X, \bar{x}}
\end{CD}
\]
is commutative.

\item
The homomorphism
$P \to M_{X, \bar{x}} \to \overline{M}_{X, \bar{x}}$ is an isomorphism.

\item
The natural morphism $g : X \to \Spec(k) \times_{\Spec(k[Q])} \Spec(k[P])$
is smooth in the usual sense.
\end{enumerate}
\end{Corollary}

\Proof
This is a corollary of Theorem~\ref{thm:log:smooth:coprime:torsion}, 
Proposition~\ref{prop:tame:good:chart} and
Proposition~\ref{prop:log:smooth:good:chart}.
\QED

Finally let us consider the following lemma,
which is needed to see that a generalized variety is a reduced scheme.

\begin{Lemma}
\label{lem:reduced:gen:semi:stable}
Let $k\lformal T_1, \ldots, T_e \rformal$ be the ring
of formal power series over $k$.
Let $A_1, \ldots, A_l$ be elements of  $\NN^e \setminus \{ 0 \}$
such that
$A_i(j)$ is either $0$ or $1$ for all $i, j$.
Let $I$ be an ideal of $k\lformal T_1, \ldots, T_e \rformal$
generated by $T^{A_1}, \ldots, T^{A_l}$.
Then, $I$ is reduced, i.e.,
$\sqrt{I} = I$.
\end{Lemma}

\Proof
We prove this by induction on $e$.
If $e=1$, our assertion is obvious, so that we assume that $e > 1$.
Let $f \in \sqrt{I}$. Then, there is $n > 0$ with
$f^n \in I$.
It is easy to see that there are 
$a_i
\in k\lformal T_1, \ldots, T_{i-1}, T_{i+1},\ldots,T_e \rformal$
and $b \in k\lformal T_1, \ldots, T_e \rformal$ with
\[
 f = a_1 + T_1 a_2 + \cdots + T_1 \cdots T_{i-1} a_i + \cdots +
 T_1 \cdots T_{e-1} a_e + T_1 \cdots T_e b.
\]
Then, $f(0, T_2, \ldots, T_e) = a_1 \in 
k\lformal T_2,\ldots,T_e \rformal$.
If $1 \in \Supp(A_i)$ for all $i$, then
\[
f(0, T_2, \ldots, T_e)^n = 0.
\]
Thus, $a_1 = 0$. In particular, $a_1 \in I$.
Otherwise, 
\[
 a_1^n = f(0, T_2, \ldots, T_e)^n \in \sum_{1 \not\in \Supp(A_i)}
T^{A_i} k\lformal T_2,\ldots,T_e \rformal.
\]
Thus, by hypothesis of induction, $a_1 \in I$.
Therefore, $(f - a_1)^n \in I$.
Note that $(f-a_1)(T_1, 0, T_3, \ldots, T_e) = T_1 a_2$.
Thus, in the same way as before, we can see that $T_1 a_2 \in I$.
Hence, $(f - a_1 - T_1 a_2)^n \in I$.
Proceeding with the same argument,
$T_1 \cdots T_{i-1} a_i \in I$ for all $i$.
On the other hand, $T_1 \cdots T_e \in I$.
Therefore, $f \in I$.
\QED

\section{Monoids of semistable type}
\label{sec:monoid:semistable:type}
In this section, we consider a monoid of semistable type.
First of all, let us give its definition.
Let $f : Q \to P$ be an integral homomorphism of fine and sharp monoids
with $Q \not= \{ 0\}$.
We say $P$ is {\em of semi-stable type}
\[
(r, l, p_1, \ldots, p_r, q_0, b_{l+1}, \ldots, b_{r})
\]
{\em over} $Q$ if the following conditions are satisfied:
\begin{enumerate}
\renewcommand{\labelenumi}{(\arabic{enumi})}
\item
$r$ and $l$ are positive integers with $r \geq l$, 
$p_1, \ldots, p_r \in P$, $q_0 \in Q \setminus \{ 0\}$, and
$b_{l+1}, \ldots, b_{r}$ are non-negative integers.

\item
$P$ is generated by $f(Q)$ and $p_1, \ldots, p_r$. 
The submonoid of $P$ generated by $p_1, \ldots, p_r$ in $P$, 
which is denoted by $N$,  is canonically isomorphic to $\NN^r$, namely,
a homomorphism $\NN^r \to N$ given by $(t_1, \ldots, t_r) \mapsto
\sum_i t_i p_i$ is an isomorphism.

\item
We set $\Delta_l, B \in \NN^r$ as follows:
\[
\Delta_l = (\underbrace{1, \ldots, 1}_l, \underbrace{0,\ldots, 0}_{r-l})
\quad\text{and}\quad
B = (\underbrace{0, \ldots, 0}_l, b_{l+1}, \ldots, b_r).
\]
Then, $\Delta_l \cdot p = f(q_0) + B \cdot p$, i.e.,
$p_1 + \cdots + p_l = f(q_0) + \sum_{i > l} b_i p_i$
(cf. Conventions and terminology \ref{subsub:log:monomial}).

\item
If we have a relation 
\[
I \cdot p = f(q) + J \cdot p\quad(I, J \in \NN^r)
\]
with $q \not= 0$,
then $I(i) > 0$ for all $i=1, \ldots, l$.
\end{enumerate}

\begin{Remark}
In the case where $l = 1$,
by using (2) of the following proposition,
we can see $P = f(Q) \times \NN p_2 \times \cdots \times \NN p_r$.
Conversely, if $P$ has a form $f(Q) \times \NN^{r-1}$ and $Q \not= \{ 0 \}$, then
$P$ is of semistable type in the following way:
Let $q_0$ be an irreducible element of $Q$ and $p_1 = f(q_0)$.
Let $e_i$ be the standard basis of $\NN^{r-1}$.
We set $p_i = (0, e_{i-1})$ for $i = 2, \ldots, r$. Then, since $Q$ is sharp,
$\NN p_1 \simeq \NN$.
Thus, the submonoid generated by $p_1, \ldots, p_r$ in $P$ is isomorphic to $\NN^r$.
Finally, let us consider a relation
$\sum_i a_i p_i = f(q) + \sum_i c_i p_i$ with $q \not= 0$.
Then,
\[
f(a_1 q_0) + \sum_{i \geq 2} a_i p_i = f(q + c_1 q_0) + \sum_{i \geq 2} c_i p_i.
\]
Thus, $a_1 q_0 = q + c_1 q_0$. Hence, if $a_1 = 0$, then $q = 0$. Therefore, $a_1 > 0$.
\end{Remark}

First, let us see elementary properties of a monoid of semistable type.

\begin{Proposition}
\label{prop:relation:coprime}
Let $f : Q \to P$ be an integral homomorphism of fine and sharp monoids.
We assume that $P$ is of semi-stable type
\[
(r, l, p_1, \ldots, p_r, q_0, b_{l+1}, \ldots, b_{r})
\]
over $Q$. Then, we have the following:
\begin{enumerate}
\renewcommand{\labelenumi}{(\arabic{enumi})}
\item
Let  $I \cdot p = f(q) + J \cdot p$ \rom{(}$I, J \in \NN^r$\rom{)} be a relation
with $q \not= 0$.
Then, $q = n q_0$ for some $n \in \NN$.
Moreover, if $\Supp(I) \cap \Supp(J) = \emptyset$,
then $I = n \Delta_l$ and $J = n B$.

\item
Let us consider two elements
\[
f(q) + T \cdot p
\quad\text{and}\quad
f(q') + T' \cdot p
\]
of $P$ such that there are $i$ and $j$ with $1 \leq i, j\leq l$ and
$T(i) = T'(j) = 0$.
If $f(q) + T \cdot p =
f(q') + T' \cdot p$, then $q = q'$ and $T = T'$.

\item
Let $U$ \rom{(}resp. $V$\rom{)} be the submonoid of $P$ generated by $p_1, \ldots, p_l$
\rom{(}resp. $f(Q)$ and $p_{l+1}, \ldots, p_r$\rom{)}.
Then, $U \simeq \NN^l$, $V \simeq Q \times \NN^{r-l}$ and the natural homomorphism
\[
 U \times_{(\Delta_l \cdot p, \ f(q_0) + B \cdot p)} V \to P
\]
is bijective \rom{(}cf. Conventions and terminology
\rom{\ref{subsub:pushout}}\rom{)}.
\end{enumerate}
\end{Proposition}

\Proof
(1) First we assume that $\Supp(I) \cap \Supp(J) = \emptyset$.
We set 
\[
n = \min \{ I(1), \ldots, I(l) \}
\quad\text{and}\quad
I' = I - n\Delta_l.
\]
Then, $I'(i) = 0$ for some $i$ with $1 \leq i \leq l$ and
$I \cdot p   = n\Delta_l \cdot p + I' \cdot p$. Thus,
\[
f(n q_0) + (n B + I') \cdot p = f(q) + J \cdot p.
\]
Therefore, since $f : Q \to P$ is integral,
there are $q_1, q_2 \in Q$ and $T \in \NN^r$ such that
$n q_0 + q_1 = q + q_2$,
\[
 (n B + I') \cdot p = f(q_1) + T \cdot p
  \quad\text{and}\quad
  J \cdot p = f(q_2) + T \cdot p.
\]
Note that $(n B + I')(i) = 0$ for some $i$ ($1 \leq i \leq l$). 
Thus, $q_1 = 0$.
Moreover, since $\{ 1, \ldots, l \} \subseteq \Supp(I)$, we have
$\Supp(J) \subseteq \{ l+1, \ldots, r \}$, so that
$q_2 = 0$.
Therefore, $q = nq_0$ and
$(n B + I') \cdot p = J \cdot p$. In particular,
$n B + I' = J$.
Note that $(n B + I')(i) = I'(i)$ and $J(i) = 0$ for $i=1, \ldots, l$.
Thus, $I'(1) = \cdots = I'(l) = 0$.
We assume that $\Supp(I') \not= \emptyset$.
Let us choose $i \in \Supp(I')$. Then, $i > l$ and
$J(i) = 0$. Thus, $n B(i) + I'(i) = 0$, which implies
$I'(i) = 0$. This is a contradiction. Hence, $I' = 0$.
Therefore, $q = nq_0$, $I = n \Delta_l$ and $J = nB$.

\smallskip
Next let us consider a general case.
We define $T \in \NN^r$ by $T(i) = \min \{ I(i), J(i) \}$, and
we set $I'  = I  - T$ and $J' = J - T$.
Then, $I' \cdot p = f(q) + J' \cdot p$ and $\Supp(I') \cap \Supp(J') = \emptyset$.
Thus, we can see $q = n q_0$ for some $n \in \NN$.

\medskip
(2) Since $f : Q \to P$ is integral,
there are $q_1, q_2 \in Q$ and $h \in \NN p_1 + \cdots + \NN p_r$
such that $T \cdot p  = f(q_1) + h$, $T' \cdot p = f(q_2) + h$ and
$q + q_1 = q' + q_2$. Here $T(i) = 0$ for some $i = 1, \ldots, l$.
Thus, $q_1 = 0$. In the same way, $q_2 = 0$. Therefore, $q = q'$.
Hence $T \cdot p  =  T' \cdot p$.

\medskip
(3) By (2), it is easy to see that $U \simeq \NN^l$ and $V \simeq Q \times \NN^{r-l}$.
Let $I, I', J, J' \in \NN^r$ such that
$\Supp(I), \Supp(I') \subseteq \{ 1, \ldots, l \}$ and
$\Supp(J), \Supp(J') \subseteq \{ l+1, \ldots, r \}$.
It is sufficient to see that if
\[
I \cdot p + f(q) + J \cdot p = I' \cdot p + f(q') + J' \cdot p
\]
for some $q, q' \in Q$, then
\[
(I \cdot p, \ f(q) + J \cdot p) \sim (I' \cdot p, \ f(q') + J' \cdot p)
\]
in $U \times_{(\Delta_l \cdot p, \ f(q_0) + B \cdot p)} V$.
We set
\[
n = \min\{ I(1), \ldots, I(l) \}
\quad\text{and}\quad
n' = \min \{ I'(1), \ldots, I'(l) \}.
\]
Moreover, we set $T = I - n \Delta_l$ and $T' = I' - n' \Delta_l$.
Then
\[
(T + J + nB) \cdot p + f(q + nq_0) = (T' + J' + n'B) \cdot p + f(q' + n'q_0).
\]
Thus, by (2), $T + J + nB = T' + J' + n'B$ and $q + nq_0 = q' +n' q_0$.
In particular, $T = T'$ and $J + nB = J' + n'B$. Therefore, since
$(\Delta_l \cdot p,\ 0) \sim (0, f(q_0) + B \cdot p)$,
\begin{align*}
(I \cdot p, \ f(q) + J \cdot p) & = ((T + n \Delta_l) \cdot p, \ f(q) + J \cdot p) \\
& \sim (T \cdot p, \ f(q + nq_0) + (J + nB) \cdot p) \\
& = (T' \cdot p, \ f(q' + n'q_0) + (J' + n'B) \cdot p) \\
& \sim ((T' + n' \Delta_l) \cdot p, \ f(q') + J' \cdot p)  \\
& = (I' \cdot p, \ f(q') + J' \cdot p).
\end{align*}
\QED

\begin{Remark}
By the above properties, $k \otimes_{k[Q]} k[P]$ is canonically isomorphic to
\[
k[X_1, \ldots, X_r]/(X_1 \cdots X_l).
\]
\end{Remark}

The converse of the above remark holds under a kind of assumptions of $P$.

\begin{Proposition}
\label{prop:nonsplit:monoids:semistable}
Let $k$ be a field and $f : Q \to P$ an integral homomorphism 
of fine and sharp monoids with $Q \not= \{ 0\}$.
Let $R$ be the completion of $k \otimes_{k[Q]} k[P]$ 
\rom{(}$k$ is a $k[Q]$-module via the canonical homomorphism
$Q \to \{ 0 \}$\rom{)} at the origin and $m$ the maximal ideal of $R$.
We assume the following:
\begin{enumerate}
\renewcommand{\labelenumi}{(\arabic{enumi})}
\item
$f : Q \to P$ does not split, i.e.,
there is no submonoid $N$ of $P$ with
$P = f(Q) \times N$.

\item
Let $R' = R\lformal T_1, \ldots, T_e\rformal $ be the ring of formal power series
over $R$ and $m'$ the maximal ideal of $R'$. Then,
$R'$ is reduced, 
$\dim_k m'/{m'}^2 = \dim R' + 1$ and
$\dim R'/K' = \dim R'$ for all minimal primes $K'$ of $R'$.
\end{enumerate}
Let $p_1, \ldots, p_r$ be all irreducible elements of $P$ which is not lying in $f(Q)$.
Let $l$ be the number of minimal primes of $R$.
Then, after renumbering $p_1, \ldots, p_r$,
$P$ is of semi-stable type
\[
 (r,l, p_1, \ldots, p_r, q_0, b_{l+1}, \ldots, b_{r})
\]
over $Q$ for some $q_0 \in Q \setminus \{ 0\}$ and 
$b_{l+1}, \ldots, b_l \in \NN$.
\end{Proposition}

\Proof
Let us consider a natural homomorphism
\[
H : Q \times \NN^r \to P
\]
given by $H(q, T) = f(q) + T \cdot p$.
Since $f : Q \to P$ is integral, the system of 
congruence relations of $H$ is generated by
\[
\left\{ I_{\lambda} \cdot p = f(q_{\lambda}) + J_{\lambda} \cdot p
\right\}_{\lambda \in \Lambda},
\]
where for each $\lambda \in \Lambda$, $q_{\lambda} \in Q$ and
$I_{\lambda}, J_{\lambda} \in \NN^r$ with $\Supp(I_{\lambda}) \cap \Supp(J_{\lambda})
= \emptyset$.
Let $\phi : k\lformal  X_1, \ldots, X_r\rformal  \to R$ be the homomorphism arising from
\[
k[\NN^r] = k \otimes_{k[Q]} k[Q \times \NN^r] \to k \times_{k[Q]} k[P].
\]
Then, the kernel of $\phi$ is generated by
\[
\left\{ X^{I_{\lambda}} - \beta(q_{\lambda}) X^{J_{\lambda}}
\right\}_{\lambda \in \Lambda},
\]
where $\beta$ is given by
\[
\beta(q) = \begin{cases}
1 & \text{if $\beta = 0$} \\
0 & \text{if $\beta \not= 0$}.
\end{cases}
\]
Let $m$ be the maximal ideal of $R$.
Then, it is easy to see that
$R$ is reduced, 
$\dim_k m/{m}^2 = \dim R + 1$ and
$\dim R/K = \dim R$ for all minimal primes $K$ of $R$.
Let $M$  be the maximal ideal of $k\lformal X_1, \ldots, X_r\rformal $. 
Here $p_i$'s are irreducible. Thus, $\deg(I_{\lambda}) \geq 2$ if $q_{\lambda} \not = 0$,
and $\deg(I_{\lambda}) \geq 2$ and $\deg(J_{\lambda}) \geq 2$
if $q_{\lambda} = 0$. Hence, $\Ker(\phi) \subseteq M^2$.
Therefore,
\[
\dim_k m/m^2 = \dim_k M/(M^2 + \Ker(\phi)) = \dim_k M/M^2 = r,
\]
which says us that $r = \dim R + 1$. Since $R$ is reduced, $\Ker(\phi) = \sqrt{\Ker(\phi)}$.
Thus, we have a decomposition
\[
\Ker(\phi) = K_1 \cap \cdots \cap K_l
\]
such that $K_i$'s are prime, $K_i \not\subseteq K_j$ for all $i \not= j$ and
each $K_i$ corresponds to a minimal prime of $R$.
Note that $\dim k\lformal  X_1, \ldots, X_r \rformal /K_i = r -1$ for each $i$.
Here $k\lformal  X_1, \ldots, X_r \rformal $ is a UFD. Thus, each $K_i$'s are generated by
an irreducible element, so that we can see that there is $f \in k\lformal X_1, \ldots, X_r\rformal $
with $\Ker(\phi) = (f)$.
Here we claim the following:

\begin{Claim}
There is $\lambda \in \Lambda$ with $q_{\lambda} \not= 0$.
\end{Claim}

We assume the contrary.  Let $N$ be a submonoid of $P$ generated by $p_i$'s.
Let us see that
\[
f(q) + n = f(q') + n' \quad (q, q' \in Q, n, n' \in N) \quad\Longrightarrow
\quad
q = q, \ n= n'.
\]
Since $f : Q \to P$ is integral, there are $q_1, q_2 \in Q$ and $n'' \in N$
such that $n = f(q_1) + n''$, $n' = f(q_2) + n''$ and $q + q_1 = q' + q_2$.
Here $q_{\lambda} = 0$ for all $\lambda \in \Lambda$. We can see
$q_1 = q_2 = 0$. Thus, $n = n' = n''$ and $q = q'$.
This observation shows us that $P = Q \times N$,
which contradicts to our assumption.

\medskip
By the above claim, $\Ker(\phi)$ contains an element of the form 
$X^{I_{\lambda}}$.
Note that $f$ is a factor of $X^{I_{\lambda}}$,
$R$ is reduced and $R$ contains $l$ minimal primes.
Thus, after renumbering $p_1, \ldots, p_r$,
we can set $f = X_1 \cdots X_l = X^{\Delta_l}$. Next we claim the following:

\begin{Claim}
$q_{\lambda} \not= 0$ for all $\lambda \in \Lambda$..
\end{Claim}

We assume that there is $\lambda \in \Lambda$ with $q_{\lambda} = 0$.
Then, $X_1 \cdots X_l$ divides
$X^{I_{\lambda}} - X^{J_{\lambda}}$.
This is impossible because $\Supp(I_{\lambda}) \cap \Supp(J_{\lambda}) =
\emptyset$.

\medskip
By the above claim,  we can see that
$N$ is isomorphic to $\NN^r$.
Moreover, $\Ker(\phi)$ is generated by
$\left\{ X^{I_\lambda} \right\}_{\lambda \in \Lambda}$.
Thus, there is $\lambda \in \Lambda$ with
$I_{\lambda} = \Delta_l$.
Hence, we have a congruence relation
$\Delta_l \cdot p  = f(q_0) + B \cdot p$.

Finally, let us consider a relation
\[
I \cdot p  = f(q) + J \cdot p
\]
with $q \not= 0$.
Then, $X^I$ is an element of $\Ker(\phi)$.
Thus, $I(i) > 0$ for all $i=1, \ldots, l$.
\QED

\section{Splitting properties of monoids over a semistable variety}
\label{sec:splitting:monoid:semistable:variety}

In this section, we consider splitting properties of monoids
over a semistable variety.
First, let us consider a log smooth monoid on a smooth variety. 

\begin{Proposition}
\label{prop:split:monoids}
Let $k$ be a field and $f : Q \to P$ an integral homomorphism 
of fine and sharp monoids \rom{(}note that $Q$ might be $\{0\}$\rom{)}.
Let $R$ be the completion of $k \otimes_{k[Q]} k[P]$
\rom{(}$k$ is a $k[Q]$-module via the canonical homomorphism
$Q \to \{ 0 \}$\rom{)} at the origin
and $R\lformal T_1, \ldots, T_e\rformal $ the ring of formal power series
over $R$.
If $R\lformal T_1, \ldots, T_e\rformal $ is regular,
then there are a nonnegative integer $r$ and
a homomorphism $g : \NN^r \to P$ such that
the homomorphism
\[
h : Q \times \NN^r \to P
\]
given by
$h(q, x) = f(q) + g(x)$ is bijective.
\end{Proposition}

\Proof
First of all, note that $R$ is regular.
Let $p_1, \ldots, p_r$ be all irreducible elements of $P$ which 
are not lying in $f(Q)$.
Then, we have a homomorphism
$g : \NN^r \to P$
given by $g(n_1, \ldots, n_r) = \sum_{i=1}^r n_i p_i$.
Thus, we get $h : Q \times \NN^r \to P$ as in the statement of our proposition.
Clearly, $h$ is surjective. Then, since $f : Q \to P$ is integral,
the congruence relation is generated by a system
\[
\left\{ I_{\lambda} \cdot p = f(q_{\lambda}) + J_{\lambda} \cdot p \right\}_{\lambda \in \Lambda},
\]
where $q_{\lambda} \in Q$ and $I_{\lambda}, J_{\lambda} \in \NN^r$
with $\Supp(I_{\lambda}) \cap \Supp(J_{\lambda}) = \emptyset$
for each $\lambda$.
Then, the kernel $K$ of
\[
k\lformal X_1, \ldots, X_r\rformal   \to R
\]
is generated by
\[
\left\{ X^{I_{\lambda}} - \beta(q_{\lambda}) X^{J_{\lambda}}
\right\}_{\lambda \in \Lambda},
\]
where $\beta$ is given by
\[
\beta(q) = \begin{cases}
1 & \text{if $q = 0$} \\
0 & \text{if $q \not= 0$}.
\end{cases}
\]
Using the fact that $p_i$'s are irreducible,
we can see that $K \subset M^2$, where $M$ is the maximal ideal of
$k\lformal X_1, \ldots, X_r\rformal $.
Let $m$ be the maximal ideal of $R$.
Then,
\[
m/m^2 = M/(M^2 + K) = M/M^2
\]
Thus, $\dim_k m/m^2 = r$. On the other hand, if we have a congruence
relation, then $K \not= \{ 0 \}$. Thus,
$\dim R < r$. Therefore, $K = \{ 0 \}$, which means
that $h$ is injective.
\QED

In order to proceed with our arguments,
let us see elementary facts of the ring
\[
 k\lformal X_1, \ldots, X_n\rformal /(X^{I_0} - X^{J_0}).
\]

\begin{Proposition}
\label{prop:unique:x:times:unit}
Let $k$ be a field and $k\lformal X_1, \ldots, X_n\rformal $ the ring of formal power series
of $n$-variables over $k$.
Let $I_0$ and $J_0$ be elements of $\NN^n$ such that $\Supp(I_0) \cap \Supp(J_0) =\emptyset$
and $I_0 \not = (0, \ldots, 0)$ and $J_0 \not= (0, \ldots, 0)$.
Here let us consider the ring 
\[
R = k\lformal X_1, \ldots, X_n\rformal /(X^{I_0} - X^{J_0}).
\]
The image of $X^I$ in $R$ is denoted by $x^I$.
Then, we have the following:
\begin{enumerate}
\renewcommand{\labelenumi}{(\arabic{enumi})}
\item
The multiplication of $X_i$ in $R$ is injective.

\item
For $I, J \in \NN^n$ and $h \in R$, if
$x^I = x^J h$ and $I \not\geq J$, then either $I \geq I_0$ or $I \geq J_0$.

\item
Let $u$ and $v$ be units of $R$. For $I, J \in \NN^n$,
if $x^I u = x^J v$, then $u = v$ and $x^I = x^J$.

\item
For $I, J \in \NN^{n}$,
let us set  $I = I' + a I_0 + b J_0$ and $J = J' + a' I_0 + b' J_0$ such that
$a, b, a', b' \in \NN$
and that $I' \not\geq I_0$, $I' \not\geq J_0$, $J' \not\geq I_0$ and $J' \not\geq J_0$
If $x^I = x^J$, then $I' = J'$ and $a+ b = a' + b'$.

\item
If $\gcm(I_0)$ and $\gcm(J_0)$ are coprime, then
$X^{I_0} - X^{J_0}$ is irreducible in $k\lformal X_1, \ldots, X_n\rformal $.
%
\end{enumerate}
\end{Proposition}

\Proof
(1) Clearly $X_i$ and $X^{I_0} - X^{J_0}$ are coprime.
We assume that $X_i g = 0$ for some $g \in R$.
Then, there is $h \in k\lformal X_1, \ldots, X_r\rformal $ such that
$X_i g = (X^{I_0} - X^{J_0})h$. Thus, $g$ is divisible by $X^{I_0} - X^{J_0}$, 
which means that $g = 0$ in $R$.

\medskip
(2)
We set $X^I - X^J h = (X^{I_0} - X^{J_0})g$.
Moreover, we set
\[
h = \sum_T a_T X^T
\quad\text{and}\quad
g = \sum_T b_T X^T.
\]
Then, we have
\[
X^I - \sum_T a_T X^{T + J} = \sum_T b_T X^{I_0 + T} - \sum_T b_T X^{J_0 + T}.
\]
Since $I \not\geq J$, the term $X^I$ does not appear in
$\sum_T a_T X^{T + J}$.
Thus, the term $X^I$ must appear in
either $\sum_T b_T X^{I_0 + T}$ or $\sum_T b_T X^{J_0 + T}$.
Thus, we get (2).

\medskip
(3)
We set 
\[
a = \max \{ k \in \NN \mid I - kI_0 \geq (0, \ldots, 0) \}
\quad\text{and}\quad
b = \max \{ k \in \NN \mid  I - k J_0 \geq (0, \ldots, 0) \}.
\]
Moreover, we set $I' = I - aI_0-bJ_0$.
Then, $I' \in \NN^n$,  $I' \not\geq I_0$ and $I' \not\geq J_0$.
In the same way, we can find $a'$ and $b'$ such that
if we set $J' = J - a'I_0 - b'J_0$, then
$J' \in \NN^n$,  $J' \not\geq I_0$ and $J' \not\geq J_0$.
Thus,
\[
x^I = x^{I'} x^{(a+b)I_0}
\quad\text{and}\quad
x^J = x^{J'} x^{(a'+ b')I_0}
\]
because $x^{I_0} = x^{J_0}$.
In order to see $u = v$, we may assume that $a' + b' \geq a + b$.
Then, by (1), we have
\[
x^{I'} = x^{J' + lI_0}(v/u),
\]
where $l = (a'+ b') - (a+ b)$. Thus, by (2), we have $I' \geq J' + l I_0$.
Since $I' \not\geq I_0$, we can see $l = 0$. Hence, $I' \geq J'$.
On the other hand, $x^{J'} = x^{I'}(u/v)$. Thus, by (2), $J' \geq I'$.
Therefore, we get $I' = J'$, so that we can obtain
$u = v$, which implies $x^I = x^J$.

\medskip
(4)
First, $x^I = x^{I'} \cdot x^{(a+b)I_0}$ and $x^J = x^{J'} \cdot x^{(a' + b')I_0}$.
Clearly, we may assume that $a' + b' \geq a + b$.
Thus, $x^{I'} = x^{J' + (a' + b' - a - b)I_0}$. Therefore, by (2),
$I' \geq J' + (a' + b' - a - b)I_0$. Here $I' \not\geq I_0$. Thus, $a + b = a' + b'$, so that
$I' \geq J'$ and $x^{J'} = x^{I'}$. By using (2) again, we have $J' \geq I'$.
Therefore, $I' = J'$.

\medskip
(5) First, we need the following lemma:

\begin{Lemma}
\label{lem:torsion:free:integral}
Let $T$ be a fine and sharp monoid such that $T^{gr}$ is torsion free.
Then, $k[T]$ and the completion $k\lformal T\rformal $ at the origin are integral domains.
\end{Lemma}

\Proof
First of all, it is well known that if $\sigma$ is a finitely generated cone in $\QQ^n$
with $\sigma \cap -\sigma = \{ 0 \}$, then
there is an isomorphism $\phi : \QQ^n \to \QQ^n$ such that
$\phi(\sigma) \subseteq \QQ_{\geq 0}^n$.
Thus, we can find an injective homomorphism $\psi : T^{gr} \to \ZZ^n$
such that $\Coker(\psi)$ is finite and $\psi(T) \subseteq \NN^n$.
Thus, $k[T] \hookrightarrow k[\NN^n] = k[X_1, \ldots, X_n]$ and
$k\lformal T\rformal  \hookrightarrow k\lformal \NN^n\rformal  = k\lformal X_1, \ldots, X_n\rformal $.
\QED

\medskip
Let us go back to the proof of Proposition~\ref{prop:unique:x:times:unit}.
Let $N$ be the monoid arising from monomials of $k[X_1, \ldots, X_n]/(X^{I_0} - X^{J_0})$.
Then, $k[N] = k[X_1, \ldots, X_n]/(X^{I_0} - X^{J_0})$. 
By the above lemma, it is sufficient to show that 
$N^{gr}$ has no torsion.
We assume the contrary, that is, $\left(x^S/x^T\right)^n = 1$ and $x^S/x^T \not= 1$,
where $\Supp(S) \cap \Supp(T) = \emptyset$ and $n > 1$.
Then, $x^{nS} = x^{nT}$. Thus, by (4),
there is $L \in \NN$ and $a, b, a', b' \in \NN$ such that
$nS = L + aI_0 + bJ_0$, $nT = L + a'I _0+ b'J_0$, $L \not\geq I_0$, $L \not\geq J_0$ and
$a + b = a' + b'$. Since $\Supp(S) \cap \Supp(T) = \emptyset$, we have $L = 0$.
Hence either $b = 0, a'=0$ or $a= 0, b' = 0$. Considering $x^T/x^S$, we may
assume that $b= 0$ and $a' = 0$. Therefore, we get $nS = aI_0$ and $nT = aJ_0$.
Here there are integers $t_1, \ldots, t_n, t'_{1}, \ldots, t'_{n}$ such that
\[
t_1 I_0(1) + \cdots + t_n I_0(n) + t'_{1} J_0(1) + \cdots t'_{n} J_0(n) = 1.
\]
Thus,
\[
a = \sum_{i=1}^n t_i aI_0(i) + \sum_{i=1}^n t'_i a J_0(i) = 
n \left(  \sum_{i=1}^n t_i S(i) + \sum_{i=1}^n t'_i T(i)  \right).
\]
Hence $a = n l$ for some $l \in \NN$.
Thus, $S = l I_0$ and $T = l J_0$.
Then, 
\[
x^S/x^T = \left( x^{I_0}/x^{J_0}\right)^l = 1.
\]
This is a contradiction.
\QED

\begin{Corollary}
\label{cor:irreducible:deccmp:XI:XJ}
We assume that $k$ is algebraically closed.
Let $I_0$ and $J_0$ be elements of $\NN^n$ such that
$\deg(I_0) \geq 1$, $\deg(J_0) \geq 1$ and $\Supp(I_0) \cap \Supp(J_0) = \emptyset$.
We set $g = \gcm(g(I_0), g(J_0))$, $I_0 = gI'_0$ and $J_0 = g J'_0$.
Then,
\[
X^{I_0} - X^{J_0} = 
(X^{I'_0} - X^{J'_0}) (X^{I'_0} - \zeta X^{J'_0})  \cdots (X^{I'_0} - \zeta^{g-1}X^{J'_0})
\]
is the irreducible decomposition of $X^{I_0} - X^{J_0}$,
where $\zeta$ is a $g$-th primitive root of the unity.
\end{Corollary}

\Proof
It is sufficient to show that $X^{I'_0} - \zeta^i X^{J'_0}$ is irreducible.
Changing coordinates $X_1, \ldots, X_n$ by
$c_1 X_1, \ldots, c_n X_n$, we can make $X^{I'_0} - X^{J'_0}$ of
$X^{I'_0} - \zeta^i X^{J'_0}$.
Thus, by (5) of Proposition~\ref{prop:unique:x:times:unit}, 
we have our corollary.
\QED

\begin{Corollary}
\label{cor:semiatable:type:XI:XJ}
We assume that $k$ is algebraically closed.
Let $I_0$ and $J_0$ be elements of $\NN^n$ such that
$\deg(I_0) \geq 1$, $\deg(J_0) \geq 1$ and $\Supp(I_0) \cap \Supp(J_0) = \emptyset$.
If 
\[
k\lformal X_1, \ldots, X_n\rformal /(X^{I_0} - X^{J_0})
\]
is isomorphic to 
the ring of the type
$k\lformal T_1, \ldots, T_e\rformal /(T_1 \cdots T_l)$
\rom{(}$l \geq 2$\rom{)},
then $\operatorname{char}(k) \not= 2$ and there are $i, j \in \{ 1, \ldots, n \}$ such that
$i \not= j$ and $X^{I_0} - X^{J_0} = X_i^2 - X_j^2$.
\end{Corollary}

\Proof
We set $g = \gcm(g(I_0), g(J_0))$, $I_0 = gI'_0$ and $J_0 = g J'_0$. 
Then, by the above corollary,
\[
X^{I_0} - X^{J_0} = 
(X^{I'_0} - X^{J'_0}) (X^{I'_0} - \zeta X^{J'_0})  \cdots (X^{I'_0} - \zeta^{g-1}X^{J'_0})
\]
is the irreducible decomposition of $X^{I_0} - X^{J_0}$,
where $\zeta$ is a $g$-th primitive root of the unity.
Since $k\lformal X_1, \ldots, X_n\rformal /(X^{I_0} - X^{J_0})$ is reduced, 
$\operatorname{char}(k)$ does not divide $g$.
Here $k\lformal T_1, \ldots, T_n\rformal /(T_1 \cdots T_l)$ has $l$-minimal primes, so that
$g = l$.
Moreover, since every irreducible component is regular,
either $X^{I'_0}$ or $X^{J'_0}$ is linear.
Clearly, we may assume that $X^{I'_0}$ is linear, namely,
$X^{I'_0} = X_i$ for some $i$.
Let $m$ be the maximal ideal of $k\lformal X_1, \ldots, X_n\rformal /(X^{I_0} - X^{J_0})$.
Let $V$ be a vector subspace of $m/m^2$ generated by
$x_i - x^{J_0}, x_i - \zeta x^{J'_0},  \ldots, x_i - \zeta^{l-1}x^{J'_0}$.
Then, we must have $\dim_k V = l$ because
\[
k\lformal X_1, \ldots, X_n\rformal /(X^{I_0} - X^{J_0})
\simeq k\lformal T_1, \ldots, T_n\rformal /(T_1 \cdots T_l).
\]
If $\deg(J'_0) \geq 2$, then $\dim_k V = 1$. This contradict to the fact $l \geq 2$.
Thus, $\deg(J'_0) = 1$, so that $X^{J'_0} = X_j$ for some $j$.
In this case, $\dim_k V \leq 2$. Therefore, $g = l = 2$.
\QED

\begin{Proposition}
\label{prop:monoid:semistable:ring}
Let $k$ be a field,  $N$ a fine and sharp monoid,  and $k\lformal N\rformal $ the completion of $k[N]$
at the origin. Let $\alpha : N \to k\lformal N\rformal $ be the canonical homomorphism.
Let $p_1, \ldots, p_r$ be all irreducible elements of $N$ and
$h : \NN^r \to N$ the natural homomorphism given by
$h(a_1, \ldots, a_r) = \sum_{i=1}^r a_i p_i$.
Let $\phi : k\lformal X_1, \ldots, X_r\rformal  \to k\lformal N\rformal $ be the homomorphism induced by $h$.
Let $R' = k\lformal N\rformal \lformal X_1, \ldots, X_e\rformal $ be the ring of formal power series over $k\lformal N\rformal $
and $m'$ the maximal ideal of $R'$.
We assume that $R'$ is reduced, $\dim_k m'/{m'}^2 = \dim R' + 1$ and
$\dim R'/K' = \dim R'$ for all minimal primes $K'$ of $R'$.
Then, we have the following.
\begin{enumerate}
\renewcommand{\labelenumi}{(\arabic{enumi})}
\item
The kernel of $\phi$ is generated by an element of the form
$X^{I_0} - X^{J_0}$
such that $I_0, J_0 \in \NN^r$,
$\deg(I_0) \geq 2$, $\deg(J_0)\geq 2$,
$\Supp(I_0)\cap \Supp(J_0) = \emptyset$ and
$\gcm(\gcm(I_0), \gcm(J_0))$ is not divisible by $\operatorname{char}(k)$.

\item
Renumbering of $p_1, \ldots, p_r$, we assume that
\[
\Supp(I_0) \subseteq \{ 1, \ldots, l \}
\quad\text{and}\quad
\Supp(J_0) \subseteq \{ l+1, \ldots, r \}.
\]
Let $U$ \rom{(}resp. $V$\rom{)} be the submonoid of $N$ generated by
$p_1, \ldots, p_l$ \rom{(}resp. $p_{l+1}, \ldots, p_r$\rom{)}.
Then, $U \simeq \NN^l$, $V \simeq \NN^{r-l}$ and
the natural homomorphism
\[
U \times_{(I_0 \cdot p,\  J_0 \cdot p)} V \to N
\]
is bijective \rom{(}cf. Conventions and terminology
\rom{\ref{subsub:pushout}}\rom{)}.
\end{enumerate}
\end{Proposition}

\Proof
(1) Let us consider all relations
\[
\left\{ I_{\lambda} \cdot p = J_{\lambda} \cdot p
\right\}_{\lambda \in \Lambda}
\]
in $N$, where $I_{\lambda}, J_{\lambda} \in \NN^r$ and
$\Supp(I_{\lambda}) \cap \Supp(J_{\lambda}) = \emptyset$
for all $\lambda$.
Then, the kernel of $\phi$ is generated by
\[
\left\{ X^{I_{\lambda}} - X^{J_{\lambda}}  \right\}_{\lambda \in \Lambda}.
\]
Let $m$ be the maximal ideal of $k\lformal N\rformal $.
Then, it is easy  to see that $k\lformal N\rformal $ is reduced, 
$\dim_k m/m^2 = \dim k\lformal N\rformal  + 1$ and
$\dim k\lformal N\rformal /K = \dim k\lformal N\rformal $ for all minimal primes $K$ of $k\lformal N\rformal $.
Let $M$ be the maximal ideal of $k\lformal  X_1, \ldots, X_r\rformal $.
Since $p_i$'s are irreducible, $\deg(I_{\lambda}) \geq 2$ and
$\deg(J_{\lambda}) \geq 2$.
Thus, $\Ker(\phi) \subseteq M^2$. Therefore,
\[
m/m^2 = M/(\Ker(\phi) + M^2) = M/M^2.
\]
Then, in the same way as in the proof of Proposition~\ref{prop:nonsplit:monoids:semistable},
there is $f \in k\lformal X_1, \ldots, X_r\rformal $ with
$\Ker(\phi) = (f)$.
We set $X^{I_{\lambda}} - X^{J_{\lambda}} = f u_{\lambda}$
for all $\lambda \in \Lambda$.
If $u_{\lambda}$ is not a unit for every $\lambda \in \Lambda$,
then $X^{I_{\lambda}} - X^{J_{\lambda}} \in f \cdot M$.
Thus, there is $\lambda \in \Lambda$ such that $u_{\lambda}$ is a unit.
Hence we get (1).

\medskip
(2) By using (4) of Proposition~\ref{prop:unique:x:times:unit},
it is easy to see that $U \simeq \NN^l$ and $V \simeq \NN^{r-l}$.
Let $I, I', J, J' \in \NN^r$ such that
\[
 \Supp(I), \Supp(I') \subseteq \{ 1, \ldots, l\}
\quad\text{and}\quad
\Supp(J), \Supp(J') \subseteq \{ l+1, \ldots, r \}.
\]
It is sufficient to see that if
$I \cdot p + J \cdot p = I' \cdot p + J' \cdot p$, then
$(I \cdot p, \ J \cdot p) \sim (I' \cdot p, \ J' \cdot p)$ in
$U \times_{(I_0 \cdot p,\  J_0 \cdot p)} V$.
We set $I = T + aI_0$, $I' = T' + a'I_0$, $J = S + b J_0$ and 
$J' = S' + b' J_0$
such that $a, a', b, b' \in \NN$ and
$T \not\geq I_0$, $T' \not\geq I_0$, $S \not\geq J_0$ and
$S' \not\geq J_0$.
Then, by (4) of Proposition~\ref{prop:unique:x:times:unit}, 
we can see that $T + S = T' + S'$ and $a + b = a' + b'$.
In particular, $T = T'$ and $S = S'$. Therefore, since $(I_0 \cdot p,\  0) \sim (0, \ J_0  \cdot p)$,
\begin{align*}
(I \cdot p,\ J \cdot p) & = ((T + a I_0) \cdot p, \ (S + b J_0) \cdot p) \sim
(T \cdot p, \ (S + (a+b)J_0) \cdot p) \\
& = (T' \cdot p, \ (S' + (a' + b')J_0) \cdot p) \sim
((T' + a'I_0) \cdot p, \ (S' + b J_0) \cdot p) \\
& = (I' \cdot p, \ J' \cdot p).
\end{align*}
\QED

Gathering all results in \S\ref{sec:monoid:semistable:type} and 
\S\ref{sec:splitting:monoid:semistable:variety},
we have the following local result of a smooth log structure
on a semistable variety.

\begin{Proposition}
\label{prop:split:monoid:semistable:variety}
Let $k$ be an algebraically closed field and 
$M_k$ a fine log structure of $\Spec(k)$.
Let $X$ be semistable varieties over $k$ and
$M_X$ a fine log structures of $X$.
We assume that $(X, M_X)$ is log smooth and integral over 
$(\Spec(k), M_k)$.
For a closed point $x \in X$, 
let $(Q \to M_k, P \to M_{X, \bar{x}}, Q \to P)$ be a good chart of
$(X, M_X) \to (\Spec(k), M_k)$ at $x$, that is,
$Q \to \overline{M}_k$ and $P \to \overline{M}_{X, \bar{x}}$ are
bijective homomorphisms of fine and sharp monoids,
$k \otimes_{k[Q]} k[P] \to \OO_{X, \bar{x}}$ is smooth and
the following diagram
\[
\begin{CD}
Q @>>> P \\
@VVV @VVV \\
M_k @>>> M_{X, \bar{x}}
\end{CD}
\]
is commutative.
Then, we have the following:
\begin{enumerate}
\renewcommand{\labelenumi}{(\arabic{enumi})}
\item
If $\mult_x(X) = 1$, that is, $x$ is a regular point,
then $Q \to P$ splits and $P \simeq Q \times \NN^r$ for some $r$.

\item
If $\mult_x(X) = 2$, then we have one of the following:
\begin{enumerate}
\renewcommand{\labelenumii}{(\arabic{enumi}.\arabic{enumii})}
\item
If $Q \to P$ does not split, then $P$ is of semistable type over $Q$.

\item
If $Q \to P$ splits, then $\operatorname{char}(k) \not= 2$ and
$\widehat{\OO}_{X, x}$ is canonically isomorphic to 
\[
k\lformal X_1, \ldots, X_r\rformal /(X_1^2 - X_2^2).
\]
More precisely, if $p_1, \ldots, p_r$ are all irreducible elements
of $P$ not lying in the image of $Q \to P$,
then $\alpha(p_i) = X_i$ modulo $X_1^2 - X_2^2$ after renumbering
$p_1, \ldots, p_r$, where
$\alpha$ is the composition of $P \to \rest{M_X}{U} \to \OO_U$.
\end{enumerate}

\item
If $\mult_x(X) \geq 3$, then $Q \to P$ does not split and
$P$ is of semistable type over $Q$.

\item
If $\mult_x(X) \geq 2$ and
$P^{gr}$ is torsion free, then $Q \to P$ does not split and
$P$ is of semistable type over $Q$.
\end{enumerate}
In particular, if $M_X$ is saturated, then, for all $x \in X$,
$P$ is a monoid of semistable type over $Q$.
\end{Proposition}

\Proof
If $x \not\in \Sing(X)$, then our assertion holds by Proposition~\ref{prop:split:monoids}.
Thus, we may assume that $x \in \Sing(X)$.

We assume that $Q \to P$ split, so that $P \simeq Q \times N$ for some $N$.
Then, 
\[
k \otimes_{k[Q]} k[P] \simeq k[N].
\]
Since $k[N] \to \OO_X$ is smooth,
$k\lformal N\rformal \lformal X_1, \ldots, X_e\rformal $ is isomorphic to the ring of the type 
$k\lformal T_1, \ldots, T_n\rformal /(T_1 \cdots T_l)$.
Thus, by Corollary~\ref{cor:semiatable:type:XI:XJ} and
Proposition~\ref{prop:monoid:semistable:ring},
$\operatorname{char}(k) \not= 2$  and $l = 2$.
Moreover, if $P^{gr}$ is torsion free, then $N^{gr}$ is torsion free.
Thus, $k\lformal N\rformal $ is an integral domain by Lemma~\ref{lem:torsion:free:integral}.
This is a contradiction.
Therefore,  if $P^{gr}$ is torsion free, then $Q \to P$ does not split.

If $Q \to P$ does not split, then we get our assertion by
Proposition~\ref{prop:nonsplit:monoids:semistable}.
\QED

\section{Log morphisms and log differential sheaves on a semistable variety}
Here we consider a log structure on a semistable variety.
Especially, we consider a uniqueness problem of
a log morphism for the fixed scheme morphism, which
is one of main results of this paper.

\begin{Theorem}
\label{thm:unqueness:log:hom:semistable}
Let $k$ be an algebraically closed field and $M_k$ a fine log structure
of $\Spec(k)$.
Let $X$ and $Y$ be semistable varieties over $k$ and
$M_X$ and $M_Y$ fine log structures of $X$ and $Y$ respectively.
We assume that $(X, M_X)$ and $(Y, M_Y)$ are 
log smooth and integral over $(\Spec(k), M_k)$.
Let $\Supp(M_Y/M_k)$ be the locus of $Y$
such that $y \in \Supp(M_Y/M_k)$ if and only if
the natural homomorphism 
$M_k \times \OO_{Y, \bar{y}}^{\times} \to M_{Y, \bar{y}}$
is not surjective.
Let $\phi : X \to Y$ be a morphism over $k$ such that
$\phi(X') \not\subseteq \Supp(M_Y/M_k)$
for any irreducible component $X'$ of $X$.
If $(\phi, h) : (X, M_X) \to (Y, M_Y)$ and
$(\phi, h') : (X, M_X) \to (Y, M_Y)$ are morphisms of log schemes over
$(\Spec(k), M_k)$, then $h = h'$.
\end{Theorem}

\Proof
This is a local question. 
Let us take a fine and sharp monoid $Q$ with $M_k = Q \times k^{\times}$.
Let $x$ be a closed point of $X$ and $y = f(x)$.
Let us choose etale local neighborhoods $U$ and $V$ at $x$ and $y$
respectively with $f(U) \subseteq V$. Moreover, shrinking $U$ and $V$ enough,
by Corollary~\ref{cor:good:pair:chart}, we may assume that
there are good charts
\[
(Q \to M_k, \pi : P \to \rest{M_X}{U}, f : Q \to P)
\]
and
\[
(Q \to M_k, \pi' : P' \to \rest{M_Y}{V}, f' : Q \to P')
\]
of $(X, M_X) \to (\Spec(k), M_k)$ and
$(Y, M_Y) \to (\Spec(k), M_k)$ at $x$ and $y$ respectively.
Let $\tilde{\pi} : P \times \OO^{\times}_{X, \bar{x}} \to M_{X, \bar{x}}$
and $\tilde{\pi}' : P' \times \OO^{\times}_{Y, \bar{y}} \to M_{Y, \bar{y}}$
be the natural homomorphisms induced by $\pi$ and $\pi'$.
Note that $\tilde{\pi}$ and $\tilde{\pi}'$ are isomorphisms.
Let $H : P' \times \OO^{\times}_{Y, \bar{y}} \to P \times \OO^{\times}_{X, \bar{x}}$
and $H' : P' \times \OO^{\times}_{Y, \bar{y}} \to P \times \OO^{\times}_{X, \bar{x}}$
be homomorphisms of monoids such that the following diagrams are commutative:
\[
\begin{CD}
P' \times \OO^{\times}_{Y, \bar{y}} @>{H}>> P \times \OO^{\times}_{X, \bar{x}} \\
@V{\tilde{\pi}'}VV @VV{\tilde{\pi}}V \\
M_{Y, \bar{y}} @>{h}>> M_{X, \bar{x}} \\
@V{\alpha'}VV @VV{\alpha}V \\
\OO_{Y, \bar{y}} @>{\phi^*}>> \OO_{X, \bar{x}}
\end{CD}
\qquad
\begin{CD}
P' \times \OO^{\times}_{Y, \bar{y}} @>{H'}>> P \times \OO^{\times}_{X, \bar{x}} \\
@V{\tilde{\pi}'}VV @VV{\tilde{\pi}}V \\
M_{Y, \bar{y}} @>{h'}>> M_{X, \bar{x}} \\
@V{\alpha'}VV @VV{\alpha}V \\
\OO_{Y, \bar{y}} @>{\phi^*}>> \OO_{X, \bar{x}}
\end{CD}
\]
Here $\alpha$ and $\alpha'$ are the canonical homomorphism.
By abuse of notation,
$\alpha \cdot \tilde{\pi}$ and $\alpha' \cdot \tilde{\pi}'$ are also denoted by
$\alpha$ and $\alpha'$. Then,
$\alpha(p, u) = \alpha(\pi(p)) \cdot u$ and
$\alpha'(p', u') = \alpha'(\pi'(p')) \cdot u'$.

First we claim the following:

\begin{Claim}
$H(0, u) = H'(0, u)$ for all $u \in \OO^{\times}_{Y, \bar{y}}$.
\end{Claim}

We set $H(0, u) = (f(q) + \sum_{i=1}^r a_i p_i, v)$,
where $p_1, \ldots, p_r$ are all irreducible elements of $P$ not lying in
$f(Q)$.
Let us consider the above commutative diagram.
Then,
\[
\phi^*(u) = \phi^*(\alpha'(0, u)) =\alpha(H(0, u)) = \beta(q) x_1^{a_1} \cdots x_r^{a_r} v,
\]
where $x_i = \alpha(p_i, 1)$ and $\beta$ is given by
\[
\beta(q) = \begin{cases}
1 & \text{if $q = 0$} \\
0 & \text{if $q \not= 0$}.
\end{cases}
\]
Since $\phi^*(u)$ is a unit in $\OO_{X, \bar{x}}$ and
$x_1, \ldots, x_r$ are not units, we have
$q = 0$ and $a_1 = \cdots = a_r = 0$. Thus, $v = \phi^*(u)$.
Hence $H(0, u) = (0, \phi^*(u))$. In the same way, we can see $H'(0, u) = (0, \phi^*(u))$.
Therefore, $H(0, u) = H'(0, u)$.

\medskip
Next we claim

\begin{Claim}
$H(f'(q), 1) = H'(f'(q), 1)$ for all $q \in Q$.
\end{Claim}

Let us consider homomorphisms
\[
\tilde{f} : Q \to M_{X, \bar{x}} \to P \times \OO^{\times}_{X, \bar{x}}
\quad\text{and}\quad
\tilde{f}' : Q \to M_{Y, \bar{y}} \to P' \times \OO^{\times}_{Y, \bar{y}}.
\]
Then, we can set $\tilde{f}(q) = (f(q), \gamma(q))$ and
$\tilde{f}'(q) = (f'(q), \gamma'(q))$.
Here, $h$ and $h'$ are homomorphisms over $M_k$. Thus the following diagrams
are commutative.
\[
\begin{CD}
P' \times \OO^{\times}_{Y, \bar{y}} @>{H}>> P \times \OO^{\times}_{X, \bar{x}} \\
@V{\tilde{\pi}'}VV @VV{\tilde{\pi}}V \\
M_{Y, \bar{y}} @>{h}>> M_{X, \bar{x}} \\
@AAA @AAA \\
Q @= Q .
\end{CD}\qquad
\begin{CD}
P' \times \OO^{\times}_{Y, \bar{y}} @>{H'}>> P \times \OO^{\times}_{X, \bar{x}} \\
@V{\tilde{\pi}'}VV @VV{\tilde{\pi}}V \\
M_{Y, \bar{y}} @>{h'}>> M_{X, \bar{x}} \\
@AAA @AAA \\
Q @= Q .
\end{CD}
\]
Hence, we can see
\[
H(f'(q), \gamma'(q)) = H'(f'(q), \gamma'(q)) = (f(q), \gamma(q)).
\]
Thus,
\begin{align*}
H(f'(q), 1) & = H((f'(q), \gamma'(q)) + (0, \gamma'(q)^{-1})) =
(f(q), \gamma(q)) + (0, \phi^*( \gamma'(q))^{-1}) \\
& = (f(q), \gamma(q) \cdot \phi^*( \gamma'(q))^{-1}).
\end{align*}
In the same way, we have $H'(f'(q), 1) = (f(q), \gamma(q) \cdot \phi^*( \gamma'(q))^{-1})$.
Thus, we get our claim.

\medskip
From now on, we consider the following four cases:
\begin{enumerate}
\renewcommand{\labelenumi}{(\Alph{enumi})}
\item
$f : Q \to P$ splits and $f' : Q \to P'$ splits.

\item
$f : Q \to P$ does not split and $f' : Q \to P'$ splits.

\item
$f : Q \to P$ splits and $f' : Q \to P'$ does not split.

\item
$f : Q \to P$ does not split and $f' : Q \to P'$ does not split.
\end{enumerate}
For each case, let $U_1, \cdots, U_l$ and $V_1, \cdots, V_{l'}$ be
all irreducible components of $U$ and $V$ respectively.
Here since $\Sing(Y) \subseteq \Supp(M_Y/M_k)$ and
$\phi(U_j) \not\subseteq \Supp(M_Y/M_k)$,
for each $j$, there is a unique $i$ with $\phi(U_j) \subseteq V_i$.
We denote this $i$ by $\sigma(j)$.
Note that we have a map
$\sigma : \{ 1, \ldots, l \} \to \{ 1, \ldots, l' \}$.
In the following, we give $p_1, \ldots, p_r \in P$ (resp. $p'_1, \ldots, p'_{r'} \in P'$) for each case
(A), (B), (C) and (D) such that $P$ (resp. $P'$) is  generated by
$f(Q)$ and $p_1, \ldots, p_r$ (resp. $f'(Q')$ and $p'_1, \ldots, p'_{r'}$). 
The last claim is the following:

\begin{Claim}
$H(p'_i, 1) = H'(p'_i, 1)$ for all $i= 1, \cdots, r'$.
\end{Claim}

For this purpose, we fix common notation for all cases.
We denote $\alpha(p_j, 1)$ by $x_j$ and $\alpha'(p'_i, 1)$ by $y_i$.
Here we set
\addtocounter{Claim}{1}
\begin{equation}
H(p'_i, 1) = \left(f(q_i) + I_i \cdot p, \ u_i\right)
\quad\text{and}\quad
H'(p'_i, 1) = \left(f(q'_i) + I'_i \cdot p, \ u'_i\right),
\end{equation}
where $I_i, I'_i \in \NN^r$, $q_i, q'_i \in Q$ and $u_i, u'_i \in \OO^{\times}_{X, \bar{x}}$.
Then, since $\alpha(H(p'_i, 1)) = \phi^*(\alpha'(p'_i, 1))$ and
$\alpha(H'(p'_i, 1)) = \phi^*(\alpha'(p'_i, 1))$,
we have
\addtocounter{Claim}{1}
\begin{equation}
\phi^*(y_i) = \beta(q_i) \cdot x^{I_i} \cdot u_i = \beta(q'_i) \cdot x^{I'_i} \cdot u'_i.
\end{equation}
Let us begin with Case A.

\medskip
(Case A): In this case, there are submonoids $N$ and $N'$ of $P$ and $P'$ respectively
such that $P = f(Q) \times N$ and $P' = f'(Q) \times N'$.
Let $p_1, \ldots, p_r$ (resp. $p'_1, \ldots, p'_{r'}$) be 
all irreducible elements of $N$ (resp. $N'$).
By Proposition~\ref{prop:split:monoid:semistable:variety},
\[
 \Supp(M_Y/M_k) = \{ y_1 = 0\} \cup \cdots \cup \{ y_{r'} = 0 \}.
\]
around $\bar{y}$. Thus,
\[
\rest{\phi^*(y_i)}{U_j} =
\rest{\beta(q_i) \cdot x^{I_i} \cdot u_i }{U_j} = 
\rest{\beta(q'_i) \cdot x^{I'_i} \cdot u'_i}{U_j} \not= 0
\]
for all $j$.  In particular, $q_i = q'_i = 0$ for all $i=1, \ldots, r'$.
Therefore,
\[
x^{I_i} \cdot u_i = x^{I'_i} \cdot u'_i
\]
for all $i$. Thus, by (3) of Proposition~\ref{prop:unique:x:times:unit},
$u_i = u'_i$ and
$x^{I_i} = x^{I'_i}$.
Note that the natural homomorphism $k[N] \to \OO_{X, \bar{x}}$
is injective. Thus, we get $I_i \cdot p = I'_i \cdot p$.

\medskip
(Case B): In this case,
there is a submonoid $N'$ of $P'$
such that $P' = f'(Q) \times N'$.
Let $p'_1, \ldots, p'_{r'}$ be 
all irreducible elements of $N'$. Moreover, by 
Proposition~\ref{prop:nonsplit:monoids:semistable},
$P$  is of semistable type
\[
(r,l, p_1, \ldots, p_r, q_0, b_{l+1}, \ldots, b_{r})
\]
over $Q$.
Renumbering $U_1, \ldots, U_l$, we may assume that $U_j$ is defined by
$x_j = 0$. 
By Proposition~\ref{prop:split:monoid:semistable:variety},
\[
 \Supp(M_Y/M_k) = \{ y_1 = 0\} \cup \cdots \cup \{ y_{r'} = 0 \}.
\]
around $\bar{y}$. Thus
\[
\rest{\phi^*(y_i)}{U_j}=
\rest{\beta(q_i) \cdot x^{I_i} \cdot u_i}{U_j} = 
\rest{\beta(q'_i) \cdot x^{I'_i} \cdot u'_i}{U_j} \not= 0
\]
for all $j$.  In particular, $q_i = q'_i = 0$ and
$I_i(j) = I'_i(j) = 0$ for $j=1, \ldots, l$.
Further since $\OO_{U_j, \bar{x}}$ is a UFD,
we can see that $I_i = I'_i$.
Moreover,
$\rest{u_i}{U_j} = \rest{u'_i}{U_j}$ for all $j$.
Thus, $u_i = u'_i$. Therefore, $H(p'_i, 1) = H'(p'_i, 1)$ for all $i = 1, \ldots, r'$.

\medskip
(Case C): There is a submonoid $N$ of $P$
such that $P = f(Q) \times N$.
Let $p_1, \ldots, p_{r}$ be 
all irreducible elements of $N$. Moreover, by 
Proposition~\ref{prop:nonsplit:monoids:semistable},
$P'$  is of semistable type
\[
(r',l', p'_1, \ldots, p'_{r'}, q'_0, b'_{l+1}, \ldots, b'_{r})
\]
over $Q$.
Renumbering $V_1, \ldots, V_{l'}$, we may assume that $V_i$ is defined by
$y_i = 0$. 
Note that
\[
 \Supp(M_Y/M_k) = \Sing(Y) \cup \{ y_{l'+1} = 0 \}\cup \cdots 
\cup \{ y_{r'} = 0 \}
\]
around $\bar{y}$. Therefore,
if $i \not= \sigma(j)$, then $\rest{\phi^*(y_i)}{U_j} \not = 0$.
Thus, we can see $q_i = q'_i = 0$ for $i \not= \sigma(j)$.

First, we consider the case where
$\sigma(1) = \cdots = \sigma(r) = s$.
Note that $s \leq l'$.
Then, for $i \not=  s$,
$q_i = q'_i = 0$.
Thus, 
$x^{I_i} \cdot u_i = x^{I'_i} \cdot u'_i$
for all $i \not= s$. Therefore, in the same way as in Case A,
we can see
\[
I_i \cdot p = I'_i \cdot p
\quad\text{and}\quad
u_i = u'_i
\]
for all $i \not= s$.
On the other hand, we have the relation
$p'_1 + \cdots + p'_{l'} = f'(q'_0) + \sum_{i > l'} b'_i p'_i$.
Therefore, we have $H(p'_s, 1) = H'(p'_s, 1)$.

Hence, we may assume that $\#(\sigma(\{ 1, \cdots, l  \})) \geq 2$.
In this case, we can conclude that 
$q_i = q'_i = 0$ for all $i$.
Therefore, in the same way as in Case A,
we can see
\[
I_i \cdot p = I'_i \cdot p
\quad\text{and}\quad
u_i = u'_i
\]
for all $i$.

\medskip
(Case D):
By Proposition~\ref{prop:nonsplit:monoids:semistable},
$P$ and $P'$ are of semistable type
\[
(r,l, p_1, \ldots, p_r, q_0, b_{l+1}, \ldots, b_{r})
\quad\text{and}\quad
(r',l', p'_1, \ldots, p'_{r'}, q'_0, b'_{l'+1}, \ldots, b'_{r'})
\]
over $Q$.
Renumbering $U_1, \ldots, U_l$ and $V_1, \ldots, V_{l'}$, we may
assume that $U_j$ is defined by $x_i = 0$ and $V_i$ is defined by
$y_i = 0$.
Note that
\[
 \Supp(M_Y/M_k) = \Sing(Y) \cup \{ y_{l'+1} = 0 \}\cup \cdots 
\cup \{ y_{r'} = 0 \}
\]
around $\bar{y}$. Therefore,
if $i \not= \sigma(j)$, then $\rest{\phi^*(y_i)}{U_j} \not = 0$.
Thus, we can see $q_i = q'_i = 0$ and $I_i(j) = I'_i(j) = 0$.
Moreover, since $\OO_{U_j, \bar{x}}$ is a UFD, considering
$\rest{\phi^*(y_i)}{U_j}$, we can see that
\[
I_i = I'_i
\quad\text{and}\quad
\rest{u_i}{U_j} = \rest{u'_i}{U_j}.
\]
Gathering the above observations, we get the following:
For all $i = 1, \cdots, r'$ and $j = 1, \ldots, l$ with $i \not= \sigma(j)$,
\addtocounter{Claim}{1}
\begin{equation}
\begin{cases}
q_i = q'_i = 0, \\
I_i(j) = I'_i(j) = 0, \\
I_i = I'_i, \\
\rest{u_i}{U_j} = \rest{u'_i}{U_j}.
\end{cases}
\end{equation}

Let us see that for all $i > l'$,
\[
q_i = q'_i = 0, \ u_i = u'_i, \ I_i = I'_i.
\]
Note that if $i > l'$, then $i \not= \sigma(j)$ for all $j = 1, \ldots, l$.
Thus, we get $q_i = q'_i = 0$ and
$I_i = I'_i$. Moreover,
$\rest{u_i}{U_j} = \rest{u'_i}{U_j}$ for all $j = 1, \ldots, l$. Thus,
$u_i = u'_i$. Therefore,
\addtocounter{Claim}{1}
\begin{equation}
H(p'_i, 1) = H'(p'_i, 1) \quad\text{for all $i > l'$}.
\end{equation}

First, we consider the case where
$\sigma(1) = \cdots = \sigma(r) = s$.
Then, for $i \not=  s$,
\[
q_i = q'_i = 0, \ I_i = I'_i.
\]
Moreover, for all $j = 1. \ldots, l$ and $i \not = s$, $\rest{u_i}{U_j} = \rest{u'_i}{U_j}$.
Therefore, $u_i = u'_i$ for $i \not= s$.
Thus, $H(p'_i, 1) = H'(p'_i, 1)$ for all $i \not= s$.
On the other hand, we have the relation
$p'_1 + \cdots + p'_{l'} = f'(q'_0) + \sum_{i > l'} b'_i p'_i$.
Therefore, we have $H(p'_s, 1) = H'(p'_s, 1)$.

Hence, we may assume that $\#(\sigma(\{ 1, \cdots, l  \})) \geq 2$.
In this case, we can conclude that 
\[
q_i = q'_i = 0, \ I_i = I'_i
\]
for all $i$.
Moreover, $\rest{u_i}{U_j} = \rest{u'_i}{U_j}$ if $i \not= \sigma(j)$.
Since $p'_1 + \cdots + p'_{l'} = f'(q'_0) + \sum_{i > l'} b'_i p'_i$,
\[
H(p'_1 + \cdots + p'_{l'}, 1) = H'(p'_1 + \cdots + p'_{l'}, 1).
\] 
Thus, considering the
$\OO^{\times}_{X, \bar{x}}$-factor, we find
\[
u_1 \cdots u_{l'} = u'_1 \cdots u'_{l'}.
\]
Moreover, if we set $S_i =  \{ 1, \ldots, l \} \setminus \sigma^{-1}(i)$, then
$S_i \cup S_{i'} = \{ 1, \ldots, l\}$
for all $i \not= i'$.
Further, if we set $v_i = u_i/u'_i$, then
\[
v_1 \cdots v_{l'} = 1
\quad\text{and}\quad
\rest{v_i}{U_j} = 1
\quad\text{for all $j \in S_i$ and all $i=1, \ldots, l'$}.
\]
Therefore, using the following Lemma~\ref{lem:unit:lemma}, 
we have $v_i = 1$ for all $i= 1, \ldots, l'$.
Hence, we can see $H(p'_i, 1) = H'(p'_i, 1)$ for $i = 1, \ldots, l'$.
\QED

\begin{Lemma}
\label{lem:unit:lemma}
Let $k$ be a fields, $R = k\lformal  X_1, \ldots, X_n\rformal /(X_1 \cdots X_l)$ and
$\Lambda = \{ 1, \ldots, l\}$.
Let $\pi_j : R \to R/X_iR$ be the canonical homomorphism for
$j \in \Lambda$.
Let $S_1, \ldots, S_s$ be subsets of $\Lambda$ with
$S_{i} \cup S_{i'} = \Lambda$ for $i \not= i'$.
Moreover, let $u_1, \ldots, u_s$ be units in $R$.
If $u_1 \cdots u_s = 1$ and, for each $i$,
$\pi_j(u_i) = 1$ for all $j \in S_i$, 
then $u_1 = \cdots = u_s = 1$.
\end{Lemma}

\Proof
If $S_{i_0} = \emptyset$ for some $i_0$, then $S_i = \Lambda$ for all $i \not= i_0$.
Thus, $u_i = 1$ for all $i \not= i_0$ because
\[
\pi_1 \times \cdots \times \pi_l : R \to R/X_1R \times \cdots \times R/X_lR
\]
is injective. Then, $u_{i_0} = 1$.
Therefore, we may assume that $S_i \not= \emptyset$ for all $i$.

For a monomial $X_1^{a_1} \cdots X_n^{a_n}$, the support with respect to $\Lambda$
is given by
\[
\Supp_{\Lambda}(X_1^{a_1} \cdots X_n^{a_n}) = \{ i \in \Lambda \mid a_i > 0 \}.
\]
For a subset $S$ of $\Lambda$, let
$\Gamma_S$ be the set of formal sums of monomials $X_1^{a_1} \cdots X_n^{a_n}$
with $\Supp_{\Lambda}(X_1^{a_1} \cdots X_n^{a_n})  = S$. 
Note that $\Gamma_{\emptyset} = k\lformal X_{l+1}, \ldots, X_n\rformal $. Then,
\[
k\lformal  X_1, \ldots, X_n \rformal  = \bigoplus_{S \subseteq \Lambda} \Gamma_S.
\]
Moreover, the natural map $\bigoplus_{S \subsetneq \Lambda} \Gamma_S \to R$
is an isomorphism as $k$-vector spaces. We denote the image of $\Gamma_S$ in
$R$ by $\overline{\Gamma}_S$. For $f_S \in \overline{\Gamma}_S$ and
$f_{S'} \in \overline{\Gamma}_{S'}$, 
$f_S \cdot f_{S'} \in \overline{\Gamma}_{S \cup S'}$ if $S \cup S' \subsetneq \Lambda$, and
$f_S \cdot f_{S'} = 0$ if $S \cup S' = \Lambda$.

Here we set $u_i = \sum_{S \subsetneq \Lambda} f_{i, S}$,
where $f_{i, S} \in \overline{\Gamma}_{S}$. Then, for all $j \in S_i$,
\[
\pi_j(u_i) = \sum_{j \not\in S \subsetneq \Lambda} f_{i, S} = 1.
\]
Thus, $f_{i, \emptyset} = 1$ and $f_{i, S} = 0$ for all $S \not= \emptyset$ with
$j \not\in S$. Therefore, if we set
\[
\Delta_i = \{ S \subsetneq \Lambda \mid S_i \subseteq S \},
\]
we can write
\[
u_i = 1 + \sum_{S \in \Delta_i} f_{i, S}.
\]
Since $S_i \cup S_{i'} = \Lambda$ ($i \not= i'$), for $S \in \Delta_i$ and 
$S' \in \Delta_{i'}$ with $i \not= i'$, we can easily see (1) $S \cup S'  = \Lambda$ and
(2) $S \not= S'$.
Thus, using the above (1), we obtain
\[
u_1 \cdots u_s = 1 + \sum_{i=1}^s \sum_{S \in \Delta_i}f_{i,S}.
\]
Moreover, using the above (2), we can find $f_{i, S} = 0$.
Thus, we get $u_i = 1$ for all $i$.
\QED

\begin{Remark}
If we do not assume the condition
\[
 \text{``$\phi(X') \not\subseteq \Supp(M_Y/M_k)$
for any irreducible component $X'$ of $X$''}
\]
in Theorem~\ref{thm:unqueness:log:hom:semistable},
then the assertion of the theorem does not hold in general.
For example, let us consider $\AAA^1_k = \Spec(k[X])$.
Let $M$ be a log structure associated with
$\alpha : \NN \times \NN \to k[X]$ given by
\[
\alpha(a, b) = \begin{cases}
X^b & \text{if $a = 0$} \\
0 & \text{if $a \not= 0$}.
\end{cases}
\]
Further, let $f : \NN \to \NN \times \NN$ be a homomorphism
defined by $f(a) = (a, 0)$. Then, $(\AAA^1_k, M)$ is log smooth and
integral over $(\Spec(k), \NN \times k^{\times})$.
Here let us consider a morphism $\phi : \AAA^1_k \to \AAA^1_k$
induced by a homomorphism $\psi : k[X] \to k[X]$ given by
$\psi(X) = 0$. Then, $\phi(\AAA^1_k) =\Supp(M/\NN \times k^{\times})$.
Moreover, we consider a homomorphism
\[
h : \NN \times \NN \to \NN \times \NN
\]
defined by $h(1, 0) = (1, 0)$ and $h(0,1) = (a_0, b_0)$ ($a_0 > 0$).
Then, it is easy to see that the following diagrams are commutative:
\[
\begin{CD}
\NN \times \NN @>{h}>> \NN \times \NN \\
@A{f}AA @AA{f}A \\
\NN @= \NN
\end{CD}\qquad\qquad
\begin{CD}
\NN \times \NN @>{h}>> \NN \times \NN \\
@V{\alpha}VV @VV{\alpha}V \\
k[X] @>{\psi}>> k[X]
\end{CD}
\]
Thus, $(\phi, h) : (\AAA^1_k, M) \to (\AAA^1_k, M)$ is a log morphism over
$(\Spec(k), \NN)$. On the other hand,
we have infinitely many choices of $a_0$ and $b_0$.
\end{Remark}

\medskip
Finally, let us consider a log differential module on 
a semistable variety.

\begin{Proposition}
\label{prop:log:differential:smooth:comp}
Let $k$ be an algebraically closed field and $M_k$ a fine log structure
of $\Spec(k)$.
Let $X$ be a semistable variety over $k$ and
$M_X$ a fine log structure of $X$.
We assume that $(X, M_X)$ is log smooth and integral over $(\Spec(k), M_k)$.
Let $\nu : \widetilde{X} \to X$ be the normalization of $X$ and
$M_{\widetilde{X}}$ the underlining log structure of $\nu^*(M_X)$, 
that is, $M_{\widetilde{X}} = \nu^*(M_X)^u$
\rom{(}cf. see Conventions and terminology \rom{\ref{subsubsec:underlining:log:structure}}\rom{)}.
Then, $(\widetilde{X}, M_{\widetilde{X}})$ 
is log smooth over $(\Spec(k), k^{\times})$ and
$\Omega^1_{\widetilde{X}}(\log(M_{\widetilde{X}}/k^{\times}))$
is isomorphic to
$\nu^* \Omega^1_X(\log(M_X/M_k))$.
\end{Proposition}

\Proof
First of all, there is a fine and sharp monoid $Q$ with
$M_k = Q \times k^{\times}$.
Let $\alpha : M_{X} \to \OO_X$ and 
$\alpha' : \nu^*(M_X) \to \OO_{\widetilde{X}}$ be the canonical homomorphisms.
For a closed point $x \in \widetilde{X}$,
let $(\pi_Q : Q \to M_k, \pi_P : P \to M_{X, \overline{\nu(x)}}, 
f : Q \to P)$
be a good chart of $(X, M_X) \to (\Spec(k), M_k)$ at $\nu(x)$.
Here we consider three cases:
\begin{enumerate}
\renewcommand{\labelenumi}{(\Alph{enumi})}
\item
$\nu(x)$ is a smooth point of $X$.

\item
$\nu(x)$ is a singular point of $X$ and $f : Q \to P$ splits.

\item
$\nu(x)$ is a singular point of $X$ and $f : Q \to P$ does not split.
\end{enumerate}

\begin{Claim}
\label{prop:log:differential:smooth:comp:claim:1}
$(\widetilde{X}, M_{\widetilde{X}}) \to (\Spec(k), k^{\times})$
is log smooth at $x$.
\end{Claim}

(Case A): In this case, $\nu(x) = x$.
Then, by Proposition~\ref{prop:split:monoid:semistable:variety},
$P = f(Q) \times \NN^r$.
Let $e_i$ be the $i$-th standard basis of $\NN^r$ and
$T_i = 1 \otimes e_i$ in $k \otimes_{k[Q]} k[P]$.
Then, $k[T_1, \ldots, T_r]_{(T_1, \ldots, T_r)} \to \OO_{X, \bar{x}}$
is smooth. Therefore, adding indeterminates $T_{r+1}, \ldots, T_n$,
we have 
\[
 h : k[T_1, \ldots, T_r, T_{r+1}, \ldots,T_n]_{(T_1, \ldots, T_n)} \to 
     \OO_{X,\bar{x}}
\]
is etale.
We set $t_i = \alpha(\pi_P(e_i))$ for $i = 1, \ldots, r$.
Then, $t_1, \ldots, t_r$ form a part of local parameters of
$\OO_{X, \bar{x}}$ because $h(T_i) = t_i$ for $i=1, \ldots, r$ and
$h$ is etale.
Moreover, $M_{\widetilde{X}, \bar{x}}$ is generated by 
$t_1, \ldots, t_r$
and $\OO^{\times}_{X,\bar{x}}$.
Thus, we get our assertion.

\medskip
(Case B): In this case, 
by Proposition~\ref{prop:split:monoid:semistable:variety}, 
$\operatorname{char}(k) \not = 2$,
$P = f(Q) \times N$ and $N$ is a monoid such that
\[
 k[N] = k[T_1, \ldots, T_r]/(T_1^2 - T_2^2).
\]
Moreover, adding indeterminates $T_{r+1}, \ldots, T_{n+1}$,
\[
 h : k[T_1, \ldots, T_r, T_{r+1}, \ldots, T_{n+1}]_{(T_1, \ldots, T_{n+1})}
 /(T_1^2 - T_2^2) \to \OO_{X, \overline{\nu(x)}}
\]
is etale.
We set $t_i = \alpha(\pi_P(\bar{T}_i))$ for $i=1, \ldots, r$.
Changing the sign of $\pi_P(\bar{T}_2)$, we may assume that
$\widetilde{X}$ at $x$ is the component corresponding to
$t_1 = t_2$. Note that
$h(\bar{T}_i) = t_i$ for $i=1, \ldots, r$.
Thus,
$M_{\widetilde{X}, \bar{x}}$ is generated by $t_2, \ldots, t_r$ 
and $\OO^{\times}_{X,\bar{x}}$, and
$t_2, \ldots, t_r$ form a part of local parameters of 
$\OO_{\widetilde{X},\bar{x}}$. This shows us our assertion. 

\medskip
(Case C): In this case, 
by Proposition~\ref{prop:split:monoid:semistable:variety},
$P$ is of semistable type
\[
(r, l, p_1, \ldots, p_r, q_0, c_{l+1}, \ldots, c_{r}) 
\]
over $Q$.
Then, we have
\[
 k \otimes_{k[Q]} k[P] \simeq k[T_1, \ldots, T_r]/(T_1 \cdots T_l).
\]
via the correspondence $1 \otimes p_i \longleftrightarrow T_i$.
Adding indeterminates $T_{r+1}, \ldots, T_{n+1}$,
we have
\[
 k[T_1, \ldots, T_r, T_{r+1}, \ldots, T_{n+1}]_{(T_1, \ldots, T_{n+1})}
 /(T_1 \cdots T_{l}) \to \OO_{X,\overline{\nu(x)}}
\]
is etale.
We denote $\alpha(\pi_P(p_i))$ by $t_i$ for $i=1, \ldots, r$.
Renumbering $p_1, \ldots, p_r$, we may assume that
the component $\widetilde{X}$ at $x$ is given by $t_1 = 0$.
Note that $h(\bar{T}_i) = t_i$ for $i=1, \ldots, r$.
Thus,
$M_{\widetilde{X}, \bar{x}}$ is generated by $t_2, \ldots, t_r$ 
and $\OO^{\times}_{X,\bar{x}}$, and
$t_2, \ldots, t_r$ form a part of local parameters of 
$\OO_{\widetilde{X},\bar{x}}$. Hence, we get our assertion.

\bigskip
Next we claim the following:

\begin{Claim}
\label{prop:log:differential:smooth:comp:claim:2}
For $a \in M_{\widetilde{X}, \bar{x}}$,
there is $b \in \nu^*(M_X)_{\bar{x}}$ with
$\alpha'(b) = a$. Moreover, $b \otimes 1$ is uniquely determined
in $\nu^*(M_X)^{gr}_{\bar{x}} \otimes_{\ZZ} \OO_{\widetilde{X},\bar{x}}$.
\end{Claim}

The existence of $b$ is obvious, so that we consider only the uniqueness of $b$.
We use the same notation as in 
Claim~\ref{prop:log:differential:smooth:comp:claim:1} for each case.

\medskip
(Case A):
We set $a = u \cdot t_1^{a_1} \cdots t_r^{a_r}$ 
($u \in \OO^{\times}_{X,\bar{x}}$
and $a_1, \ldots, a_r\in \NN$).
In order to see the uniqueness of $b$, we set
$b = (f(q), b_1, \ldots, b_r, v)$
($q \in Q$, $b_1, \ldots, b_r \in \NN$ and 
$v \in \OO^{\times}_{X, \bar{x}}$).
Then, $\alpha'(b) = \beta(q)\cdot v \cdot t_1^{b_1} \cdots t_r^{b_r}$,
where $\beta$ is given by
\[
 \beta(q) = \begin{cases}
 1 & \text{if $q = 0$} \\
 0 & \text{if $q \not= 0$}.
\end{cases}
\]
Thus, $q = 0$, $v = u$ and $(b_1, \ldots, b_r) = (a_1, \ldots, a_r)$.

\medskip
(Case B):
We can set $a = u \cdot t_2^{a_2} \cdots t_r^{a_r}$
($u \in \OO^{\times}_{\widetilde{X},\bar{x}}$
and $a_2, \ldots, a_r\in \NN$).
Moreover, we set $b = 
(f(q), \bar{T}_1^{b_1} \cdot \bar{T}_2^{b_2} \cdots \bar{T}_r^{b_r}, v)$
($q \in Q$, $b_1, \ldots, b_r \in \NN$ and 
$v \in \OO^{\times}_{\widetilde{X}, \bar{x}}$).
Then, $\alpha'(b) = \beta(q) \cdot v \cdot t_2^{b_1 + b_2} \cdot t_3^{b_3}
\cdots t_r^{b_r}$.
Thus, 
\[
 \text{$q = 0$, $v = u$, $a_2 = b_1 + b_2$ and 
$(b_3, \ldots, b_r) = (a_3, \ldots, a_r)$.}
\]
Therefore, for $b' = (f(q'), \bar{T}_1^{b'_1} \cdot 
\bar{T}_2^{b'_2} \cdots \bar{T}_r^{b'_r}, v')$,
if $\alpha'(b) = \alpha'(b') = a$, then
\[
 b = b' + (0, (\bar{T}_2/\bar{T}_1)^c, 1)
\]
in $\nu^*(M_X)^{gr}_{\bar{x}}$ for some $c \in \ZZ$.
Here $\operatorname{char}(k) \not = 2$ and $(\bar{T}_2/\bar{T}_1)^2 = 1$.
Hence, $b \otimes 1 = b' \otimes 1$ in
$\nu^*(M_X)^{gr}_{\bar{x}} \otimes_{\ZZ} \OO_{\widetilde{X}, \bar{x}}$.

\medskip
(Case C):
We set $a = u \cdot t_2^{a_2} \cdots t_r^{a_r}$ 
($u \in \OO^{\times}_{\widetilde{X},\bar{x}}$
and $a_2, \ldots, a_r\in \NN$).
Let us see the uniqueness of $b$.
Let us set
$b = (f(q) + \sum_{i=1}^r b_i p_i, v)$
($q \in Q$, $b_1, \ldots, b_r \in \NN$ and 
$v \in \OO^{\times}_{\widetilde{X}, \bar{x}}$).
Then, $\alpha'(b) = \beta(q)\cdot v \cdot t_1^{b_1} \cdots t_r^{b_r}$.
Thus, $q = 0$, $v = u$, $b_1 = 0$ and
$(b_2, \ldots, b_r) = (a_2, \ldots, a_r)$.

\bigskip
By Claim~\ref{prop:log:differential:smooth:comp:claim:2},
there is a natural homomorphism
\[
 \gamma : \Omega^1_{\widetilde{X}}(\log(M_{\widetilde{X}}/k^{\times}))
 \to
 \Omega^1_{\widetilde{X}}(\log(\nu^*(M_X)/M_k)).
\]
Moreover, we have a natural homomorphism
\[
 \gamma' : \nu^*(\Omega^1_{X}(\log(M_X/M_k))) \to  
 \Omega^1_{\widetilde{X}}(\log(\nu^*(M_X)/M_k)).
\]

\begin{Claim}
$\gamma$ and $\gamma'$ are isomorphisms.
\end{Claim}

(Case A):
In this case, $\gamma'$ is an isomorphism around $x$.
We set $t_j = h(T_j)$ for $j=r+1, \ldots, n$.
Then, $d\log(t_1), \ldots, d\log(t_r), dt_{r+1}, \ldots, dt_{n}$
form a basis of 
$\Omega^1_{\widetilde{X}, \bar{x}}(\log(M_{\widetilde{X}}/k^{\times}))$.
Moreover, 
$d\log(e_1), \ldots, d\log(e_r), dt_{r+1}, \ldots, dt_{n}$
form a basis of
$\Omega^1_{\widetilde{X},\bar{x}}(\log(\nu^*(M_X)/M_k))$.
On the other hand,
$\gamma(d\log(t_i)) = d\log(e_i)$ for $i = 1, \ldots, r$
and $\gamma(d t_j) = d t_j$ for $j = r+1, \ldots, n$.
Thus, $\gamma$ is an isomorphism around $x$.

\medskip
(Case B):
We set $t_j = h(\bar{T}_j)$ for $j = r+1, \ldots, n+1$.
Then, 
\[
 d\log(t_2), \ldots, d\log(t_r), dt_{r+1}, \ldots, dt_{n+1}
\]
form a basis of 
$\Omega^1_{\widetilde{X}, \bar{x}}(\log(M_{\widetilde{X}}/k^{\times}))$.
Moreover,
$\gamma(d\log(t_i)) = d\log(\bar{T}_i)$ for $i = 2, \ldots, r$
and $\gamma(d t_j) = d t_j$ for $j = r+1, \ldots, n+1$.
Let $N'$ be the submonoid of $N$ generated by
$\bar{T}_2, \ldots, \bar{T}_r$.
Then, we can see that
$N^{gr} = {N'}^{gr} \times \langle \bar{T}_1/\bar{T}_2 \rangle$,
$(\bar{T}_1/\bar{T}_2)^2 = 1$ and
$N' \simeq \NN^{r-1}$.
Thus, if we set $M' = f(Q)\times N' \times
\OO^{\times}_{\widetilde{X}, \bar{x}}$,
then the natural homomorphism
\[
 \Omega^1_{\widetilde{X},\bar{x}}(\log(M'/M_k)) \to
 \Omega^1_{\widetilde{X},\bar{x}}(\log(\nu^*(M_X)/M_k))
\]
is an isomorphism because $\operatorname{char}(k) \not = 2$.
Moreover, $M'$ is log smooth over $M_k$.
Therefore, 
$\Omega^1_{\widetilde{X},\bar{x}}(\log(\nu^*(M_X)/M_k))$
is a free $\OO_{\widetilde{X},\bar{x}}$-module whose basis
is
\[
 d\log(\bar{T}_2), \ldots, d\log(\bar{T}_r), d\log(t_{r+1}), \ldots, d\log(t_{n+1}).
\]
Thus, $\gamma$ is an isomorphism.
On the other hand,
we can choose
\[
 d\log(\bar{T}_2), \ldots, d\log(\bar{T}_r), d\log(t_{r+1}), \ldots, d\log(t_{n+1})
\]
as a basis of
$\nu^* \Omega^1_X(\log(M_X/M_k))_{\bar{x}}$.
Thus, $\gamma'$ is also an isomorphism.

\medskip
(Case C):
We set $t_j = h(\bar{T}_j)$ for $j = r+1, \ldots, n+1$.
Then, 
\[
 d\log(t_2), \ldots, d\log(t_r), dt_{r+1}, \ldots, dt_{n+1}
\]
forms a basis of 
$\Omega^1_{\widetilde{X}, \bar{x}}(\log(M_{\widetilde{X}}/k^{\times}))$.
Moreover,
$\gamma(d\log(t_i)) = d\log(p_i)$ for $i = 2, \ldots, r$
and $\gamma(d t_j) = d t_j$ for $j = r+1, \ldots, n+1$.
Let $P'$ be the submonoid of $P$ generated by
$f(Q)$ and $p_2, \ldots, p_r$.
Then, since 
\[
 p_1 = -(p_2 + \cdots + p_l) + f(q_0) + \sum_{i > l} c_i p_i,
\]
we have ${P'}^{gr} = P^{gr}$.
Thus, if we set $M' = P' \times \OO^{\times}_{\widetilde{X}, \bar{x}}$,
then the natural homomorphism
\[
 \Omega^1_{\widetilde{X},\bar{x}}(\log(M'/M_k)) \to
 \Omega^1_{\widetilde{X},\bar{x}}(\log(\nu^*(M_X)/M_k))
\]
is an isomorphism.
Moreover, since $P' = f(Q) \times \NN^{r-1}$,
we can see $M'$ is log smooth over $M_k$.
Therefore, 
$\Omega^1_{\widetilde{X},\bar{x}}(\log(\nu^*(M_X)/M_k))$
is a free $\OO_{\widetilde{X},\bar{x}}$-module whose basis
is
\[
 d\log(p_2), \ldots, d\log(p_r), d\log(t_{r+1}), \ldots, d\log(t_{n+1}).
\]
Thus, $\gamma$ is an isomorphism.
On the other hand,
\[
 d\log(p_2), \ldots, d\log(p_r), d\log(t_{r+1}), \ldots, d\log(t_{n+1})
\]
is a basis of
$\nu^* \Omega^1_X(\log(M_X/M_k))_{\bar{x}}$.
Thus, $\gamma'$ is also an isomorphism.
\QED

\renewcommand{\theTheorem}{\arabic{section}.\arabic{subsection}.\arabic{Theorem}}
\renewcommand{\theClaim}{\arabic{section}.\arabic{subsection}.\arabic{Theorem}.\arabic{Claim}}
\renewcommand{\theequation}{\arabic{section}.\arabic{subsection}.\arabic{Theorem}.\arabic{Claim}}

\section{Geometric preliminaries}

\setcounter{Theorem}{0}
\subsection{Relative rational maps}
\label{subsec:relative:rat:map}

Let $k$ be an algebraically closed field,
$X$ and $Y$ a proper algebraic varieties over $k$, and
$T$ a reduced algebraic scheme over $k$.
Let $\Phi : X \times_k T \dasharrow Y \times_k T$
be a relative rational map over $T$, namely,
there is a dense open set $U$ of $X \times_k T$ such that
$\Phi$ is defined over $U$, $\Phi : U \to Y \times_k T$
is a morphism over $T$ and for all $t \in T$,
$U \cap (X \times \{ t\}) \not= \emptyset$.
In this subsection, we consider the following proposition.

\begin{Proposition}
\label{prop:subset:constructible}
\begin{enumerate}
\renewcommand{\labelenumi}{(\arabic{enumi})}
\item
$\{ t \in T \mid \text{$\rest{\Phi}{X \times \{ t \}}$ is dominant} \}$ is closed.

\item
$\{ t \in T \mid \text{$\rest{\Phi}{X \times \{ t \}}$ is separably dominant} \}$ is
locally closed.

\item
We assume that $X$ is normal.
Let $D_X$ and $D_Y$ be reduced divisors on $X$ and $Y$
respectively.
For a rational map $\phi : X \dasharrow Y$, we denote by $X_{\phi}$ 
the maximal open set over which $\phi$ is defined.
Then, 
\[
\left\{ t \in T \mid \text{$(\rest{\Phi}{X \times \{ t \}})^{-1}(D_Y) 
\subseteq D_X$ on $X_{\rest{\Phi}{X \times \{ t \}}}$} \right\}
\]
is constructible.

\item
Let $Z$ be a subvariety of $Y$. Then, $\{ t \in T \mid
\rest{\Phi}{X \times \{ t\}}(X) \subseteq Z\}$ is closed.

\item
Let $h : F \to G$ be a homomorphism of locally free sheaves
on $X \times_k T$ such that
$h_t : F_t \to G_t$ is not zero for every $t \in T$. 
Then, 
\[
 \{ t \in T \mid \text{the image
of $h_t : F_t \to G_t$ is rank one}\}
\]
is closed.
\end{enumerate}
\end{Proposition}

\Proof
(1) Let $Z$ be the closure of $\Phi(U)$ and $p : Z \to T$ the projection
induced by $Y \times_k T \to T$.
Since $Z$ is proper over $T$, it is well know that the function $T \to \ZZ$ given by
$t \mapsto \dim Z_t$ is upper semicontinuous.
Moreover, $\dim Z_t \leq \dim Y$ and the equality hold if and only if $Z_t = Y$.
Thus, we get (1).

\medskip
(2) By virtue of (1), we may assume that $\rest{\Phi}{X \times \{ t \}}$ is dominant
for all $t \in T$. In this case, we need to prove that it is open.
Then, this can be easily checked by Lemma~\ref{lem:upper:semicont:rank} and
the following fact:
Let $L$ be a finitely generated field over a field $K$. Then,
$\dim_L \Omega^1_{L/K} \geq \trdeg_K(L)$ and the equality holds if and only if
$L$ is separable over $K$.

\medskip
(3) First we assume that $T$ is normal. 
We may assume that $U$ is maximal.
Then, since $X \times_k T$ is normal,
for all $t \in T$, $\codim(X \times \{ t \} \setminus U) \geq 2$.
Thus, $(\rest{\Phi}{X \times \{ t \}})^{-1}(D_Y) \subseteq D_X$ on 
$X_{\rest{\Phi}{X \times \{ t \}}}$ if and only if
$(\rest{\Phi}{(X \times \{ t \})\cap U})^{-1}(D_Y) \subseteq D_X$.
Here we set $W = \Phi^{-1}(D_Y \times_k T) \setminus D_X \times_k T$ on $U$.
Let $q : W \to T$ be the projection induced by $X \times_k T \to T$.
Then, $t \not\in q(W)$ if and only if
$(\rest{\Phi}{(X \times \{ t \})\cap U})^{-1}(D_Y) \subseteq D_X$.

Next we consider a general case. Let $\pi : \widetilde{T} \to T$ be the normalization
of $T$.
Then,
\begin{multline*}
\left\{ t \in T \mid \text{$(\rest{\Phi}{X \times \{ t \}})^{-1}(D_Y) 
\subseteq D_X$ on $X_{\rest{\Phi}{X \times \{ t \}}}$} \right\} \\
 = \pi \left(
\left\{ \tilde{t} \in \widetilde{T} \mid 
\text{$(\rest{\Phi}{X \times \{  \tilde{t} \}})^{-1}(D_Y) 
\subseteq D_X$ 
on $X_{\rest{\Phi}{X \times \{  \tilde{t} \}}}$} \right\} \right)
\end{multline*}
Thus, we get (3).

\medskip
(4)
Let $W$ be the Zariski closure of
$\Phi^{-1}(Z \times_k T_1)$.
Then, $\rest{\Phi}{X \times \{ t\}}(X) \subseteq Z$ if and only if
$X \times \{ t \} = W_t$.
Since $W$ is proper over $T_1$,
it is well known that the function $T_1 \to \ZZ$ given by
$t \mapsto \dim W_t$ is upper semicontinuous.
Moreover, $\dim W_t \leq \dim X$ and the equality hold 
if and only if $W_t = X$.
Thus, $T$ is closed.

\medskip
(5)
Let $K$ be the function field of $X$. Let us consider
homomorphisms $F \otimes_k K \to G \otimes_k K$.
Since $h_t \not= 0$ for all $t \in T$,
we have (5) by Lemma~\ref{lem:upper:semicont:rank}.
\QED

\begin{Lemma}
\label{lem:upper:semicont:rank}
Let $K[X_1, \ldots, X_r]$ be the $r$-variable polynomial ring over 
a field $K$ and $k$ an algebraically closed  
subfield of $K$. 
Let $I$ be an ideal of $k[X_1, \ldots, X_r]$ and
$A(X_1, \ldots, X_r)$
an $n \times m$-matrix whose entries are elements of 
\[
 K[X_1, \ldots, X_r]/I K[X_1, \ldots, X_r].
\]
Then, the function given by
\[
 k^r \supseteq V(I) \ni (t_1, \ldots, t_r) \mapsto 
 \rank A(t_1, \ldots, t_r) \in \ZZ
\]
is lower semi-continuous, where
\[
 V(I) = \{ (x_1, \ldots, x_r) \in k^r \mid
f(x_1, \ldots, x_r) = 0\ \forall f \in I \}.
\]
\end{Lemma}

\Proof
Clearly we may assume that $I = \{ 0\}$.
Considering minors of the matrix $A(X_1, \ldots, X_r)$,
it is sufficient to see the following claim:

\begin{Claim}
For $f_1, \ldots, f_l \in K[X_1, \ldots, X_r]$, the set
\[
 \{ (x_1, \ldots, x_r) \in k^r \mid f_1(x_1, \ldots, x_r) =
  \cdots = f_l(x_1, \ldots, x_r) = 0 \}
\]
is closed.
\end{Claim}

Replacing
$K$ by a field generated by coefficients of
$f_1, \ldots, f_l$ over $k$,
we may assume that
$K$ is finitely generated over $k$. Since $k$ is algebraically closed,
$K$ is separated over $k$.
Thus, there are $T_1, \ldots, T_s$ of $K$ such that
$T_1, \ldots, T_s$ are algebraically independent over $k$ and
$K$ is a finite separable extension over $k(T_1, \ldots, T_s)$.
By taking the Galois closure of $K$ over $k(T_1, \ldots, T_s)$,
we may assume that $K$ is a Galois extension over $k(T_1, \ldots, T_s)$.
For $f = \sum_I a_I X^I \in K[X_1, \ldots, X_r]$ and
$\sigma \in \Gal(K/k(T_1, \ldots, T_s))$, we denote
$\sum_I \sigma(a_I) X^I$ by $f^{\sigma}$.
Here, we set
\[
 F_i = \prod_{\sigma \in \Gal(K/k(T_1, \ldots, T_s))} f_i^{\sigma}
\]
for $i = 1, \ldots, l$.
Then, $F_1, \ldots, F_l \in k(T_1, \ldots, T_l)[X_1, \ldots, X_r]$
and, for $(x_1, \ldots, x_r) \in k^r$,
\[
 F_i(x_1, \ldots, x_r) = 0
\quad\Longleftrightarrow\quad
f_i(x_1, \ldots, x_r) = 0
\]
for $i = 1, \ldots, l$. Indeed,
if $F_i(x_1, \ldots, x_r) = 0$, then
$f_i^{\sigma}(x_1, \ldots, x_r) = 0$ for some $\sigma \in
\Gal(K/k(T_1, \ldots, T_s))$, which implies
\[
 0 = \sigma^{-1}(f_i^{\sigma}(x_1, \ldots, x_r)) = f_i(x_1, \ldots, x_r).
\]
By the above observation, we may assume that $K = k(T_1, \ldots, T_s)$.
By multiplying some 
$\phi(T_1, \ldots, T_r) \in k[T_1, \ldots, T_s]$ to $f_i$,
we may further assume that
\[
 f_1, \ldots, f_l \in k[T_1, \ldots, T_s][X_1, \ldots, X_r].
\]
We set
\[
 f_i = \sum_{J} c_{i, J}T^J\qquad(c_{i, J} \in k[X_1, \ldots, X_r])
\]
for $i=1, \ldots, l$.
Then, for $(x_1, \ldots, x_r) \in k^r$,
\[
 f_i(x_1, \ldots, x_r) = 0
\quad\Longleftrightarrow\quad
 c_{i, J}(x_1, \ldots, x_r) = 0\ \ \forall J.
\]
Thus,
\begin{multline*}
 \{ (x_1, \ldots, x_r) \in k^r \mid f_i(x_1, \ldots, x_r) = 0 \ \
 \forall i\} \\
 = \{ (x_1, \ldots, x_r) \in k^r \mid c_{i,J}(x_1, \ldots, x_r) = 0 \ \ 
     \forall i, J \}.
\end{multline*}
Therefore, we get the claim.
\QED

\setcounter{Theorem}{0}
\subsection{Geometric trick for finiteness}

Let $k$ be an algebraically closed field.
Let $X$ be a proper normal variety over $k$ and
$Y$ a proper algebraic variety over $k$.
Let $E$ be a vector bundle on $X$
and $H$ a line bundle on $Y$.
We assume that there is a dense open set $Y_0$ of $Y$
such that $H^0(Y, H) \otimes_k \OO_Y \to H$ is surjective over $Y_0$.
Let $\phi : X \dasharrow Y$ be a dominant rational map over $k$.
Let $X_{\phi}$ be the maximal open set of $X$ over which $\phi$ is defined.
We also assume that
there is a non-trivial homomorphism
$\theta : \phi^*(H) \to \rest{E}{X_{\phi}}$.
Then, since $\codim(X \setminus X_{\phi}) \geq 2$, we have 
a sequence of homomorphisms
\[
 H^0(Y, H) \to H^0(X_{\phi}, \phi^*(H)) \to H^0(X_{\phi}, E) = H^0(X, E).
\]
We denote the composition of the above homomorphisms by
\[
 \beta(\phi, \theta) : H^0(Y, H) \to H^0(X, E).
\]
Then, we have the following.

\begin{Lemma}
\label{lemma:composition:rational:maps}
Let $L$ be the image of
\[
\begin{CD}
 H^0(Y, H) \otimes_k \OO_X @>{\beta(\phi, \theta) \otimes_k \operatorname{id}}>> 
 H^0(X, E) \otimes_k \OO_X @>>> E.
\end{CD}
\]
Then, the rank of $L$ is one and
the rational map
\[
 \phi' : X \dasharrow \PP(H^0(Y, H))
\]
induced by $H^0(Y, H) \otimes_k \OO_X \to L$
is the composition of rational maps
\[
 X \overset{\phi}{\dasharrow} Y 
   \overset{\phi_{\vert H \vert}}{\dasharrow}
 \PP(H^0(Y, H)),
\]
namely, $\phi' = \phi_{\vert H \vert} \cdot \phi$.  
\end{Lemma}

\Proof
Considering the following commutative diagram:
\[
 \begin{CD}
  H^0(Y, H) \otimes_k \OO_{X_{\phi}} @>{\beta(\phi, \theta) \otimes_k \operatorname{id}}>> 
  H^0(X, E) \otimes_k \OO_{X_{\phi}} \\
  @VVV @VVV \\
  \phi^*(H) @>{\theta}>> \rest{E}{X_{\phi}}.
 \end{CD}
\]
We can see that $\theta$ gives rise to an isomorphism
\[
 \rest{\phi^*(H)}{X_{\phi} \cap \phi^{-1}(Y_0)} 
 \overset{\sim}{\longrightarrow} \rest{L}{X_{\phi} \cap \phi^{-1}(Y_0)}.
\]
Moreover, the rational map $X_{\phi} \dasharrow \PP(H^0(Y, H))$ given
by $H^0(Y, H) \otimes_k \OO_{X_{\phi}} \to \phi^*(H)$ 
is $\phi_{\vert H \vert} \cdot \phi$.
Thus, the rational map $\phi' : X \dasharrow \PP(H^0(Y, H))$
induced by $H^0(Y, H) \otimes_k \OO_X \to L$
is nothing more than the composition of rational maps
\[
 X \overset{\phi}{\dasharrow} Y 
   \overset{\phi_{\vert H \vert}}{\dasharrow}
 \PP(H^0(Y, H)).
\]
\QED

From now on, we assume that $H$ is very big, that is,
the morphism $Y_0 \to \PP(H^0(Y,H))$ induced by
$H^0(Y, H) \otimes_k \OO_{Y_0} \to \rest{H}{Y_0}$ is a birational morphism.
Let $\mathcal{C}$ be a subset of $\Rat_k(X, Y)$
(the set of all rational maps of $X$ into $Y$).
We assume that for all $\phi \in \mathcal{C}$,
\begin{enumerate}
\renewcommand{\labelenumi}{(\arabic{enumi})}
\item
$\phi$ is a dominant rational map, and

\item
we can attach a non-trivial homomorphism
$\theta_{\phi} : \phi^*(H) \to \rest{E}{X_{\phi}}$ to $\phi$,
where $X_{\phi}$ is the maximal Zariski open set of $X$
over which $\phi$ is defined.
\end{enumerate}
As before, we have an homomorphism
\[
 \beta(\phi, \theta_{\phi}) : H^0(Y, H) \to H^0(X, E).
\]
We denote the class of $\beta(\phi, \theta_{\phi})$ in 
$\PP(\Hom_k(H^0(Y, H), H^0(X, E))^{\vee})$ by $\gamma(\phi)$.

\begin{Lemma}
\label{lem:injectibity:gamma}
For $\phi, \psi \in \mathcal{C}$,
if $\gamma(\phi) = \gamma(\psi)$, then $\phi = \psi$.
\end{Lemma}

\Proof
By our assumption, there is $a \in k^{\times}$ with
$a \beta(\phi) = \beta(\psi)$.
Hence, we have the following commutative diagram:
\[
 \begin{CD}
 H^0(Y, H) \otimes_k \OO_X 
 @>{\beta(\phi, \theta_{\phi}) \otimes_k \operatorname{id}}>> 
 H^0(X, E) \otimes_k \OO_X @>>> E \\
 @| @VV{\times a}V @VV{\times a}V \\
 H^0(Y, H) \otimes_k \OO_X 
 @>{\beta(\psi, \theta_{\psi}) \otimes_k \operatorname{id}}>> 
 H^0(X, E) \otimes_k \OO_X @>>> E
 \end{CD}
\]
Let $L_{\phi}$ (resp. $L_{\psi}$) be the image
of $H^0(Y, H) \otimes_k \OO_X \to E$ in terms of $\beta(\phi, \theta_{\phi})$
(resp. $\beta(\psi, \theta_{\psi})$).
Then, the above diagram gives rise to a commutative diagram
\[
 \begin{CD}
 H^0(Y, H) \otimes_k \OO_X @>>> L_{\phi} \\
 @| @VV{\times a}V \\
 H^0(Y, H) \otimes_k \OO_X @>>> L_{\psi}.
 \end{CD}
\]
Let $\phi' : X \dasharrow \PP(H^0(Y,H))$ and
$\psi' : X \dasharrow \PP(H^0(Y,H))$
be the rational maps
induced by $H^0(Y, H) \otimes_k \OO_X \to L_{\phi}$ and
$H^0(Y, H) \otimes_k \OO_X \to L_{\psi}$ respectively.
Then, by the above diagram, we can see $\phi' = \psi'$.
Hence, we get our lemma by Lemma~\ref{lemma:composition:rational:maps}.  
\QED

Next we consider the following proposition.

\begin{Proposition}
\label{prop:relative:rational:map:constant}
Let $T$ be a connected proper normal variety over $k$, and
\[
\Phi : X \times_k T \dasharrow Y \times_k T
\]
be a relative rational map over $T$
\rom{(}cf. Conventions and terminology 
\rom{\ref{subsubsec:rational:map}}\rom{)}.
Let $f : X \times_k T \to T$ and $g : Y \times_k T \to T$ be the projections to the
second factor respectively, and let $p : X \times_k T \to X$ and $q : Y \times_k T \to Y$ be the
projections to the first factor respectively.
We assume that there is an open set $T_0$ of $T$ and
a non-trivial homomorphism $\Theta : \Phi^*(q^*(H)) \to \rest{p^*(E)}{U}$  such that,
for all $t \in T_0$,
$\rest{\Phi}{X \times \{ t\}} \in \mathcal{C}$ and the class of $\beta(\Phi_t, \Theta_t)$  in
$\PP(\Hom_k(H^0(Y, H), H^0(X, E))^{\vee})$ is $\gamma(\Phi_t)$,
where $U$ is the maximal open set over which $\Phi$ is defined.
Then, there is $\phi \in \mathcal{C}$ such that $\Phi = \phi \times \operatorname{id}_T$.
\end{Proposition}

\Proof
Since $X \times_k T$ is normal,
we may assume that $\codim((X \times_k T) \setminus U) \geq 2$.
Here we have a homomorphism
\[
H^0(Y, H) \otimes_k \OO_T = g_*(q^*(H)) \to
(\rest{f}{U})_* (\Phi^*(q^*(H))) \overset{\Theta}{\longrightarrow}
(\rest{f}{U})_*(p^*(E)).
\]
We claim that the natural homomorphism
$f_*(p^*(E)) \to (\rest{f}{U})_*(p^*(E))$ is an isomorphism.
Indeed, if $W$ is an open set of $T$, then
\[
 (\rest{f}{U})_*(p^*(E))(W)  = H^0(U \cap (X \times_k W), p^*(E)).
\]
Note that $\codim((X \times_k W) \setminus U \cap (X \times_k W)) \geq 2$.
Thus, $H^0(U \cap (X \times_k W), p^*(E)) = H^0(X \times_k W, p^*(E))$.
Hence we get a homomorphism
\[
\beta : H^0(Y, H) \otimes_k \OO_T \to H^0(X, E) \otimes \OO_T.
\]
Here, $T$ is proper and irreducible. Hence, there is
$\beta_0 \in \Hom_k(H^0(Y, H), H^0(X, E))$ such that
$\beta = \beta_0 \otimes \operatorname{id}$.
This means that $\beta(\Phi_t, \Theta_t) = \beta_0$.
Thus, by Lemma~\ref{lem:injectibity:gamma},
there is $\phi \in \mathcal{C}$ such that $\Phi_t = \phi$ for all $t \in T_0$.
Therefore, we get our proposition.
\QED

Finally, let us see the following proposition.

\begin{Proposition}
\label{prop:universal:rat:map}
There is a closed subset $T$ of $\PP(\Hom_k(H^0(Y, H), H^0(X, E))^{\vee})$ and
a relative rational map $\Phi : X \times_k T \dasharrow Y \times_k T$ over $T$
such that if we consider $\gamma : \mathcal{C} \to \PP(\Hom_k(H^0(Y, H), H^0(X, E))^{\vee})$,
then $\gamma(\mathcal{C}) \subseteq T$ and $\rest{\Phi}{X \times \{ \gamma(\phi) \}} = \phi$.
\end{Proposition}

\Proof
We set $P = \PP(\Hom_k(H^0(Y, H), H^0(X, E))^{\vee})$.
Then, there is the canonical homomorphism
\[
 \Hom_k(H^0(Y, H), H^0(X, E))^{\vee} \otimes_k \OO_P \to \OO_P(1),
\]
which gives rise to a universal homomorphism
\[
 \beta : H^0(Y, H) \otimes_k \OO_P(-1) \to H^0(X, E) \otimes_k \OO_P,
\]
that is, for all $t \in P$, the class of
\[
 \beta_t : H^0(Y, H) \otimes_k (\OO_P(-1) \otimes \kappa(t)) \to H^0(X, E)
\]
in $P$ coincides with $t$, where $\kappa(t)$ is the residue field of $\OO_P$
at $t$.
Here we consider the composition of homomorphisms
\[
 h : H^0(Y, H) \otimes_k \OO_P(-1) \otimes_k \OO_X 
\overset{\beta \otimes \operatorname{id}}{\longrightarrow}
 H^0(X, E) \otimes_k \OO_P \otimes_k \OO_X \to E \otimes_k \OO_P
\]
on $X \times_k P$.
Then, by (5) of Proposition~\ref{prop:subset:constructible},
if $T_1$ is the set of all $t \in P$ 
such that the image of $h_t$ is of rank $1$, then $T_1$ is closed.
Let $L$ be the image of
\[
 \rest{h}{T_1} : H^0(Y, H) \otimes_k \OO_{T_1}(-1) \otimes_k \OO_{X} \to
  E \otimes_{k} \OO_{T_1}.
\]
Then, we have the surjective homomorphism
\[
 H^0(Y, H) \otimes_k \OO_{X \times_k T_1} \to
 L \otimes_{\OO_{X \times_k T_1}}  \OO_{X \times_k T_1}(1).
\]
Let $U_1$ be the maximal Zariski open set of $X \times_k T_1$
such that $L$ is invertible over $U_1$.
Here, note that, for all $t \in T_1$, 
$U_1 \cap \left(X \times_k \{ t_1 \} \right) \not= \emptyset$.
Thus, we get a relative rational map
\[
 \Phi_1 : X \times_k T_1 \dasharrow \PP(H^0(Y,H)) \times_k T_1
\]
over $T_1$ (cf. Conventions and terminology 
\ref{subsubsec:rational:map}).
Let $Y_1$ be the closure of the image of $\phi_{\vert H \vert}(Y)$.
By (4) of Proposition~\ref{prop:subset:constructible}, 
that
\[
T = \{ t \in T_1 \mid (\Phi_1)_t(X) \subseteq Y_1 \}
\]
is closed.
Hence we obtain a relative rational map
\[
 \Phi_2 : X \times_k T \dasharrow Y_1 \times_k T,
\]
which gives rise to a relative rational map
\[
 \Phi : X \times_k T \dasharrow Y \times_k T.
\]
By our construction, this rational map has the following properties:
For all $t \in T$, let $\beta_t : H^0(Y, H) \to H^0(X, E)$ be the
homomorphism modulo $k^{\times}$ corresponding to $t \in P$, and
$L_t$ the image of 
\[
H^0(Y, H) \otimes \OO_X \to H^0(X, E) \otimes \OO_X \to E.
\]
Here, the rank of $L_t$ is one. Thus, we have a rational map
$\phi_t : X \dasharrow \PP(H^0(Y, H))$ induced by
$H^0(Y, H) \otimes \OO_X \to L_t$.
Then, $\phi_t(X) \subseteq Y_1$ and the following diagram is commutative:
\[
\xymatrix{
X \ar[rr]^{\rest{\Phi}{X \times \{ t \}}} \ar[rd]_{\phi_t} &  &  Y \ar[ld]^{\phi_{|H|}} \\
& Y_1 & 
}
\]
Therefore, by Lemma~\ref{lemma:composition:rational:maps},
$\Phi : X \times_k T \dasharrow Y \times_k T$ is our desired relative 
rational map.
\QED

\renewcommand{\theTheorem}{\arabic{section}.\arabic{Theorem}}
\renewcommand{\theClaim}{\arabic{section}.\arabic{Theorem}.\arabic{Claim}}
\renewcommand{\theequation}{\arabic{section}.\arabic{Theorem}.\arabic{Claim}}
\section{Log smooth case over the trivial log structure}

Let $k$ be an algebraically closed field and let
$X$ and $Y$ be proper normal algebraic varieties over $k$.
Let $D_X$ and $D_Y$ be reduced divisors on $X$ and $Y$ respectively.
Let $M_X$ and $M_Y$ be fine log structures of
$X$ and $Y$ respectively such that
\[
 M_X = {j_X}_*(\OO^{\times}_{X \setminus D_X}) \cap \OO_X
\quad\text{and}\quad
 M_Y \subseteq {j_Y}_*(\OO^{\times}_{Y \setminus D_Y}) \cap \OO_Y,
\]
where $j_X$ and $j_Y$ are natural inclusion maps
$X \setminus D_X \hookrightarrow X$ and
$Y \setminus D_Y \hookrightarrow Y$ respectively.
Then, for a rational map $\phi : X \dasharrow Y$,
$\phi$ extends to $(X, M_X) \to (Y, M_Y)$ if
$\phi^{-1}(D_Y) \subseteq D_X$.
We assume that $(X, M_X)$ and $(Y, M_Y)$ are
log smooth over $(\Spec(k), k^{\times})$.
Note that if $X$ is smooth over $k$,
then the log smoothness of $(X, M_X)$
over $(\Spec(k), k^{\times})$
guarantees that $M_X = {j_X}_*(\OO^{\times}_{X \setminus D_X}) \cap \OO_X$
for $D_X = \Supp(M_X/\OO^{\times}_X)$ (cf. Proposition~\ref{prop:split:monoid:semistable:variety}).
Moreover, we assume that $(Y, M_Y)$ is of log general type, namely,
$\det \Omega^1_Y(\log(M_Y/k^{\times}))$ is big.
Thus, there is a positive integer $m$ such that
$\det \Omega^1_Y(\log(M_Y/k^{\times}))^{\otimes m}$
is very big.
Here we set 
\[
 H = \det \Omega^1_Y(\log(M_Y/k^{\times}))^{\otimes m}
\quad\text{and}\quad
 E = \Sym^m(\wedge^{\dim Y} \Omega^1_X(\log(M_X/k^{\times}))).
\]
Then, if $\phi : (X, M_X) \dashrightarrow (Y, M_Y)$ is a rational map, 
then we have a natural homomorphism
\[
 \theta_{\phi} : \phi^*(H) \to \rest{E}{X_{\phi}},
\]
where $X_{\phi}$ is the maximal open set over which $\phi$
is defined. Moreover, if $\phi$ is separably dominant,
then $\theta_{\phi}$ is non-trivial.
Let $\SDRat((X, M_X), (Y, M_Y))$ be the set of separably dominant 
rational maps $(X, M_X) \dashrightarrow (Y, M_Y)$ over 
$(\Spec(k), k^{\times})$.

\begin{Theorem}
\label{thm:finite:dominant:rat:map:general:type:smooth}
$\SDRat((X, M_X), (Y, M_Y))$ is finite.
\end{Theorem}

\Proof
First we need the following lemma.

\begin{Lemma}
\label{lem:log:map:const}
Let $T$ be a smooth proper curve over $k$ and
$\Phi : X \times_k T \dasharrow Y \times_k T$ a relative 
rational map over $T$ \rom{(}cf. Conventions and terminology
\rom{\ref{subsubsec:rational:map}}\rom{)}.
If there is a non-empty open set $T_0$ of $T$ such that for all $t \in T_0$,
$\Phi_t$ is separably dominant and $\Phi_t^{-1}(D_Y) \subseteq D_X$, 
then there is a rational map
$\phi : X \dasharrow Y$ with $\Phi = \phi \times \operatorname{id}_T$.
\end{Lemma}

\Proof
First of all,
by Proposition~\ref{prop:subset:constructible},
for all $t \in T$, $\rest{\Phi}{X \times \{t \}} : X \dasharrow Y$ is dominant.
Let us take a effective divisor $D$ on $X$ such that
\[
\rest{\Phi}{X \times \{t\}}^{-1}(D_Y)
\subseteq D_X \cup D
\]
for all $t \in T \setminus T_0$. By using de-Jong's alteration \cite{DeJ}, 
there are a smooth proper variety $X'$ and
a separable and generically finite morphism $\mu : X' \to X$ such that
$\mu^{-1}(D_X \cup D)$ is a normal crossing divisor on $X'$.
Let $D_{X'} = \mu^{-1}(D_X \cup D)$ and
$M_{X'} = {j_{X'}}_*(\OO^{\times}_{X' \setminus D_{X'}}) \cap \OO_{X'}$,
where $j_{X'} : X' \setminus D_{X'} \to X'$ is the natural inclusion map.
Then, $(X', M_{X'})$ is log smooth over $(\Spec(k), k^{\times})$.
We set $\Phi' = \Phi \cdot (\mu \times \operatorname{id}_T)$. Then,
for all $t \in T$, $\rest{\Phi'}{X \times \{ t\}}^{-1}(D_Y)
\subseteq D_{X'}$.
Moreover, for all $t \in T_0$,
$\rest{\Phi'}{X \times \{ t\}}$ is separably dominant.
Thus, in order to prove our lemma,  we may assume that
for all $t \in T$,  $\rest{\Phi}{X \times \{ t\}}^{-1}(D_Y) \subseteq D_X$.

Let $f : X \times_k T \to T$ and $g : Y \times_k T \to T$ be 
the projections to the
second factor respectively, and let $p : X \times_k T \to X$ and $q : Y \times_k T \to Y$ be the
projections to the first factor respectively.
Let $U$ be the maximal open set over which $\Phi$ is defined.
Then, we have a rational map
$(X \times_k T, p^*(M_X)) \dashrightarrow (Y \times_k T, q^*(M_Y))$ and
$(X \times_k T, p^*(M_X))$ and $(Y \times_k T, q^*(M_Y))$
are log smooth over $(T, \OO^{\times}_T)$.
Thus, there is a non-trivial homomorphism
\[
\Theta : \Phi^*(q^*(H)) \to 
\rest{p^*(E)}{U}.
\]
Therefore, we get our lemma 
by Proposition~\ref{prop:relative:rational:map:constant}.
\QED

\medskip
Let us go back to the proof of 
Theorem~\ref{thm:finite:dominant:rat:map:general:type:smooth}.
If $\phi \in \SDRat((X, M_X), (Y, M_Y))$, 
then we have the non-trivial homomorphism
\[
\theta_{\phi} : \phi^{*}(H) \to \rest{E}{X_{\phi}}.
\]
Thus, by Proposition~\ref{prop:universal:rat:map},
there is a closed subset $T$ of 
\[
 \PP(\Hom_k(H^0(Y, H), H^0(X, E))^{\vee})
\]
and a relative rational map 
$\Phi : X \times_k T \dasharrow Y \times_k T$ over $T$
such that if we consider 
\[
 \gamma : \SDRat((X, M_X), (Y, M_Y))
\to \PP(\Hom_k(H^0(Y, H), H^0(X, E))^{\vee}),
\]
then
\[
 \gamma(\SDRat((X, M_X), (Y, M_Y))) \subseteq T
\]
and $\rest{\Phi}{X \times \{ \gamma(\phi) \}} = \phi$.
Note that $\gamma$ is injective by Lemma~\ref{lem:injectibity:gamma}.
Let $T_1$ be the set of all $t \in T$ such that
$\rest{\Phi}{X \times \{ t \}}$ is separably dominant and
$\rest{\Phi}{X \times \{ t \}}^{-1}(D_Y) 
 \subseteq D_X$.
Then, by Proposition~\ref{prop:subset:constructible},
$T_1$ is constructible. Let $T_2$ be the Zariski closure of $T_1$.
If $\dim T_2 = 0$, then we have done, so that we assume that $\dim T_2 > 0$.
Then, there is a proper smooth curve $C$ and $\pi : C \to T_2$ such that
the generic point of $C$ goes to $T_1$ via $\pi$.
Moreover, we have a rational map $\Psi : X \times_k C \dasharrow Y \times_k C$ induced
by $X \times_k T_2 \dasharrow Y \times_k T_2$.
By our construction, there is an open set $C_0$ of $C$ such that
for all $t \in C_0$, $\rest{\Psi}{X \times_k C_0}$ is separably dominant and
$\rest{\Psi}{X \times \{ t \}}^{-1}(D_Y) \subseteq D_X$. Thus, by Lemma~\ref{lem:log:map:const},
there is a rational map $\psi : X \dasharrow Y$ with
$\Psi = \psi \times \operatorname{id}$.
We choose $x_1, x_2\in C$ with $\pi(x_1) \not= \pi(x_2)$ and $\pi(x_1), \pi(x_2) \in T_1$.
Then, we have $\phi_1, \phi_2 \in \SDRat((X, M_X), (Y, M_Y))$ 
with $\gamma(\phi_1) = \pi(x_1)$ and
$\gamma(\phi_2) = \pi(x_2)$. Since $\gamma$ is injective, $\phi_1 \not= \phi_2$.
On the other hand,
\[
\psi = \rest{\Psi}{X \times_k \{ x_i\}} = \rest{\Phi}{X \times_k \{ \pi(x_i) \}} = \phi_i
\]
for each $i$.
This is a contradiction.
\QED

\section{The proof of the main theorem}

In this section, let us consider the proof of the main theorem of this paper.

\begin{Theorem}
\label{thm:finite:dominant:rat:map:general:type:semi:stable}
Let $k$ be an algebraically closed field and
$M_k$ a fine log structure of $\Spec(k)$.
Let $X$ and $Y$ be proper semistable varieties over $k$,
and let $M_X$ and $M_Y$ be fine log structures of $X$ and $Y$ respectively.
We assume that $(X, M_X)$ and $(Y, M_Y)$ are integral and smooth over
$(\Spec(k), M_k)$.
If $(Y, M_Y)$ is of log general type, then
the set of all separably dominant rational maps
$(X, M_X) \dasharrow (Y, M_Y)$ over $(\Spec(k), M_k)$ 
defined in codimension one is finite
\rom{(}see Conventions and terminology \rom{\ref{subsubsec:rational:map}}\rom{)}.
\end{Theorem}

\Proof
First we need the following lemma:

\begin{Lemma}
\label{lem:big:irreducible:comp}
Let $Y$ be a semistable variety over $k$ and
$H$ a line bundle on $Y$. 
Let $Y'$ be an irreducible component of the normalization
of $Y$ and $\mu : Y' \to Y$ the natural morphism.
If $H$ is big, then
$\mu^*(H)$ is big.
\end{Lemma}

\Proof
Let $m$ be a positive integer $m$ such that
$H^{\otimes m}$ is very big.
Let $V$ be the image of 
$H^0(Y, H^{\otimes m}) \to H^0(Y',\mu^*(H^{\otimes m}))$.
Then, we have the following diagram:
\[
 \xymatrix{
 Y' \ar[r]^{\mu} \ar@{-->}[rrd] \ar@{-->}[rrdd] & Y
 \ar@{-->}[r] & \PP(H^0(Y, H^{\otimes m})) \\
                    &          & \PP(V) \ar@{^{(}->}[u] \\
                    &          & \PP(H^0(Y', f^*(H^{\otimes m}))) \ar@{-->}[u]
 }
\]
Let $Y_1$ and $Y_2$ be the image of
$Y' \dasharrow \PP(V)$ and 
$Y' \dasharrow \PP(H^0(Y', \mu^*(H^{\otimes m})))$ respectively.
Then,
\[
 k(Y') = k(Y_1) \subseteq k(Y_2) \subseteq k(Y').
\]
Thus, we can see that 
$Y' \dasharrow Y_2$
is birational.
\QED

\bigskip
Let us go back to the proof of 
Theorem~\ref{thm:finite:dominant:rat:map:general:type:semi:stable}.
Let $X_1, \ldots, X_r$ and $Y_1, \ldots, Y_s$ be irreducible components
of the normalizations of $X$ and $Y$ respectively.
Moreover, let $f_i : X_i \to X$ and $g_j : Y_j \to Y$ be the canonical
morphisms.
We set $M_{X_i} = f_i^*(M_X)^u$ and $M_{Y_j} = g_j^*(M_Y)^u$
(cf. see Conventions and terminology \ref{subsubsec:underlining:log:structure}).
Then, by Proposition~\ref{prop:log:differential:smooth:comp},
$(X_i, M_{X_i})$ and $(Y_j, M_{Y_j})$ are integral and log smooth
over $(\Spec(k), k^{\times})$. Further, by Proposition~\ref{prop:log:differential:smooth:comp}
again, 
\[
\Omega^1_{X_i}(\log(M_{X_i})) = f_i^*(\Omega^1_{X}(\log(M_X/M_k)))
\]
and
\[
\Omega^1_{Y_j}(\log(M_{Y_j})) = g_j^*(\Omega^1_{Y}(\log(M_Y/M_k))).
\]
Thus, by the above lemma,
$(Y_j, M_{Y_j})$ is of log general type for every $j$.
We denote the set of all separably dominant rational maps
$(X, M_X) \dasharrow (Y, M_Y)$ defined in codimension one over
$(\Spec(k), M_k)$ by 
\[
 \SDRat((X, M_X), (Y, M_Y)).
\]
Moreover, the set of all separably dominant rational maps
$(X_i, M_{X_i}) \dasharrow (Y_j, M_{Y_j})$ over $(\Spec(k), k^{\times})$
is denoted by 
\[
 \SDRat((X_i, M_{X_i}), (Y_j, M_{Y_j})).
\]
Then, we have a natural map
\[
\Psi : \SDRat((X, M_{X}), (Y, M_{Y})) \longrightarrow
\coprod_{\sigma \in S(r,s)}
\prod_{i=1} ^r\SDRat((X_i, M_{X_i}), (Y_{\sigma(i)}, M_{Y_{\sigma(i)}}))
\]
as follows. Here $S(r, s)$ is the set all maps from $\{ 1, \ldots, r\}$ to $\{ 1, \ldots, s \}$.
Let $(\phi,  h) \in  \SDRat((X, M_{X}), (Y, M_{Y}))$.
Then, for each $i$, there is a unique $\sigma(i)$ such that the Zariski closure of 
$\phi(X_i)$ is $Y_{\sigma(i)}$.
Then, we have $(\rest{\phi}{X_i}, h_i) : (X_i, M_{X_i}) \to (Y_{\sigma(i)}, M_{Y_{\sigma(i)}})$.
By Theorem~\ref{thm:finite:dominant:rat:map:general:type:smooth},
$\SDRat((X_i, M_{X_i}), (Y_j, M_{Y_j}))$'s are finite sets.
Therefore, it is sufficient to see that $\Psi$ is injective.
Let us pick up $(\phi, h), (\phi' , h') \in \SDRat((X, M_{X}), (Y, M_{Y}))$ with
$\Psi(\phi) = \Psi(\phi')$. 
Then, clearly, $\phi = \phi'$. Thus, by Theorem~\ref{thm:unqueness:log:hom:semistable},
we have $h = h'$.
\QED

\renewcommand{\thesection}{}

\renewcommand{\theTheorem}{A.\arabic{Theorem}}

\renewcommand{\theClaim}{A.\arabic{Theorem}.\arabic{Claim}}

\renewcommand{\theequation}{A.\arabic{Theorem}.\arabic{Claim}}
\setcounter{Theorem}{0}

\section*{Appendix}

In this appendix, we consider several results, which are 
well known facts for researchers of log geometry.
It is however difficult to find references, so that
for reader's convenience, we prove them here.
First, let us consider irreducible elements of a fine and sharp monoid.

\begin{Proposition}
\label{prop:fine:sharp:irreducible}
Let $P$ be a fine and sharp monoid.
Then, $P$ is generated by irreducible elements and
there are finitely many irreducible elements of $P$.
\end{Proposition}

\Proof
In this proof, the binary operation of $P$ is written by product.
We define a vector subspace $M$ of $\QQ[P]$ to be
\[
M = \bigoplus_{x \in P \setminus \{ 1 \}} \QQ x.
\]
Here we claim $M$ is a maximal ideal of $\QQ[P]$.
For $x \in P$ and $x' \in P \setminus \{ 1 \}$,
we have $x \cdot x' \in P \setminus \{ 1 \}$ because $P$ is sharp.
This shows us that
$M$ is an ideal. Moreover, $\QQ[P]/M \simeq \QQ$.
Thus, we get the claim. We set $R = \QQ[P]_M$ 
(the localization at $M$) and $m = M\QQ[P]_M$.
Note that $\bigcap_{n \geq 0} m^n = \{ 0\}$ because
$R$ is a noetherian local ring.
Moreover, since $P$ is integral,
the natural map $P \to R$ is injective and
$x \not= 0$ in $R$ for all $x \in P$.

For $x \in P$, we define $\deg(x)$ to be
\[
\deg(x) = \max \{ n \in \NN \mid x \in m^n \}.
\]
Then, it is easy to see that $\deg(x) = 0$ if and only if
$x = 1$ and $\deg(x \cdot y) \geq \deg(x) + \deg(y)$ for
$x, y \in P$.
We say $x$ is decomposable by irreducible elements if
there are irreducible elements $p_1, \ldots, p_s$ such that
$x = p_1 \cdots p_s$.
Here we set 
\[
\Sigma = \{ x \in P \setminus \{ 1\}
\mid \text{$x$ is not decomposable by irreducible elements} \}.
\]
We would like to show $\Sigma = \emptyset$.
We assume the contrary. Let us choose $x \in \Sigma$ such that
$\deg(x)$ is minimal in $\{ \deg(y) \mid y \in \Sigma \}$.
Then, $x$ is not irreducible, so that we have a decomposition
$x = y \cdot z$ ($y \not= 1$ and $z \not= 1$).
Then, $\deg(x) \geq \deg(y) + \deg(z)$, $\deg(y) \not= 0$ and $\deg(z) \not= 0$.
Thus, $\deg(y), \deg(z) < \deg(x)$, which implies 
$y, z \not\in \Sigma$. Therefore, $y$ and $z$ are decomposable by
irreducible elements. Thus, so does $x$. This is a contradiction.

Next, let us see that we have only finitely many irreducible elements.
Since $P$ is finitely generated, there is a surjective homomorphism
$h : \NN^n \to P$.  Let $p$ be an irreducible element of $P$.
Let us choose $I \in \NN^n$ such that $h(I) = p$ and
$\deg(I)$ is minimal in $\{ \deg(J) \mid h(J) = p \}$.
Here we claim that $I$ is irreducible in $\NN^n$.
We suppose $I = I' + I''$ ($I' \not= 0$ and $I'' \not= 0$).
Then, $h(I') \cdot h(I'') = p$. Here $p$ is irreducible.
Thus, either $h(I') = 1$ or $h(I'') = 1$, which means that
either $h(I') = p$ or $h(I'') = p$.
This is a contradiction because $\deg(I'), \deg(I'') < \deg(I)$. 
Therefore, $I$ is irreducible.
Note that an irreducible element of $\NN^n$ has a form
$(0, \ldots, 1, \ldots, 0)$.
Hence, we have only finitely many irreducible elements.
\QED

Finally, let us consider two propositions concerning
the existence of a good chart of a smooth log morphism
(cf. \cite{OgusDeRham}).

\begin{Proposition}
\label{prop:tame:good:chart}
Let $(\phi, h) : (X, M_X) \to (Y, M_Y)$ be
a morphism of log schemes with
fine log structures. Let $x \in X$ and $y = \phi(x)$.
We assume the following:
\begin{enumerate}
\renewcommand{\labelenumi}{(\arabic{enumi})}
\item
The homomorphism 
$\bar{h}_x : \overline{M}_{Y, \bar{y}} \to \overline{M}_{X, \bar{x}}$ 
induced by $h_x : M_{Y, \bar{y}} \to M_{X, \bar{x}}$ is
injective and the torsion part of 
$\Coker({\bar{h}}^{gr}_x : \overline{M}^{gr}_{Y, \bar{y}} \to 
\overline{M}^{gr}_{X, \bar{x}})$
is a finite group of order invertible in $\OO_{X, \bar{x}}$.

\item
There is a splitting homomorphism
$s_y : \overline{M}_{Y, \bar{y}} \to M_{Y, \bar{y}}$ of
the natural homomorphism
$p_y : M_{Y, \bar{y}} \to \overline{M}_{Y, \bar{y}}$, that is,
$p_y \cdot s_y = \operatorname{id}_{\overline{M}_{Y, \bar{y}}}$.
\end{enumerate}
Then, there is a splitting homomorphism
$s_x : \overline{M}_{X, \bar{x}} \to M_{X, \bar{x}}$ of
the natural homomorphism $p_x : M_{X, \bar{x}} \to \overline{M}_{X, \bar{x}}$
such that $p_x \cdot s_x = \operatorname{id}_{\overline{M}_{X, \bar{x}}}$
and the following diagram is commutative:
\[
 \begin{CD}
 \overline{M}_{Y, \bar{y}} @>{\bar{h}_x}>> \overline{M}_{X, \bar{x}} \\
 @V{s_x}VV @VV{s_y}V \\
 M_{Y, \bar{y}} @>{h_x}>> M_{X, \bar{x}}
 \end{CD}
\]
\end{Proposition}

\Proof
First of all, note that $\Coker(\OO^{\times}_{X, \bar{x}} \to
\phi^*(M_Y)_{\bar{x}}) = \overline{M}_{Y, \bar{y}}$. Moreover,
\[
 s'_y : \overline{M}_{Y, \bar{y}} \overset{s_y}{\longrightarrow} 
 M_{Y, \bar{y}} \to \phi^*(M_Y)_{\bar{x}}
\]
gives rise to a splitting homomorphism of
$\phi^*(M_Y)_{\bar{x}} \to \overline{M}_{Y, \bar{y}}$.

Let us consider the following commutative diagram:
\[
 \begin{CD}
 0 @>>> \OO^{\times}_{X, \bar{x}} @>>> \phi^*(M_Y)^{gr}_{\bar{x}} @>>>
 \overline{M}^{gr}_{Y, \bar{y}} @>>> 0 \\
 @. @| @VVV @VVV @. \\
 0 @>>> \OO^{\times}_{X, \bar{x}} @>>> M^{gr}_{X, \bar{x}} @>>> 
 \overline{M}^{gr}_{X, \bar{x}} @>>> 0,
 \end{CD}
\]
which gives rise to
\[
 \begin{CD}
 \Hom(\overline{M}^{gr}_{X, \bar{x}}, M^{gr}_{X, \bar{x}}) @>>>
 \Hom(\overline{M}^{gr}_{X, \bar{x}}, \overline{M}^{gr}_{X, \bar{x}}) 
 @>{\delta_1}>>
 \Ext^1(\overline{M}^{gr}_{X, \bar{x}}, \OO^{\times}_{X, \bar{x}}) \\
 @VVV @VV{\gamma_1}V  @VV{\lambda}V \\
 \Hom(\overline{M}^{gr}_{Y, \bar{y}}, M^{gr}_{X, \bar{x}}) @>>>
 \Hom(\overline{M}^{gr}_{Y, \bar{y}}, \overline{M}^{gr}_{X, \bar{x}}) 
 @>{\delta_2}>>
 \Ext^1(\overline{M}^{gr}_{Y, \bar{y}}, \OO^{\times}_{X, \bar{x}}) \\
 @AAA @AA{\gamma_2}A @| \\
 \Hom(\overline{M}^{gr}_{Y, \bar{y}}, \phi^*(M_Y)^{gr}_{\bar{x}}) @>>>
 \Hom(\overline{M}^{gr}_{Y, \bar{y}}, \overline{M}^{gr}_{Y, \bar{y}}) 
 @>{\delta_3}>>
 \Ext^1(\overline{M}^{gr}_{Y, \bar{y}}, \OO^{\times}_{X, \bar{x}}).
 \end{CD}
\]
By using the diagram
\[
 \begin{CD}
 \overline{M}^{gr}_{Y, \bar{y}} @>{\bar{h}^{gr}_x}>> 
 \overline{M}^{gr}_{X, \bar{x}} \\
 @| @| \\
 \overline{M}^{gr}_{Y, \bar{y}} @>{\bar{h}^{gr}_x}>> 
 \overline{M}^{gr}_{X, \bar{x}},
 \end{CD}
\]
we can see that 
$\gamma_1(\operatorname{id}_{\overline{M}^{gr}_{X, \bar{x}}}) =
\bar{h}^{gr}_x$ and
$\gamma_2(\operatorname{id}_{\overline{M}^{gr}_{Y, \bar{y}}}) =
\bar{h}^{gr}_x$.
Note that the exact sequence
\[
 0 \to \OO^{\times}_{X, \bar{x}} \to \phi^*(M_Y)^{gr}_{\bar{x}} \to
 \overline{M}^{gr}_{Y, \bar{y}} \to 0
\]
splits by ${s'}^{gr}_y$.
Thus,
\[
 \lambda(\delta_1(\operatorname{id}_{\overline{M}^{gr}_{X, \bar{x}}}))=
 \delta_2(\gamma_1(\operatorname{id}_{\overline{M}^{gr}_{X, \bar{x}}}))=
 \delta_2(\gamma_2(\operatorname{id}_{\overline{M}^{gr}_{Y, \bar{y}}}))=
 \delta_3(\operatorname{id}_{\overline{M}^{gr}_{Y, \bar{y}}}) = 0. 
\]
On the other hand, by our assumption, we can see
that 
\[
 \Ext^1(\overline{M}_{X, \bar{x}}/\overline{M}_{Y, \bar{y}},\ 
\OO_{X, \bar{x}}) = 0.
\]
Thus, we obtain that
$\lambda$ is injective. Therefore, 
$\delta_1(\operatorname{id}_{\overline{M}^{gr}_{X, \bar{x}}}) = 0$.
Hence, we have a splitting homomorphism
$s : \overline{M}^{gr}_{X, \bar{x}} \to M^{gr}_{X, \bar{x}}$ of
$M^{gr}_{X, \bar{x}} \to \overline{M}_{X, \bar{x}}$.

Here we claim that $s(\overline{M}_{X, \bar{x}}) \subseteq 
M_{X,\bar{x}}$.
Indeed, let us choose $a \in \overline{M}_{X, \bar{x}}$.
Then, there is $b \in M_{X, \bar{x}}$ with $p_x(b) = a$.
Since $p_x(s(a)) = a$, there is $c \in \OO^{\times}_{X, \bar{x}}$
such that $s(a) = b + c$ in $M^{gr}_{X, \bar{x}}$. 
Here $b, c \in M_{X,\bar{x}}$, which implies $s(a) \in M_{X, \bar{x}}$.

Therefore, we get a diagram
\[
 \begin{CD}
 \overline{M}_{Y, \bar{y}} @>{\bar{h}_x}>> \overline{M}_{X, \bar{x}} \\
 @V{s_y}VV @VV{s}V \\
 M_{Y,\bar{y}} @>{h_x}>> M_{X, \bar{x}}.
 \end{CD}
\]
Our problem is that the above diagram is not necessarily commutative.
By our assumption, for all $a \in \overline{M}_{Y, \bar{y}}$,
there is a unique $u \in \OO^{\times}_{X, \bar{x}}$
such that $s(\bar{h}_x(a)) + u = h_x(s_y(a))$.
We denote this $u$ by $\mu(a)$. Thus, we have a homomorphism
$\mu^{gr} : \overline{M}^{gr}_{Y, \bar{y}} \to 
\OO^{\times}_{X,\bar{x}}$.
Here we consider an exact sequence
\[
 0 \to \overline{M}^{gr}_{Y, \bar{y}} \to
       \overline{M}^{gr}_{X, \bar{x}} \to
       \overline{M}^{gr}_{X, \bar{x}}/\overline{M}^{gr}_{Y, \bar{y}}
    \to 0,
\]
which gives rise to
\[
 \Hom(\overline{M}^{gr}_{X, \bar{x}}, \OO^{\times}_{X, \bar{x}}) \to
 \Hom(\overline{M}^{gr}_{Y, \bar{y}}, \OO^{\times}_{X, \bar{x}}) \to
 \Ext^1(\overline{M}^{gr}_{X, \bar{x}}/\overline{M}^{gr}_{Y, \bar{y}},
        \OO^{\times}_{X, \bar{x}}) = \{0\}.
\]
Thus, there is $\nu \in
\Hom(\overline{M}^{gr}_{X, \bar{x}}, \OO^{\times}_{X, \bar{x}})$
with $\nu \cdot \bar{h}^{gr}_x = \mu^{gr}$.
Here we set $s_x = s + \nu$. Then,
\[
 s_x(\bar{h}_x(a)) = s(\bar{h}_x(a)) + \nu(\bar{h}_x(a)) =
 s(\bar{h}_x(a)) + \mu(a) = h_x(s_y(a)).
\]
Thus, we get our desired $s_x$.
\QED

\begin{Proposition}
\label{prop:log:smooth:good:chart}
Let $(\phi, h) : (X, M_X) \to (Y, M_Y)$ be
a smooth morphism of log schemes with fine log structures.
Let us fix $x \in X$ and $y = \phi(x)$.
We assume that there are \rom{(a)} etale neighborhoods
$U$ and $V$ of $x$ and $y$ respectively, 
\rom{(b)} charts $\pi_P : P \to \rest{M_X}{U}$ and 
$\pi_Q : Q \to \rest{M_Y}{V}$, and
\rom{(c)} a homomorphism $f : Q \to P$
with the following properties:
\begin{enumerate}
\renewcommand{\labelenumi}{(\arabic{enumi})}
\item
$\phi(U) \subseteq V$.

\item 
The induced homomorphism $P \to \overline{M}_{X, \bar{x}}$ and
$Q \to \overline{M}_{Y, \bar{y}}$ are bijective.


\item
The following diagram is commutative:
\[
 \begin{CD}
  Q @>{f}>> P \\
  @V{\pi_Q}VV @VV{\pi_P}V \\
  \rest{M_Y}{V} @>{h}>> \rest{M_X}{U}.
 \end{CD}
\]
\end{enumerate}
Then, the canonical morphism
$g : X \to Y \times_{\Spec(\ZZ[Q])} \Spec(\ZZ[P])$
is smooth around $x$ in the classical sense.
\end{Proposition}

\Proof
We consider the natural homomorphism
\[
 \alpha : \Coker(Q^{gr} \to P^{gr}) \otimes_{\ZZ} \OO_{X, \bar{x}}
          \to \Omega^1_{X/Y, \bar{x}}(\log(M_X/M_Y)).
\]
Let us begin with the following claim:

\begin{Claim}
$\alpha$ is injective and gives rise to a direct summand
of
\[
 \Omega^1_{X/Y, \bar{x}}(\log(M_X/M_Y)).
\]
\end{Claim}

In the same way as in \cite[(3.13)]{KatoLog},
we can construct a chart $\pi_{P'} : P' \to M_{X, \bar{x}}$
and an injective  homomorphism $f' : Q \to P'$
with the following properties:
\begin{enumerate}
\renewcommand{\labelenumi}{(\roman{enumi})}
\item
The torsion part of  $\Coker(Q^{gr} \to {P'}^{gr})$
is a finite group of order invertible in $\OO_{X, \bar{x}}$.

\item
The following diagram is commutative:
\[
 \begin{CD}
 Q @>{f'}>> P' \\
 @V{\pi_Q}VV @VV{\pi_{P'}}V \\
 M_{Y, \bar{y}} @>>> M_{X, \bar{x}}.
 \end{CD}
\]

\item
The natural homomorphism
\[
 \alpha' : \Coker(Q^{gr} \to {P'}^{gr})\otimes_{\ZZ} \OO_{X, \bar{x}} \to
 \Omega^1_{X/Y, \bar{x}}(\log(M_X/M_Y))
\]
is an isomorphism.
Moreover, there are $t_1, \ldots, t_r \in {P'}$ such that
a subgroup generated by $t_1, \ldots, t_r$ in $\Coker(Q^{gr} \to {P'}^{gr})$
is a free group of rank $r$ and its index
in $\Coker(Q^{gr} \to {P'}^{gr})$ is invertible in $\OO_{X, \bar{x}}$. 
In particular,\
\[
 d\log(\pi_{P'}(t_1)), \ldots, d\log(\pi_{P'}(t_r))
\] 
form a free basis of $\Omega^1_{X/Y, \bar{x}}(\log(M_X/M_Y))$.
\end{enumerate}
Considering the commutative diagram
\[
 \begin{CD}
 Q @>{\sim}>> \overline{M}_{Y, \bar{y}} @<{\sim}<< Q \\
 @V{f'}VV @V{\bar{h}_x}VV @VV{f}V \\
 P' @>>> \overline{M}_{X, \bar{x}} @<{\sim}<< P,
 \end{CD}
\]
we have a surjective homomorphism
$\lambda : P' \to P$ with $\lambda \cdot f' = f$.
Thus, we obtain the natural surjective homomorphism
\[
 \beta : \Coker(Q^{gr} \to {P'}^{gr}) \otimes_{\ZZ} \OO_{X, \bar{x}}\to
         \Coker(Q^{gr}\to P^{gr}) \otimes_{\ZZ}  \OO_{X, \bar{x}}.
\]
Hence, we have the following commutative diagram:
\[
 \xymatrix{
   \Coker(Q^{gr} \to {P'}^{gr}) \otimes_{\ZZ} \OO_{X, \bar{x}} 
   \ar[r]^{\ \ \sim}_{\ \ \alpha'} \ar[d]^{\beta} &   
   \Omega^1_{X/Y, \bar{x}}(\log(M_X/M_Y)) \\
   \Coker(Q^{gr} \to P^{gr}) \otimes_{\ZZ} \OO_{X, \bar{x}} \ar[ur]_{\alpha} 
 } 
\]
In order to see the claim, it is sufficient to see that
$\gamma = \beta \cdot {\alpha'}^{-1} \cdot \alpha$ is an automorphism on
$\Coker(Q^{gr} \to P^{gr}) \otimes_{\ZZ} \OO_{X, \bar{x}}$
because $(\beta \cdot {\alpha'}^{-1}) \cdot (\alpha \cdot \gamma^{-1}) =
\operatorname{id}$.
Here we set $\pi_{P'}(t_i) = p_i u_i$ ($p_i \in P$, 
$u_i \in \OO^{\times}_{X, \bar{x}}$) for $i=1, \ldots, r$.
Let us consider 
the natural surjective homomorphism
\begin{multline*}
 \theta :
 \Omega^1_{X/Y, \bar{x}}(\log(M_X/M_Y)) \otimes_{\ZZ} \kappa(\bar{x})\to \\
 \Coker(\overline{M}^{gr}_{Y, \bar{y}} \to 
 \overline{M}^{gr}_{X,\bar{x}}) \otimes_{\ZZ} \kappa(\bar{x})
 \simeq \Coker(Q^{gr} \to P^{gr}) \otimes_{\ZZ} \kappa(\bar{x})
\end{multline*}
given by $d\log(a) \mapsto a \otimes 1$ as in \cite[(3.13)]{KatoLog}.
This is nothing more than 
$(\beta \cdot {\alpha'}^{-1}) \otimes \kappa(\bar{x})$.
Indeed, 
\[
\begin{cases}
 (\beta \cdot {\alpha'}^{-1})(d\log(\pi_{P'}(t_i)))
= \beta(t_i) = p_i \\
\theta(d\log(\pi_{P'}(t_i))) = t_i = p_i \mod \OO^{\times}_{X, \bar{x}}.
\end{cases}
\]
On the other hand, we have the natural map
\[
 \alpha \otimes \kappa(\bar{x}) :
 \Coker(Q^{gr} \to P^{gr}) \otimes_{\ZZ} \kappa(\bar{x}) \to
 \Omega^1_{X/Y, \bar{x}}(\log(M_X/M_Y)) \otimes_{\ZZ} \kappa(\bar{x})
\]
given by $a \otimes 1 \mapsto d\log(a)$,
which is a section of $\theta$.
Therefore, $\gamma \otimes \kappa(\bar{x}) = \operatorname{id}$.
Thus, by Nakayama's lemma, $\gamma$ is surjective, so that
$\gamma$ is an isomorphism by \cite[Theorem~2.4]{MatComm}.

\medskip
We set $X' = Y \times_{\Spec(\ZZ[Q])} \Spec(\ZZ[P])$.
Let $\psi : X' \to \Spec(\ZZ[P])$ be the canonical morphism and
$M_P$ the canonical log structure on $\Spec(\ZZ[P])$.
We set $M_{X'} = \psi^*(M_P)$. Let $o$ the origin of $\Spec(\ZZ[P])$
and $x' = (y, o)$. Then, $M_{X', \bar{x}'} = 
\OO^{\times}_{X', \bar{x}'} \times P$.
Here, $\Omega^1_{X'/Y, \bar{x}'}$ is generated by 
$\{ d(1 \otimes x) \}_{x \in \ZZ[P]_{\bar{o}}}$. 
Thus, there is a natural surjective homomorphism
\[
 \Coker(Q^{gr} \to P^{gr})   \otimes_{\ZZ} \OO_{X', \bar{x}'}
\to \Omega^1_{X'/Y, \bar{x}'}(\log(M_{X'}/M_Y)).
\]
Therefore, we have a surjective homomorphism
\[
\Coker(Q^{gr} \to P^{gr})  \otimes_{\ZZ}  \OO_{X, \bar{x}}
\to g^*(\Omega^1_{X'/Y, \bar{x}'}(\log(M_{X'}/M_Y))).
\]
Thus, by the claim,
\[
 g^*(\Omega^1_{X'/Y, \bar{x}'}(\log(M_{X'}/M_Y))) \to
 \Omega^1_{X/Y, \bar{x}}(\log(M_{X}/M_Y))
\]
is injective and $g^*(\Omega^1_{X'/Y, \bar{x}'}(\log(M_{X'}/M_Y)))$
is a direct summand of
\[
 \Omega^1_{X/Y, \bar{x}}(\log(M_{X}/M_Y)).
\]
Therefore, by \cite[Proposition~(3.12)]{KatoLog},
$g$ is a smooth log morphism. Moreover, note that
$g^*(M_{X'}) = M_X$. Thus, $g$ is smooth in the classical
sense.
\QED

\end{document}